\documentclass[nonblindrev]{informs3} 

\DoubleSpacedXI 

\usepackage{algorithm}
\usepackage{algpseudocode}
\usepackage{amsmath}
\usepackage{tikz}
\usepackage{mathdots}
\usepackage{yhmath}
\usepackage{cancel}
\usepackage{color}
\usepackage{siunitx}
\usepackage{array}
\usepackage{multirow}
\usepackage{amssymb}
\usepackage{tabularx}
\usepackage{extarrows}
\usepackage{booktabs}
\usetikzlibrary{fadings}
\usetikzlibrary{patterns}
\usetikzlibrary{shadows.blur}
\usetikzlibrary{shapes}
\usepackage[font=small]{caption}
\usepackage[font=small]{subcaption}
\usepackage{dsfont}
\usepackage{mathtools}
\usepackage[T1]{fontenc}

\usepackage{algorithm}
\usepackage{algpseudocode}

\usepackage{xr-hyper}
\usepackage{hyperref}
\usepackage{fix-cm}

\usepackage{pifont}
\newcommand{\cmark}{\ding{51}}%
\newcommand{\xmark}{\ding{55}}

\newcommand{\boldsamewidth}[1]{%
    \pdfliteral direct {2 Tr 0.3 w} 
     #1%
    \pdfliteral direct {0 Tr 0 w}%
}

\usepackage{endnotes}
\let\enotesize=\normalsize
\def\notesname{Endnotes}%
\def\makeenmark{$^{\theenmark}$}
\def\enoteformat{\rightskip0pt\leftskip0pt\parindent=1.75em
\leavevmode\llap{\theenmark.\enskip}}

\usepackage{natbib}
 \bibpunct[, ]{(}{)}{,}{a}{}{,}%
 \def\bibfont{\small}%
\TheoremsNumberedThrough     
\ECRepeatTheorems

\EquationsNumberedThrough    

\newcommand*{\QEDA}{\hfill\ensuremath{\blacksquare}}
\newcommand*{\HOUSE}{\hfill\ensuremath{\blacktriangle}}

\newcommand{\bea}{\begin{equation*}\begin{aligned}}
\newcommand{\eea}{\end{aligned}\end{equation*}}

\newcommand{\R}{\mathbb{R}}

\newcommand{\wh}{\widehat}
\newcommand{\mbb}{\mathbb}

\newcommand{\PP}{\mbb P}

\newcommand{\QQ}{\mbb Q}

\newcommand{\MI}{\mathcal{A}}
\newcommand{\mi}{\alpha}

\DeclareMathOperator{\poly}{poly}

\DeclareMathOperator{\VAR}{-VAR}

\newcommand{\opt}{^\star}
\newcommand{\eps}{\varepsilon}

\newcommand{\BARY}{\mathbb{P}^\star}
\newcommand{\pibar}{\Bar{\pi}}

\newcommand{\cvar}{\eta}

\usepackage{xcolor}
\usepackage{amsmath}

\begin{document}


\RUNAUTHOR{Rychener, Esteban-P\'erez, Morales and Kuhn}

\RUNTITLE{Wasserstein DRO with Heterogeneous Data Sources}

\TITLE{Wasserstein Distributionally Robust Optimization with Heterogeneous Data Sources}

\ARTICLEAUTHORS{%
\AUTHOR{Yves Rychener}
\AFF{Risk Analytics and Optimization Laboratory,  \'Ecole polytechnique f\'ed\'erale de Lausanne, Lausanne CH-1015,  Switzerland \EMAIL{yves.rychener@alumni.epfl.ch}} 
 \AUTHOR{Adri\'an Esteban-P\'erez}
 \AFF{Department of Computer and Systems Sciences \& Digital Futures,
Stockholm University,
P.O.Box 1073, SE-164 25 Kista, Sweden \EMAIL{adrian.esteban@dsv.su.se}}
\AUTHOR{Juan M. Morales}
\AFF{Department of Mathematical Analysis, Statistics and Operations Research and Applied
Mathematics (Area of Statistics and Operations Research), University of Málaga, Málaga 29071, Spain \EMAIL{juan.morales@uma.es}}
\AUTHOR{Daniel Kuhn}
\AFF{Risk Analytics and Optimization Laboratory,  \'Ecole polytechnique f\'ed\'erale de Lausanne, Lausanne CH-1015, Switzerland \EMAIL{daniel.kuhn@epfl.ch}}
} 

    \ABSTRACT{We study decision problems under uncertainty, where the decision-maker has access to $K$ data sources that carry {\em biased} information about the underlying risk factors. The biases are measured by the mismatch between the risk factor distribution and the $K$ data-generating distributions with respect to an optimal transport (OT) distance. In this situation the decision-maker can exploit the information contained in the biased samples by solving a distributionally robust optimization (DRO) problem, where the ambiguity set is defined as the intersection of $K$ OT neighborhoods, each of which is centered at the empirical distribution on the samples generated by a biased data source. We show that if the decision-maker has a prior belief about the biases, then the out-of-sample performance of the DRO solution can improve with~$K$---irrespective of the magnitude of the biases. We also show that, under standard convexity assumptions, the proposed DRO problem is computationally tractable if either~$K$ or the dimension of the risk factors is kept constant.}

\KEYWORDS{data-driven decision-making; distributionally robust optimization; optimal transport} 

\maketitle

\section{Introduction}\label{sec:intro}
Stochastic optimization is the canonical framework for modeling decision problems under uncertainty \citep{shapiro2021lectures}. A basic single-stage stochastic program seeks a decision~$\theta\in\Theta\subseteq \R^n$ that minimizes the expected value $\mathbb{E}_\PP[\ell(\theta, \xi)]$ of an uncertainty-affected loss function~$\ell(\theta, \xi)$ with respect to the distribution~$\PP$ of the random vector~$\xi\in\Xi\subseteq \R^d$. However, in most real decision problems, the distribution~$\PP$ is unobservable, implying that an essential input of the stochastic optimization model is unknown. When the decision-maker has access to training samples~$\{\widehat \xi_i\}_{i=1}^N$ from~$\PP$, then the stochastic program can be addressed with the sample average approximation (SAA) \citep[\S~5]{shapiro2021lectures}, which minimizes the expected loss~$\mathbb{E}_{\widehat \PP} [\ell(\theta, \xi)]$ under the empirical distribution $\widehat{\mathbb{P}}=\frac{1}{N}\sum_{i=1}^{N} \delta_{\widehat{\xi}_{i}}$, or with methods from distributionally robust optimization (DRO) \citep{kuhn2024distributionally}, which minimize the worst-case expected loss $\sup_{ \QQ\in\mathcal P} \mathbb{E}_{\QQ} [\ell(\theta, \xi)]$ with respect to all distributions in an ambiguity set~$\mathcal P$. Intuitively, $\mathcal P$ should contain all distributions that are sufficiently likely to have generated the training samples. A popular choice is to define~$\mathcal P$ as a ball around~$\widehat{\mathbb{P}}$ with respect to a Wasserstein distance \citep{mohajerin2018data, zhao2018data, blanchet2019quantifying, kuhn2019wasserstein, gao2023wasserstein}. If the radius of this ball scales with~$N^{-\frac{1}{2}}$, then the optimal value of the DRO problem constitutes, in a precise statistical sense, a safe estimator for the optimal value of the underlying true stochastic program \citep{ref:blanchet2021statistical, gao2022finite}. The optimal solution of the SAA problem typically performs well {\em in sample} (that is, under the empirical distribution~$\widehat\PP$) but poorly out of sample (that is, under the unknown true distribution~$\PP$). This phenomenon is often referred to as the optimizer's curse \citep{smith2006optimizer}. Empirical evidence suggests that DRO mitigates the optimizer's curse.

Unfortunately, all existing approaches to data-driven optimization become ineffective at small sample sizes~$N$. This is troubling because training data is scarce or even absent in many relevant decision problems under uncertainty. For example, a manufacturer introducing a new product has no historical sales data to predict the product's demand distribution, which would help to plan production. Similarly, an investor trading in new securities has no historical market data to predict the securities' return distribution, which would help to design a portfolio strategy. In addition, a retailer expanding into a new market has no historical data on customer behavior to inform assortment planning, a logistics manager has no historical data on resource availabilities after a supply chain disruption to inform inventory decisions, and government authorities have no historical data on how people react to a new pandemic, which would help to determine effective health interventions. Broadly speaking, whenever an organization faces structural disruptions or implements strategic changes, it enters into a new regime characterized by a lack of relevant data.

In the remainder of the paper we refer to the distribution~$\PP$ of the uncertain problem parameters~$\xi$ as the {\em target distribution}, and we assume that there is insufficient (or no) data from~$\PP$. In this situation, one could try to leverage data from one or several {\em source distributions}~$\PP_k$, $k\in[K]$, that are not too dissimilar from~$\PP$. For example, a retailer expanding into a new market with an unknown demand distribution may have access to demand data from other markets with a comparable customer structure. Similarly, an investor trading in new securities with an unknown return distribution may have access to return data from other securities issued by similar companies. 

Suppose from now on that there are~$N_k\geq 1$ training samples from~$\PP_k$, which can be used to construct an empirical distribution~$\widehat\PP_k$ that approximates~$\PP_k$, $k\in[K]$. In addition, assume temporarily that~$\PP_1$ coincides with the target distribution~$\PP$.  This implies that~$\widehat\PP_1$ differs from~$\PP$ only due to statistical errors. In contrast, for~$k\geq 2$, the difference between~$\hat \PP_k$ and~$\PP$ can be decomposed into a statistical error from~$\widehat\PP_k$ to~$\PP_k$ and a distribution shift from~$\PP_k$ to~$\PP$. We now discuss different approaches for estimating the optimal value of the true stochastic program. A na\"ive approach would be to ignore all biased training samples and to solve the SAA problem under~$\widehat \PP_1$ only. As we assumed the target data to be scarce, however, this results in a fragile estimator with an unacceptably high variance. Alternatively, one could rely on data pooling, which amounts to solving the SAA problem under the mixture distribution $\sum_{k=1}^K N_k/(\sum_{k'=1}^K N_{k'}) \widehat \PP_k$.  We thus adopt the standard convention that, when datasets are merged, each observation is assigned equal weight. If there is abundant source data,  the resulting estimator has a small sampling variability. As the pooled dataset is dominated by biased samples, however, this estimator typically suffers from an unacceptably large bias. In particular, it is likely to overestimate the variance of~$\PP$ (see Example~\ref{ex:normal-barycenters} below). A better---geometry-aware---method for aggregating the empirical source distributions is to calculate their Wasserstein barycenter \citep{agueh2011barycenters}. Recall that the barycenter~$\overline \xi =\sum_{k=1}^K \lambda_k\xi_k$ of~$K$ points~$\xi_k\in\mathbb R^d$ with weights~$\lambda_k\geq 0$, $k\in[K]$, is the unique solution of $\min_{\xi\in\mathbb R^d} \sum_{k=1}^K\lambda_k \|\xi-\xi_k\|_2^2$. Analogously, the $p$-Wasserstein barycenter of $K$ distributions~$\widehat\PP_k$ with weights~$\lambda_k\geq 0$, $k\in[K]$, is defined as a solution~$\overline \PP$ of $\min_{\QQ\in\mathcal P(\Xi)} \sum_{k=1}^K\lambda_k W_p({\QQ},\widehat\PP_k)^p$, where~$W_p$ denotes the $p$-Wasserstein distance for some~$p\geq 1$. Wasserstein barycenters are widely used in machine learning \citep{claici2018stochastic,dognin2019wasserstein,montesuma2021wasserstein,srivastava2018scalable,zhuang2022kmeans} and computer vision \citep{barre2020averaging,solomon2015convolutional}. The statistical properties of empirical Wasserstein barycenters are studied in \citep{legouicloubes2017existence, legouic2022fastalexbary,panaretos2020invitation}, and efficient methods for their computation are proposed in \citep{altschuler2021barycenterComplexity, cuturi14wasserstein-barycenters}. Any $p$-Wasserstein barycenter~$\overline\PP$ of the empirical source distributions is {\em discrete} and can thus be identified with an aggregated dataset. Solving the SAA problem under~$\overline\PP$ yields a more promising estimator than data pooling. The following example suggests indeed that, unlike na\"ive mixtures, Wasserstein barycenters preserve stylized features of the source distributions.

\begin{example}[Mixtures versus Wasserstein Barycenters]
\label{ex:normal-barycenters}
Assume that all source distributions are univariate Gaussians with different means and identical variances, that is, $\PP_k=\mathcal N(\mu_k,\sigma^2)$ for all~$k\in[K]$. If $\sum_{k=1}^K\lambda_k=1$ and $\lambda_k\geq0$ for all~$k\in[K]$, then the mixture distribution $\sum_{k=1}^K\lambda_k\PP_k$ fails to be Gaussian, and it has mean~$\sum_{k=1}^K\lambda_k\mu_k$ and variance $\sigma^2+\sum_{k=1}^K\lambda_k\mu^2_k-(\sum_{k=1}^K\lambda_k\mu_k)^2> \sigma^2$. Using \cite[Theorem~2.1]{gelbrich1990formula}, one can further show that the 2-Wasserstein barycenter of the source distributions with weights~$\lambda_k$, $k\in[K]$, is the Gaussian distribution with mean~$\sum_{k=1}^K\lambda_k\mu_k$ and variance $\sigma^2$. Thus, the 2-Wasserstein barycenter captures shared features of the source distributions such as their normality and identical variance~$\sigma^2$. In contrast, the mixture distribution is structurally different. It is multimodal, and its variance exceeds~$\sigma^2$. In applications where the source distributions are believed to be structurally similar to the target distribution~$\PP$, the 2-Wasserstein barycenter may thus be the more plausible model for~$\PP$ than the mixture of the source distributions.~\HOUSE
\end{example}

Despite the encouraging insights from Example~\ref{ex:normal-barycenters}, the 2-Wasserstein barycenter of two empirical source distributions is highly sensitive to data perturbation \citep[Example~2]{zhuang2022kmeans}. In addition, its variance provides only a {\em biased} estimator for the variance of the 2-Wasserstein barycenter of the true source distributions (see Proposition~\ref{prop:moment-reducing-barycenter} below). Even worse, the 1-Wasserstein barycenter either coincides with one of the source distributions or is highly degenerate (see Corollary~\ref{cor:barycenterK2} below).

So far we have discussed several SAA models, which differ only in how the source data enters the empirical distribution. Recall that SAA methods are generally susceptible to the optimizer's curse and may thus lead to disappointment in out-of-sample tests. Following \citet{mohajerin2018data}, it is therefore natural to robustify each of these SAA models against all distributions in a Wasserstein ball centered at the corresponding empirical distribution. If this empirical distribution converges (in an appropriate sense) to the target distribution, then the resulting DRO models offer the usual statistical guarantees; see {\em e.g.}, \citep[\S~3]{kuhn2019wasserstein}. In particular, if the target distribution~$\PP$ is known to coincide with the $2$-Wasserstein barycenter of the source distributions~$\PP_k$ with weights~$\lambda_k$, $k\in[K]$, and if the empirical distribution is set to the $2$-Wasserstein barycenter of the {\em empirical} source distributions~$\widehat\PP_k$ with the {\em same} weights~$\lambda_k$, $k\in[K]$, then the optimal value of the DRO problem provides an asymptotically consistent upper confidence bound on the optimal value of the true stochastic program \citep[Theorems~4.7 and~4.8]{lau2022wassersteinbarydro}. However, this guarantee is only available if the weights~$\lambda_k$, $k\in[K]$, are known.  In addition, as we will see in the numerical experiments of Section~\ref{sec:experiments}, the solution of the DRO problem inherits all shortcomings of the Wasserstein barycenter at the center of its ambiguity set---such as outlier sensitivity, degeneracy and biasedness. More worryingly, we demonstrate in this paper that the Wasserstein barycenter-based DRO model suffers from instability, that is, its optimal value as well as its optimal solution can change discontinuously under Wasserstein perturbations of the empirical source distributions.

In view of the challenges outlined above, we propose a new approach to data-driven decision-making with~$K>1$ data sources. Recall that any useful source distribution~$\PP_k$, $k\in[K]$, must be close to the target distribution~$\PP$.  It is therefore reasonable to assume that the decision-maker has a prior belief regarding the discrepancy between~$\PP$ and~$\PP_k$, as measured by a Wasserstein distance. As the empirical distribution~$\widehat\PP_k$ approximates~$\PP_k$, any such belief naturally translates into an upper bound on the Wasserstein distance between~$\PP$ and~$\widehat\PP_k$. Hence, the target distribution should reside in the intersection of~$K$ Wasserstein balls, each of which is centered at one of the empirical source distributions. This prompts us to introduce the class of {\em multi-source DRO} models whose ambiguity sets are obtained by intersecting {\em multiple} Wasserstein balls associated with {\em different data sources}. 

The main contributions of this paper can be summarized as follows.

\begin{itemize}
\item We introduce multi-source DRO and show that it is dually related to the problem of computing the Wasserstein barycenter of the $K$ empirical source distributions. We also prove that multi-source DRO is stable under mild conditions, that is, its optimal value as well as its optimal solution change continuously under Wasserstein perturbations of the empirical source distributions.

\item We show that if the loss function is piecewise concave in~$\xi\in\mathbb R^d$, then the worst-case expected loss over the intersection of the Wasserstein balls associated with the $K$ data sources matches the optimal value of a finite convex program. However, this convex program involves exponentially many variables and constraints and is NP-hard. Leveraging recent results on Wasserstein barycenters by \citet{altschuler2021barycenterComplexity}, we then prove that it becomes tractable if~$d$ or~$K$ is constant.

\item We develop frequentist as well as Bayesian performance guarantees for multi-source DRO. Specifically, we provide guidance for calibrating the Wasserstein radii such that the multi-source DRO problem provides an upper confidence bound on the stochastic program under the target distribution. We find that a {\em frequentist} decision-maker without any information about the distances between the target and the source distributions does {\em not} benefit from source data. In contrast, we show that a Bayesian decision-maker with prior beliefs about the distances between the target and the source distributions can obtain stronger performance guarantees by using source data.

\item Numerical experiments focusing on portfolio selection and assortment planning show that multi-source DRO can outperform various single-source DRO schemes on synthetic as well as real data.
\end{itemize}

All computational results of this paper are established within a general framework based on optimal transport (OT) discrepancies, encompassing Wasserstein distances as a special case. The statistical and numerical analyses as well as the stability results, however, focus on Wasserstein distances. Our tractability results imply that solving a multi-source DRO problem is not fundamentally harder than computing the Wasserstein barycenter of the source distributions. Thus, our multi-source DRO problems belong to the same complexity class as the single-source DRO problems by \citet{lau2022wassersteinbarydro}, which use the Wasserstein barycenter of the source distributions as the nominal distribution. However, multi-source DRO is less sensitive to perturbations of the data.

Intersections of $K=2$ ambiguity sets constructed from source and target data were first studied in the context of distributionally robust linear regression~\citep{taskesen2021domainadaptation}. This work critically relies on the quadratic nature of the least-squares loss, which implies that only the first and second moments of the covariates and the responses are needed to evaluate the expected prediction loss. More recently and independently of our work, intersections of $K=2$ Wasserstein ambiguity sets were studied in the context of logistic regression \citep{awasthi2022dro-join}. The dual DRO problem derived in this work is a special case of our Theorem~\ref{th_strong_duality}. As the logistic loss fails to be (piecewise) concave, this dual problem does not admit an exact reformulation as a finite convex program. \cite{selvi2024intersection} study a generalization of the model by \citet{awasthi2022dro-join} that not only accounts for $K=2$ different data sources but also for adversarial attacks. They leverage tools from adjustable robust optimization to solve their model approximately. Concurrent to our work, \citet{wang2024contextual} introduce a contextual DRO model that provides protection against covariate shifts. The underlying ambiguity set is given by an intersection of $K=2$ 1-Wasserstein balls, whose centers correspond to different parametric and nonparametric estimators of the relevant conditional distribution. These estimators are all constructed from one single data source. \citet{wang2024contextual} also provide a duality result for a DRO problem with $K\geq 2$ Wasserstein balls, which is similar to our Theorem~\ref{th_strong_duality}. \citet{tanoumand2023data} explore the intersection of a Wasserstein ball with an ambiguity set derived from a goodness-of-fit test based on the linear-convex ordering of random vectors. Their rationale for intersecting ambiguity sets is to reduce the conservativeness of the DRO solutions. However, they construct both ambiguity sets from a single dataset, and therefore they do not cater for multi-source DRO. Intersections of ambiguity sets also play a role in the analysis of market equilibria when agents optimize coherent risk measures \citep{ralph2011pricing}. Indeed, it can be shown that an equilibrium exists if the ambiguity sets underlying the coherent risk measures of all agents have a nonempty intersection.

Decision-making under uncertainty with multiple data sources is reminiscent of domain adaptation in machine learning \citep{zhao2020multi, farahani2021brief, montesuma2021wasserstein} and also resembles federated learning with heterogeneous clients \citep{mohri2019agnostic, rofederated2021}. Domain adaptation leverages data from one or several source distributions to improve the performance of a machine learning model on the target distribution \citep{singhal2023domain}. This can be achieved, for example, by using tools from optimal transport. Specifically, an optimal transport map between the source and the target distribution can be used to convert a predictor trained on the source domain to one on the target domain \citep{courty2014domain,courty2016optimal,redko2017theoretical}. However, these methods are tailored to prediction problems and cannot be readily extended to decision problems.

Heterogeneous data sources are also common in federated learning, where different devices share data to train a machine learning model in a distributed manner---often under privacy constraints (see, {\em e.g.}, \citep{zhang2021survey} for a survey). Distributionally robust federated learning models trained on a mixture of the client distributions with uncertain mixture weights are studied in \citep{deng2020distributionally,sagawa2019distributionally,xiong2023distributionally,wang2023dromultisource}. A distributionally robust model with a Wasserstein ball around the empirical mixture distribution is proposed in \citep{nguyen2022generalization}.

Data-driven optimization with heterogeneous data sources has not yet received much attention in operations research and management science. It is known that if many unrelated data-driven decision problems must be solved simultaneously, each having only access to a small amount of data, then shrinking the empirical distribution of each individual problem towards an anchor distribution constructed from the pooled data can improve upon na\"ive SAA approaches \citep{gupta2022data}. However, as all empirical distributions are shrunk towards the same anchor, this approach is difficult to justify in the presence of distribution shifts. Also, we focus on a single decision problem.

The rest of the paper is structured as follows. Section~\ref{sec:barycenters} introduces optimal transport barycenters as natural generalizations of Wasserstein barycenters. Section~\ref{sec:multi-source-dro} then formally introduces the multi-source DRO problem and shows that it is equivalent to a finite convex program of exponential size if the loss function is piecewise concave. Section~\ref{sec:tractability} shows that this convex program is tractable if either~$d$ or~$K$ is constant. Section~\ref{sec:stat_guarantees} establishes statistical performance guarantees for multi-source DRO, and Section~\ref{sec:experiments} reports on numerical experiments in the context of portfolio optimization.  Appendix~\ref{app:stability} develops a stability theory for multi-source DRO, while all proofs are collected in Appendix~\ref{sec:proofs}. Finally, Appendix~\ref{sec:scenario-relationships} provides supplementary details on the statistical guarantees in Section~\ref{sec:stat_guarantees}, and Appendix~\ref{sec:additional-experiments} reports on additional numerical experiments.

\paragraph{Notation.}
For any Borel set $\Xi\subseteq\mathbb{R}^d$, we use $\mathcal{M}_+(\Xi)$, $\mathcal{P}(\Xi)$ and $\mathcal{P}_p(\Xi)$ to denote the convex cone of all finite Borel measures on $\Xi$, the convex set of all probability distributions in $\mathcal{M}_+(\Xi)$ and the convex set of all probability distributions in $\mathcal{P}(\Xi)$ with finite $p$-th moment, respectively. For any $\mathbb{P}\in\mathcal{P}(\Xi)$ and Borel measurable transformation $\Phi:\Xi\rightarrow\Xi'$ between Borel sets $\Xi\subseteq\mathbb{R}^d$ and $\Xi'\subseteq\mathbb{R}^{d'}$, we denote by $\Phi_{\#}\mathbb{P}$ the pushforward distribution of $\mathbb{P}$ under $\Phi$. Thus, if $\xi\in\Xi$ follows~$\mathbb{P}$, then $\Phi(\xi)\in\Xi'$ follows $\Phi_{\#}\mathbb{P}$. 
The conjugate $f^*$ of a convex function $f:\mathbb{R}^d\rightarrow [-\infty,\infty]$ is defined through $f^*(u)=\sup_{\xi\in\mathbb{R}^d}\langle u, \xi\rangle -f(\xi)$. The indicator function $\iota_\Xi$ of a closed convex set $\Xi\subseteq\mathbb{R}^d$ is defined through $\iota_\Xi(\xi)=0$ if $\xi\in\Xi$, $\iota_\Xi(\xi)=\infty$ otherwise. The conjugate of $\iota_\Xi$ is termed the \emph{support function} of $\Xi$ and is denoted by $\sigma_\Xi$. The \emph{perspective} of a closed convex function $f(\xi)$ is defined as $tf(\xi/t)$ for $t\geq 0$, where $0f(\xi/0)$ is interpreted as $\sigma_{{\rm dom}(f^*)}(\xi)$. For any $n\in\mathbb N$ we set $[n]=\{1,\ldots,n\}$.

\section{Multi-Margin OT and OT Barycenters}\label{sec:barycenters}
This section introduces OT barycenters as natural generalizations of Wasserstein barycenters and shows that they can be computed by solving multi-margin OT problems. Throughout the paper, we let $\Xi\subseteq\mathbb{R}^d$ be a  closed convex support set. The family of multi-margin transportation plans associated with $K$ probability distributions $\mathbb{P}_1,\ldots,\mathbb{P}_K\in\mathcal{P}(\Xi)$ is defined as
$$
\Pi(\mathbb{P}_1, \ldots, \mathbb{P}_K)=\left\{\pi\in\mathcal{P}(\Xi^{K})\ :\ \pi\ \text{has marginals}\ \mathbb{P}_1, \ldots, \mathbb{P}_K,\ \text{respectively}\right\}
$$
\citep{pass2015multimargin}. If $K=2$, then any $\pi\in\Pi(\mathbb{P}_1, \mathbb{P}_2)$ is simply referred to as a transportation plan or a coupling of $\mathbb{P}_1$ and $\mathbb{P}_2$.
Optimal transport seeks a transportation plan of minimum cost with respect to a cost function $c:\Xi\times\Xi\rightarrow \mathbb{R}_+$. We adopt the following standard assumption regarding~$c$.

\begin{assumption}[Transportation Cost Function]\label{ass:assumptions-c}
The function $c(\xi_1, \xi_2)$ is lower semi-con\-tinuous, obeys the identity of indiscernibles (i.e., $c(\xi_1,\xi_2)=0$ if and only if~$\xi_1=\xi_2$) and is lower bounded by $d^p(\xi_1,\xi_2)$ for some metric $d$ on $\Xi$ with compact sublevel sets and some exponent~$p\geq 1$.
\end{assumption}

The OT cost associated with two distributions $\mathbb{P}_1, \mathbb{P}_2\in\mathcal{P}(\Xi)$ is then defined as 
$$
C(\mathbb{P}_1, \mathbb{P}_2)=\min_{\pi\in\Pi(\mathbb{P}_1, \mathbb{P}_2)} \int_{\Xi^2} c(\xi_1, \xi_2) {\rm d}\pi(\xi_1, \xi_2).
$$
Intuitively, the OT problem on the right hand side of the above expression seeks the cheapest transportation plan $\pi$ for morphing $\mathbb{P}_1$ into $\mathbb{P}_2$, where the cost of moving unit mass from $\xi_1$ to $\xi_2$ amounts to $c(\xi_1,\xi_2)$. The OT problem is always solvable ({\em i.e.}, an optimal solution exists) under Assumption~\ref{ass:assumptions-c} \citep[Theorem~4.1]{villani2009OT}.
If $c(\xi_1,\xi_2)=d^p(\xi_1,\xi_2)$ for some metric $d$ on $\Xi$ and some exponent $p\geq 1$, then $C(\mathbb{P}_1, \mathbb{P}_2)^{1/p}$ is termed the $p$-Wasserstein distance and is denoted by~$W_p(\mathbb{P}_1, \mathbb{P}_2)$.

We are now ready to define the OT barycenter of multiple distributions. It generalizes the notion of the Wasserstein barycenter introduced by~\cite{agueh2011barycenters}.
\begin{definition}[OT Barycenter]\label{def:wasserstein-barycenter}
An OT barycenter $\BARY$ of $K$ distributions $\mathbb{P}_1,\ldots,\mathbb{P}_K\in\mathcal{P}(\Xi)$ with associated weights $\lambda_1,\ldots,\lambda_K\geq 0$ is any solution of the minimization problem
\begin{equation}\label{eq:barycenter}\tag{OT-BC}
\min_{\QQ\in\mathcal{P}(\Xi)} \sum_{k=1}^K\lambda_k C({\QQ}, \mathbb{P}_k).
\end{equation}
If $C({\QQ}, \mathbb{P}_k)=W_p^p({\QQ}, \mathbb{P}_k)$, then the OT barycenters are referred to as $p$-Wasserstein barycenters.
\end{definition}

\begin{assumption}[OT Cost]\label{ass:barycenter-existence}
    There exists $\xi_0\in\Xi$ with $\int_\Xi c(\xi_0, \xi_k)\,{\rm d}\mathbb{P}_k(\xi_k)<\infty$ for all $k\in[K]$. 
\end{assumption}

The following lemma shows that OT barycenters must exist under Assumptions~\ref{ass:assumptions-c} and~\ref{ass:barycenter-existence}.

\begin{lemma}[OT Barycenters]\label{lemma:barycenter:existence}
Assumptions~\ref{ass:assumptions-c} and~\ref{ass:barycenter-existence} ensure that problem~\eqref{eq:barycenter} is solvable.
\end{lemma}

The next theorem shows that OT barycenters can be found by solving multi-margin OT problems. It extends an existing result for 2-Wasserstein barycenters \citep{agueh2011barycenters, anderes2016discrete} to general OT barycenters. We include a short proof to keep the paper self-contained.

\begin{theorem}[OT Barycenters and Multi-Margin Transportation Plans]\label{thm:barycenter:pushforward-reformulation}
Suppose that Assumptions~\ref{ass:assumptions-c} and~\ref{ass:barycenter-existence} hold, and define $\phi:\Xi^K\rightarrow [0,\infty]$ and $\Phi:\Xi^K\rightarrow\Xi$ as the optimal value function and a measurable selector of the set-valued solution map of the minimization problem 
\begin{equation}\label{eq:barycenter:inner-min}
    \min_{\xi\in\Xi}\;\sum_{k=1}^K \lambda_k c(\xi, \xi_k)
\end{equation}
parametrized by $\xi_1,\ldots, \xi_K\in\Xi$. Then there is an OT barycenter of $\mathbb{P}_1,\ldots, \mathbb{P}_K\in\mathcal{P}(\Xi)$ representable as $\BARY=\Phi_{\#}\pi\opt$, where $\pi\opt$ is an optimal solution of the multi-margin OT problem
\begin{equation}\label{eq:barycenter:outer-min}
    \min_{\pi\in\Pi(\mathbb{P}_1,\ldots, \mathbb{P}_K)} \int_{\Xi^K} \phi(\xi_1,\ldots,\xi_K)\,{\rm d}\pi(\xi_1, \ldots,\xi_K).
\end{equation}
If the measurable selector $\Phi$ is unique, then every OT barycenter is representable in this form.
\end{theorem}

\cite{agueh2011barycenters} prove that the 2-Wasserstein barycenter  is unique if the distributions $\mathbb P_1,\ldots,\mathbb P_K$ assign probability zero to all sets of Hausdorff dimension less than~$d$. Conversely, Theorem~\ref{thm:barycenter:pushforward-reformulation} suggests that the OT barycenter is {\em not} unique if either the multi-margin transportation plan $\pi^\star$ or the minimizer map $\Phi$ fails to be unique. We illustrate each case with an example. These examples are based on discrete distributions, which are supported on sets of Hausdorff dimension~$0$.

\begin{example}[Non-Unique Multi-Margin Transportation Plan]\label{ex:nonunique-bary-pi} The OT bary\-center of the distributions $\mathbb{P}_1 = \frac{1}{2}(\delta_{(1,1)} + \delta_{(0,0)})$ and $\mathbb{P}_2 = \frac{1}{2}(\delta_{(0,1)} + \delta_{(1,0)})$ in $\mathcal P_2(\mathbb{R}^2)$ with $\lambda_1=\lambda_2=1$ and $c(\xi,\xi')=\|\xi-\xi'\|_2^2$ is not unique. Indeed, while the minimizer map $\Phi(\xi_1,\xi_2)=\frac{1}{2}(\xi_1+\xi_2)$ of problem~\eqref{eq:barycenter:inner-min} is unique, both $\pi^{\star 1} = \frac{1}{2}(\delta_{((1,1),(1,0))} + \delta_{((0,0),(0,1))})$ as well as $\pi^{\star 2} = \frac{1}{2}(\delta_{((1,1),(0,1))} + \delta_{((0,0),(1,0))})$ solve the multi-margin OT problem~\eqref{eq:barycenter:outer-min}. By Theorem~\ref{thm:barycenter:pushforward-reformulation}, both $\BARY{}^1 = \frac{1}{2}(\delta_{(1,0.5)} + \delta_{(0,0.5)})$ as well as $\BARY{}^2 = \frac{1}{2}(\delta_{(0.5,1)} + \delta_{(0.5,0)})$ thus constitute minimizers of the OT barycenter problem~\eqref{eq:barycenter}.
\HOUSE
\end{example}

\begin{example}[Non-Unique Minimizer Map]\label{ex:nonunique-bary-phi}
The OT barycenter of two distributions $\mathbb{P}_1\neq \mathbb{P}_2$ in $\mathcal P_2(\mathbb{R}^2)$ with $\lambda_1=\lambda_2=1$ and $c(\xi,\xi')=\|\xi-\xi'\|_2$ is not unique. Indeed, while the multi-margin transportation plan $\pi^\star$ that solves~\eqref{eq:barycenter:outer-min} may or may not be unique in this situation, both $\Phi^{\star 1}(\xi_1,\xi_2)=\xi_1$ and $\Phi^{\star 2}(\xi_1,\xi_2)=\xi_2$ constitute minimizer maps for problem~\eqref{eq:barycenter:inner-min}. As $\pi^\star\in\Pi(\mathbb{P}_1,\mathbb{P}_2)$, Theorem~\ref{thm:barycenter:pushforward-reformulation} thus implies that both $\mathbb{P}_1=\Phi^{\star 1}_{\#}\pi\opt$ and $\mathbb{P}_2=\Phi^{\star 2}_{\#}\pi\opt$ solve the OT barycenter problem~\eqref{eq:barycenter}.
\HOUSE
\end{example}

To close this section, we show that if distances between distributions are measured by a Wasser\-stein distance of order~1 or~2, then the resulting OT barycenters may have undesirable properties. We first show that the 1-Wasserstein barycenter of $K=2$ distributions coincides with the distribution that has the larger weight and is degenerate if both distributions have the same weight.

\begin{corollary}[1-Wasserstein Barycenter of Two Distributions]\label{cor:barycenterK2}
If $\mathbb{P}_1, \mathbb{P}_2\in\mathcal{P}(\Xi)$ satisfy $\int_{\Xi} \|\xi\|\rm d\mathbb{P}_k(\xi)<\infty$ for $k=1,2$ and if $d(\xi,\xi')=\|\xi-\xi'\|$, then a 1-Wasserstein barycenter of $\mathbb{P}_1$ and $\mathbb{P}_2$ with respective weights $\lambda_1, \lambda_2\geq0$ is given by $\mathbb{P}_1$ if $\lambda_1\geq\lambda_2$ and by $\mathbb{P}_2$ if $\lambda_1\leq\lambda_2$.
\end{corollary}

Any reasonable aggregation of two distributions~$\PP_1$ and~$\PP_2$ should differ from both; selecting either~$\PP_1$ or~$\PP_2$ as the aggregate would be unsatisfactory. Corollary~\ref{cor:barycenterK2} thus shows that some $1$-Wasserstein barycenters fail to yield a meaningful aggregation of two distributions. In contrast, $2$-Wasserstein barycenters are widely used for aggregating distributions in computer graphics \citep{solomon2015convolutional}, machine learning \citep{dognin2019wasserstein, montesuma2021wasserstein}, Bayesian inference \citep{srivastava2018scalable} and statistics \citep{claici2018stochastic, barre2020averaging} etc.
In all of these applications, however, the reference distributions $\mathbb{P}_1,\ldots,\mathbb{P}_K$ are typically unknown. Given $N_k$ independent samples $\{\widehat{\xi}_{k,j}\}_{j=1}^{N_k}$ from~$\mathbb{P}_k$, $k=1,\ldots, K$, one can approximate the unobservable distribution $\mathbb{P}_k$ with the corresponding empirical distribution $\widehat{\mathbb{P}}_k=\frac{1}{N_k}\sum_{j=1}^{N_k}\delta_{\widehat{\xi}_{k,j}}$. The unknown 2-Wasserstein barycenter $\BARY$ of $\mathbb{P}_1,\ldots,\mathbb{P}_K$ can thus be estimated by the empirical 2-Wasserstein barycenter $\widehat{\mathbb{P}}$ of $\widehat{\mathbb{P}}_1,\ldots,\widehat{\mathbb{P}}_K$, which is itself a discrete distribution. The next proposition shows that this estimator is biased. In this proposition, we use $\mathbb{P}\VAR(\xi)=\mathbb E_{\mathbb P}[\|\xi-\mathbb E_{\mathbb P}[\xi]\|_2^2]$ to denote the trace of the covariance matrix of~$\xi\in\mathbb R^d$ under the distribution~$\mathbb P\in\mathcal P_2(\Xi)$. For brevity, we will sometimes refer to $\mathbb{P}\VAR(\xi)$ as the variance of~$\xi$. It is easy to verify that $\mathbb{P}\VAR(\xi)=\min_{\xi_0\in\Xi} \mathbb E_{\mathbb P}[\|\xi-\xi_0\|_2^2]$, where the minimization problem over~$\xi_0$ is solved by $\mathbb E_\mathbb P[\xi]$. This formula shows that $\mathbb{P}\VAR(\xi)$ is concave in~$\mathbb P$.

\begin{proposition}[Empirical OT Barycenters are Biased]\label{prop:moment-reducing-barycenter}
    Set $c(\xi_1,\xi_2)=\|\xi_1-\xi_2\|_2^2$ and $\Xi=\mathbb R^d$, and select $\lambda_1,\ldots, \lambda_K>0$. Define $\mathbb P^\star$ as the OT barycenter of $K$ mutually different continuous distributions $\mathbb{P}_1,\ldots,\mathbb P_K\in\mathcal P_2(\Xi)$. Similarly, define $\widehat{\mathbb{P}}$ as an OT barycenter of $\widehat{\mathbb{P}}_1,\ldots,\widehat{\mathbb P}_K$, where $\widehat{\mathbb{P}}_k$ denotes the empirical distribution corresponding to $N_k$ random samples from~$\mathbb P_k$, $k\in[K]$. Then, we have $\mathbb{E}[\widehat{\mathbb{P}}\VAR(\xi)]< \mathbb{P}^\star\VAR(\xi)$, where the expectation is with respect to all samples.
\end{proposition}

As the transportation cost function $c(\xi_1,\xi_2)$ is quadratic and the distributions $\mathbb P_1,\ldots,\mathbb P_K$ are square-integrable, Assumptions~\ref{ass:assumptions-c} and~\ref{ass:barycenter-existence} both hold under the conditions of Proposition~\ref{prop:moment-reducing-barycenter}. Hence, $\mathbb P^\star$ exists thanks to Lemma~\ref{lemma:barycenter:existence}, and it is unique thanks to~\cite[Proposition~3.5]{agueh2011barycenters}. The empirical OT barycenter $\widehat{\mathbb{P}}$ trivially exists because $\widehat{\mathbb{P}}_1,\ldots,\widehat{\mathbb P}_K$ are discrete.

Having covered the basic theory of OT barycenters, we now turn to data-driven decision-making. We will see that the tools for constructing and analyzing OT barycenters developed in Section~\ref{sec:barycenters} will be instrumental for constructing data-driven DRO models with heterogeneous data sources.

\section{Multi-Source DRO}
\label{sec:multi-source-dro}
From now on we consider a stochastic optimization problem of the form
\begin{equation}
    \label{eq:sp}
    \min_{\theta\in\Theta} \int_{\Xi}\ell(\theta, \xi)\, {\rm d}\mathbb{P}(\xi),
\end{equation}
which minimizes the expected value of an uncertainty-affected loss function $\ell:\Xi\times\Theta\rightarrow \mathbb{R}$ over a feasible set $\Theta\subseteq \mathbb R^n$. We henceforth assume that $\ell(\theta, \xi)$ is lower semicontinuous in $\theta$ for every fixed $\xi\in\Xi$ and upper semicontinuous in~$\xi$ for every fixed $\theta\in\Theta$. We also assume that $\Theta$ is closed. In most applications of stochastic programming, the target distribution~$\PP$, which is needed to evaluate the expected loss of any fixed decision~$\theta$, is unknown. Given samples from~$\PP$, however, one can approximate~$\PP$ with an estimated distribution~$\widehat\PP$ constructed from the data. Solving~\eqref{eq:sp} under~$\widehat\PP$ instead of~$\PP$ may lead to overfitting effects if data is scarce. That is, the optimal solution computed under~$\widehat \PP$ may underperform under~$\PP$, a phenomenon known as the {\em optimizer's curse} \citep{smith2006optimizer}. As a possible remedy, it has been proposed to solve a DRO problem of the form
\begin{equation}
    \label{eq:dro}
    \inf_{\theta\in\Theta}~\sup_{{\QQ} \in \mathbb{B}_\varepsilon(\widehat{\mathbb{P}})}~\int_{\Xi}\ell(\theta, \xi)\, {\rm d}{\QQ}(\xi),
\end{equation}
where $\mathbb{B}_\varepsilon(\widehat{\mathbb{P}})=\{{\QQ}\in\mathcal{P}(\Xi):\ C({\QQ},\widehat{\mathbb{P}})\leq\varepsilon\}$ denotes the set of all distributions whose OT distance to~$\widehat{\mathbb{P}}$ is at most~$\varepsilon\geq 0$ \citep{pflug2012}. The DRO problem~\eqref{eq:dro} pretends that an evil adversary can morph the estimated distribution~$\widehat\PP$ into some undesirable shape at a finite transportation budget~$\varepsilon$. It turns out that preparing for destructive actions of a fictitious adversary can mitigate the optimizer's curse. The statistical and computational properties of OT-based DRO are studied in \citep{mohajerin2018data, zhao2018data, blanchet2019quantifying, gao2023wasserstein}.

In this paper we address more challenging decision-making situations in which data from the target distribution~$\PP$ is extremely scarce or even unavailable. Instead, we assume to have access to data from~$K>1$ source distributions~$\PP_k\in\mathcal P(\Xi)$, $k\in [K]$, which are believed to be similar to~$\PP$. Each source distribution~$\PP_k$ can thus be approximated by an estimated distribution~$\widehat\PP_k\in\mathcal P(\Xi)$ constructed from the respective source data. In this situation, one might solve the DRO problem~\eqref{eq:dro} for a reference distribution~$\widehat \PP$ that is defined as a mixture or as an OT barycenter of the reference distributions~$\widehat\PP_k$, $k\in[K]$ \citep{lau2022wassersteinbarydro}. However, mixtures often fail to retain stylized features of the target distribution (see Example~\ref{ex:normal-barycenters}), and OT barycenters are plagued by the shortcomings highlighted in Section~\ref{sec:barycenters}. The following example shows that DRO problems with OT ambiguity sets around OT barycenters suffer from instability and are thus vulnerable to perturbations.

\begin{example}[Instability of DRO with OT Barycenters]
Consider a two-item newsvendor problem with $\Theta =[0,1]^2$, $\Xi=\mathbb R^2$ and loss function $\ell(\theta,\xi) = |\langle {\rm e}_1, (\theta-\xi)\rangle | + 2|\langle {\rm e}_2, \theta- \xi\rangle| + \langle v,\theta\rangle$, where $\rm e_1$ and $\rm e_2$ denote the standard basis vectors in~$\mathbb R^2$, and {$v=(\frac{1}{2}, 0)$}. Assume that there are $K=2$ estimated source distributions affected by a perturbation parameter~$\kappa\in[-1,1]$, that is,
\[
	\widehat{\mathbb P}_1^\kappa = \textstyle \frac{1}{2} \delta_{(\kappa,0)} + \frac{1}{2} \delta_{(1,1)} \quad \text{and} \quad \widehat{\mathbb P}_2^\kappa = \frac{1}{2} \delta_{(\kappa,1)} + \frac{1}{2} \delta_{(1,0)}.
\]
If we introduce { $\underline{\mathbb P}^{\kappa} = \frac{1}{2}\,\delta_{(\kappa,\,1/2)} + \frac{1}{2}\,\delta_{(1,\,1/2)}$ and $\overline{\mathbb P}{}^{\kappa} =\frac{1}{2}\,\delta_{((\kappa+1)/2,\,0)} + \frac{1}{2}\,\delta_{((\kappa+1)/2,\,1)}$} as shorthands, then one readily verifies that all equally weighted 2-Wasserstein barycenters of $\widehat{\mathbb P}_1^\kappa$ and $\widehat{\mathbb P}_2^\kappa$ are given by
\[
	\widehat{\mathbb{P}}^{\kappa} = \left\{\begin{array}{ll}
	\underline{\mathbb P}^{\kappa} & \text{if }\kappa<0, \\[-1ex]
	\eta\cdot \underline{\mathbb P}^{0}+(1-\eta)\cdot \overline{\mathbb P}{}^{0} ~~ \text{for any }\eta\in[0,1] ~~ & \text{if }\kappa=0,\\[-1ex]
	\overline{\mathbb P}{}^{\kappa}  & \text{if }\kappa>0.
	\end{array}\right.
\]
In particular, the barycenter is unique for~$\kappa\neq 0$. Consider now the distributionally robust newsvendor problem that minimizes the worst-case expected loss over all distributions in a $2$-Wasserstein ball of radius~$\varepsilon\geq 0$ around the barycenter~$\widehat{\mathbb P}^{\kappa}$, and define its optimal value function as $V(\kappa,\varepsilon) =\min_{\theta\in\Theta} \max_{{ \QQ} \in \mathbb{B}_\varepsilon (\widehat{\mathbb{P}}^{\kappa})} \int_{\Xi}\ell(\theta, \xi)\, {\rm d}{ \QQ}(\xi)$. By dualizing the inner maximization problem using \cite[Theorem~4.18]{kuhn2024distributionally} and by solving the resulting minimization problem in closed form, we find
\[
	V(\kappa,\varepsilon) = \min_{\theta\in\Theta} \int_{\Xi}\ell(\theta, \xi)\, {\rm d}\widehat{\mathbb{P}}^{\kappa}(\xi) +\sqrt{5} \cdot \varepsilon = \left\{ \begin{array}{ll}
	{ \frac{1-\kappa}{2}} +\sqrt{5} \cdot \varepsilon & \text{if } \kappa<0, \\[-1ex]
    { \frac{5+\kappa}{4}}
    +\sqrt{5} \cdot \varepsilon & \text{if } \kappa > 0.
	\end{array}\right.
\]
The factor $\sqrt{5}$ emerges because $\|\nabla_\xi \ell(\theta,\xi)\|_2= \sqrt{5}$ whenever the gradient exists.  One can further show that the optimal order quantity for the first product is uniquely given by $\theta_1^\star(\kappa,\varepsilon) = 0$ if $\kappa<0$ and by $\theta_1^\star(\kappa,\varepsilon) =\frac{1+\kappa}{2}$ if~$\kappa > 0$. Hence, neither the optimal value nor the optimal order quantity  admits a continuous extension at~$\kappa=0$. This example shows that DRO problems with ambiguity sets around Wasserstein barycenters may fail to be stable.\footnote{For simplicity, we assumed that $\Xi=\mathbb R^2$. However, instability persists even if we require~$\Xi$ to be convex and compact.} Small perturbations in the estimated source distributions can dramatically change the optimal solution, which is troubling for modelers.~\HOUSE
\end{example}

The flaws of traditional DRO models prompt us to introduce the multi-source DRO problem
\begin{equation}
\tag{MS-DRO}\label{eq:original-problem}
    \inf_{\theta\in\Theta} ~\sup_{{ \QQ}\in\cap_{k=1}^K\mathbb{B}_{\varepsilon_k}(\widehat{\mathbb{P}}_k)} ~ \int_{\Xi}\ell(\theta, \xi)\,{\rm d}{ \QQ}(\xi),
\end{equation}
which seeks to minimize the worst-case expected loss across all distributions in the intersection of~$K$ OT ambiguity sets. This ambiguity set contains all distributions that can be obtained by reshaping the reference distribution~$\widehat{\mathbb{P}}_k$ at a transportation cost of at most~$\varepsilon_k$ for every~$k\in[K]$. In contrast to the standard DRO problem~\eqref{eq:dro} with an ambiguity set around an OT barycenter, the multi-source DRO problem~\eqref{eq:original-problem} enjoys an attractive stability property. Indeed, as we rigorously prove in Appendix~\ref{app:stability}, its solutions change {\em continuously} under Wasserstein perturbations of the estimated source distributions $\widehat{\mathbb P}_k$, $k\in[K]$, provided that mild compactness and continuity conditions hold.

\begin{remark}[Guaranteeing Nonempty Intersections] Problem~\eqref{eq:original-problem} is only meaningful if the OT ambiguity sets $\mathbb{B}_{\varepsilon_k}(\widehat{\mathbb{P}}_k)$, $k\in[K]$, have a nonempty intersection. In practice, this can be ensured by selecting sufficiently large radii~$\varepsilon_k$, $k\in[K]$, using the following procedure. First, compute an OT barycenter~$\widehat{\mathbb{P}}$ of the reference distributions~$\widehat{\mathbb{P}}_k$ with respective weights~$\lambda_k>0$, $k\in[K]$. Recall from Lemma~\ref{lemma:barycenter:existence}
that~$\widehat{\mathbb{P}}$ exists whenever Assumptions~\ref{ass:assumptions-c} and~\ref{ass:barycenter-existence} hold. Next, define $\underline\varepsilon_k = C(\widehat{\mathbb{P}}_k, \widehat{\mathbb{P}})$, and select any~$\varepsilon_k \geq \underline\varepsilon_k$ for every $k\in[K]$. This choice guarantees that the ambiguity set in~\eqref{eq:original-problem} contains at least the barycenter~$\widehat{\mathbb{P}}$ and is thus nonempty. By~\cite{altschuler2021barycenterComplexity}, both the barycenter~$\widehat{\mathbb{P}}$ and the minimal radii~$\underline\varepsilon_k$ can be computed efficiently if all reference distributions are discrete and if either the number~$K$ of ambiguity sets or the dimension~$d$ of~$\xi$ is fixed.~\HOUSE
\end{remark}
 
In the remainder of this section we identify conditions under which the multi-source DRO problem~\eqref{eq:original-problem} is equivalent to a finite convex program susceptible to numerical solution. To this end, we first dualize the inner maximization problem over~{$\QQ$} to obtain an equivalent minimization problem that can be combined with the outer minimization over~$\theta$. For ease of exposition, we hide the dependence of the loss function on $\theta$ in the remainder of the theoretical analysis and only reintroduce it in Section~\ref{sec:experiments}. The inner maximization problem in~\eqref{eq:original-problem} thus simplifies to the following {\em primal uncertainty quantification problem}
\begin{equation}\tag{P-UQ}\label{eq:oringal-problem:inner-max}
    \sup_{{ \QQ}\in\cap_{k=1}^K\mathbb{B}_{\varepsilon_k}(\widehat{\mathbb{P}}_k)} ~\int_{\Xi}\ell(\xi)\, {\rm d}{ \QQ}(\xi). 
\end{equation}
Our convex reformulation results rely on the assumption that all reference distributions are discrete.

\begin{assumption}[Reference Distributions]
\label{ass:discrete-reference-distributions}
The $k$-th reference distribution is representable as 
$\widehat{\mathbb{P}}_k=\sum_{j=1}^{N_k}p_{k,j} \delta_{\widehat{\xi}_{k,j}}$ with atoms $\{\widehat{\xi}_{k,j}\}_{j=1}^{N_k}$ and strictly positive probabilities $\{p_{k,j}\}_{j=1}^{N_k}$ for all~$k\in[K]$.
\end{assumption}

In view of Assumption~\ref{ass:discrete-reference-distributions} it is useful to define $\MI=\times_{k=1}^K [N_k]$ as a set of $K$-dimensional multi-indices. Each $\alpha\in\mathcal A$ uniquely identifies a combination of atoms from the $K$ discrete reference distributions. In addition, we will also need the following growth condition on the loss function.
\begin{assumption}[Growth Condition]\label{ass:growth-condition}
    There exist $\xi'\in\Xi$, $g'>0$ and $p'\geq 1$  such that $\ell(\xi)\leq g'(1+d^{p'}(\xi,\xi'))$ for all $\xi\in\Xi$, where~$d$ is the same metric as in Assumption~\ref{ass:assumptions-c}.
\end{assumption}
We now derive the dual of the uncertainty quantification problem~\eqref{eq:oringal-problem:inner-max}. Going forward, we use~$\mathcal E$ to denote the set of all~$\varepsilon=(\varepsilon_1,\ldots,\varepsilon_K)\in \mathbb R_+^K$ for which the ambiguity set in~\eqref{eq:oringal-problem:inner-max} is nonempty.

\begin{theorem}[Strong Duality]\label{th_strong_duality}
If Assumptions~\ref{ass:assumptions-c}, \ref{ass:discrete-reference-distributions} and~\ref{ass:growth-condition} hold with $p'\leq p$, then the primal uncertainty quantification problem~\eqref{eq:oringal-problem:inner-max} admits the dual
\begin{equation}\tag{D-UQ}
\label{eq:dro_dual}
\begin{array}{c@{\quad}l@{\quad}l}
  \inf & \displaystyle \sum_{k=1}^K \varepsilon_k \lambda_k+\sum_{k=1}^K\sum_{j=1}^{N_k}p_{k,j} \, \gamma_{k,j}\\
  \rm{s.t.} & \lambda_k \in \mathbb{R}_{+},\ \gamma_k\in\mathbb{R}^{N_k} & \forall k\in[K]\\
  &\displaystyle \sup_{\xi\in\Xi}~ \ell(\xi)-\sum_{k=1}^K\lambda_k c(\xi,\wh \xi_{k, \mi_k}) \leq \sum_{k=1}^K\gamma_{k,\mi_k} & \forall \mi \in \MI.
\end{array}
\end{equation}
Strong (i.e., gap-free) duality holds for all~$\varepsilon\notin\partial\mathcal{E}$.
\end{theorem}

If each reference distribution has the same number of atoms, that is, if $N_k=N$ for all $k\in[K]$, then $|\MI|=N^K$. Hence, the number of constraints in problem~\eqref{eq:dro_dual} grows exponentially with~$K$.

Problem~\eqref{eq:dro_dual} can be viewed as a robust optimization problem with exponentially many robust constraints indexed by $\mi\in\MI$. Each constraint involves an embedded maximization problem over all uncertainty realizations~$\xi\in\Xi$, which is generically hard. Under the following convexity assumption, however, the robust constraints can be reformulated in terms of finitely many convex constraints.

\begin{assumption}[Convexity]
\label{assump_piecew_concave}
The support set~$\Xi$ is convex, the transportation cost function $c(\xi,\xi')$ is convex in $\xi$ for every fixed~$\xi'\in \Xi$, and the loss function is representable as $\ell(\xi)=\max_{l\in[L]}\ell_l(\xi)$ for some concave component functions $\ell_l:\Xi\rightarrow\mathbb R$, $l\in[L]$.
\end{assumption}

Assumption~\ref{assump_piecew_concave} is standard in DRO and is satisfied in many applications; see \cite[Section~5]{mohajerin2018data}. 
The next theorem follows from now standard reformulation techniques developed in~\citep{bental2015deriving, zhen2021mathematical} and is therefore stated without proof.

\begin{theorem}[Reformulation as a Finite Convex Program]\label{th_piecew_concave}
If Assumption~\ref{assump_piecew_concave} holds, then the dual uncertainty quantification problem~\eqref{eq:dro_dual} is equivalent to the finite convex program
\begin{equation}
\label{eq:convex-reformulation}
\begin{array}{c@{\quad}l@{~}l}
  \inf & \displaystyle \sum_{k=1}^K \varepsilon_k \lambda_k + \sum_{k=1}^K\sum_{j=1}^{N_k} p_{k,j}\, \gamma_{k,j}\\
      \rm{s.t.} & \lambda_k\in\mathbb{R}_+,\ \gamma_{k}\in\mathbb{R}^{N_k},\ u_{\mi,l}\in\mathbb{R}^d,\ v_{\mi,l}\in\mathbb{R}^d,\ w_{\mi,l,k}\in\mathbb{R}^d &\forall \mi \in \MI,\ \forall k\in[K],\ \forall l\in[L]\\
      & \displaystyle [-\ell_l]^*(u_{\mi,l}) +\sigma_{\Xi}(v_{\mi,l}) +\sum_{k=1}^K \lambda_kc^{*1}\left(\frac{w_{\mi,l,k}}{\lambda_k},\wh \xi_{k,\mi_k}\right) \leq \sum_{k=1}^K \gamma_{k,\mi_k} &\forall \mi \in \MI,\ \forall l\in[L]\\
      & \displaystyle u_{\mi,l} + v_{\mi,l}+\sum_{k=1}^Kw_{\mi,l,k}=0 &\forall \mi \in \MI,\ \forall l\in[L],
\end{array}
\end{equation}
where $c^{*1}(u,\wh \xi)$ is the  convex conjugate of the function $c(\xi,\wh\xi)$ with respect to its first argument~$\xi$.
\end{theorem}

We emphasize that problem~\eqref{eq:convex-reformulation} is indeed convex because the conjugate of the convex component function $-\ell_l$, the support function of the convex set~$\Xi$ and the perspective of the partial conjugate of the convex cost function~$c$ are all convex. However, the number of decision variables and constraints of problem~\eqref{eq:convex-reformulation} scales with~$|\MI|$ and is thus exponential in the number~$K$ of reference distributions.

The following corollary derives a simplification of problem~\eqref{eq:convex-reformulation} under polyhedrality conditions.

\begin{corollary}[Piecewise Affine Loss Function]\label{cor:pw-affine}
Suppose that $\Xi=\{ \xi : C\xi \leq g \}$ for some $C\in\mathbb{R}^{m\times d}$ and $g\in\mathbb{R}^m$, $c(\xi,\xi')=\|\xi-\xi'\|^p$ for some norm $\| \cdot \|$ on $\mathbb R^d$ and exponent~$p\in [1,\infty]$ and that $\ell(\xi)=\max_{l\in[L]}\langle a_l,\xi\rangle +b_l$ for some $a_l\in\mathbb R^d$ and $b_l\in\mathbb R$, $l\in[L]$. Then, problem~\eqref{eq:convex-reformulation} simplifies to
\[
\begin{array}{c@{\quad}l@{~}l}
    \inf & \displaystyle \sum_{k=1}^K \varepsilon_k \lambda_k +\sum_{k=1}^K\sum_{j=1}^{N_k}p_{k,j}\, \gamma_{k,j}\\
    \rm{s.t.} & \lambda_k\in\mathbb{R}_+,\ \gamma_{k}\in\mathbb{R}^{N_k},\ z_{\mi,l}\in\mathbb{R}_{+}^m,\ w_{\mi,l,k}\in\mathbb{R}^d &\forall \mi \in \MI,\ \forall k\in [K],\ \forall l\in[L] \\
    & \displaystyle b_l+ \sum_{k=1}^K \left( \langle w_{\mi,l,k},\wh \xi_{k,\mi_k}\rangle+\varphi(q) \,\lambda_k \left\| \dfrac{w_{\mi,l,k}}{\lambda_k}\right\|^{q}_{*}\right)  +\langle z_{\mi,l},g \rangle\leq \sum_{k=1}^K \gamma_{k,\mi_k} &\forall \mi \in \MI,\ \forall l\in[L]\\
    & \displaystyle a_l-C^T z_{\mi,l}=\sum_{k=1}^Kw_{\mi,l,k} &\forall \mi \in \MI,\ \forall l\in[L]
\end{array}
\]
where $\| \cdot \|_*$ is the norm dual to $\| \cdot \|$, $q\in [1,\infty]$ is the conjugate index of $p$ satisfying $1/p+1/q=1$, and $\varphi(q)=(q-1)^{q-1}/q^q$ for $q\in(1,\infty)$,  $\varphi(1)=1$ and $\varphi(\infty)=0$.
\end{corollary}

Note that Assumption~\ref{assump_piecew_concave} is automatically satisfied under the conditions of Corollary~\ref{cor:pw-affine}. Note also that Corollary~\ref{cor:pw-affine} extends \citep[Corollary~5.1]{mohajerin2018data} and \citep[\S~6]{zhen2021mathematical} to~$K>1$. The proof of Corollary~\ref{cor:pw-affine} relies on now standard techniques and is thus omitted.

Next, we show that problem~\eqref{eq:oringal-problem:inner-max} admits a worst-case distribution with sparse discrete support.

\begin{proposition}[Sparse Worst-Case Distribution]\label{prop:sparse-sol}
If Assumptions~\ref{ass:assumptions-c}, \ref{ass:discrete-reference-distributions} and~\ref{ass:growth-condition} hold, then problem~\eqref{eq:oringal-problem:inner-max} is solved by a discrete distribution supported on at most $1+\sum_{k=1}^K N_k$ points. 
\end{proposition}

Proposition~\ref{prop:sparse-sol} can be seen as a generalization of \cite[Corollary~2]{gao2023wasserstein}, which focuses on the special case when~$K=1$; see also \citep[Theorem~4]{yue2021linearWasserstein}. It is perhaps surprising that problem~\eqref{eq:oringal-problem:inner-max} admits a worst-case distribution with only $\mathcal O(\sum_{k=1}^KN_k)$ atoms even though the convex reformulation~\eqref{eq:convex-reformulation} of~\eqref{eq:oringal-problem:inner-max} involves $\mathcal O(\prod_{k=1}^K N_k)$ variables and constraints.

To close this section, we highlight the striking similarity between our dual uncertainty quantification problem~\eqref{eq:dro_dual} and the following dual of the OT barycenter problem~\eqref{eq:barycenter} with discrete marginals $\mathbb P_k=\widehat{\mathbb P}_k$, $k\in[K]$, which is derived in \citep[p.~5]{altschuler2021barycenterComplexity}.
\begin{equation}
    \label{eq:dual-barycenter-problem}
    \sup_{\gamma_{k}\in\mathbb{R}^{N_k} \, \forall k\in[K]} \left\{ \sum_{k=1}^K\sum_{j=1}^{N_k} p_{k,j}\, \gamma_{k,j} \;:\; \inf_{\xi\in\Xi}~\sum_{k=1}^K \lambda_k c(\xi,\widehat{\xi}_{k,\mi_k})\geq \sum_{k=1}^K \gamma_{k,\mi_k} ~ \forall\mi\in\MI \right\}
\end{equation}
This problem is equivalent to \eqref{eq:dro_dual} if the Lagrange multipliers $\lambda_k$, $k\in[K]$, are identified with the (fixed) weights of the reference distributions in~\eqref{eq:barycenter} and if the loss function~$\ell(\xi)$ is set to~0.

\section{Computational Complexity of Multi-Source DRO}
\label{sec:tractability}
We now demonstrate that the primal uncertainty quantification problem~\eqref{eq:oringal-problem:inner-max}, which constitutes a critical component of multi-source DRO, is generically intractable. We then show that~\eqref{eq:oringal-problem:inner-max} becomes tractable under widely used convexity assumptions if either the number~$K$ of OT ambiguity sets or the dimension~$d$ of the uncertain parameter~$\xi$ is kept constant. The main results of this section are inspired by~\citet{altschuler2021barycenterComplexity}. We thus relegate their proofs to Appendix~\ref{sec:proofs}.

We first prove that the following feasibility version of~\eqref{eq:oringal-problem:inner-max} is already NP-hard.

\begin{center}
\vspace{2ex}
\fbox{\parbox{0.98\columnwidth}{ {\centering
\vspace{0.5ex}\textsc{UQ Feasibility}\\}
		\textbf{Instance.} Discrete distributions~$\widehat{\mathbb P}_k\in\mathcal P(\mathbb R^d)$ and radii $\varepsilon_k\geq 0$ for all $k\in[K]$. \\[0.5ex]
		\textbf{Goal.} For $\Xi=\mathbb R^d$, $\ell(\xi)=0$ and $c(\xi,\xi')=\|\xi-\xi'\|_2^2$, decide if the primal uncertainty quantification problem~\eqref{eq:oringal-problem:inner-max} is feasible and, if not, find an infeasibility certificate given by a recession direction $(\lambda^\infty, \gamma^\infty)$ of the dual feasible set along which the objective function of~\eqref{eq:dro_dual} strictly decreases.}}
        \vspace{2ex}
\end{center}

Note that if $(\lambda^\infty, \gamma^\infty)$ is an infeasibility certificate, then $(\lambda^0,\gamma^0)+ t\cdot (\lambda^\infty, \gamma^\infty)$ is feasible in~\eqref{eq:dro_dual} for every~$(\lambda^0,\gamma^0)$ feasible in~\eqref{eq:dro_dual} and for every~$t\geq 0$, and
$\sum_{k=1}^K \varepsilon_k\lambda^\infty_k +\sum_{k=1}^K\sum_{j=1}^{N_k}p_{k,j} \, \gamma^\infty_{k,j}<0$. The \textsc{UQ Feasibility} problem is well-defined by virtue of Theorem~\ref{th_strong_duality}, which asserts that~\eqref{eq:oringal-problem:inner-max} and~\eqref{eq:dro_dual} are strong duals for all radii~$\varepsilon_k\geq 0$, $k\in[K]$. Indeed, this implies that an infeasibility certificate $(\lambda^\infty, \gamma^\infty)$ is guaranteed to exist whenever~\eqref{eq:oringal-problem:inner-max} is infeasible.

\begin{theorem}[NP-Hardness]\label{thm:nphard-informal}
The \textsc{UQ Feasibility} problem is NP-hard.
\end{theorem}

We now show that~\eqref{eq:oringal-problem:inner-max} becomes tractable under standard regularity conditions if either~$K$ or~$d$ is fixed. This finding mirrors the known tractability results for OT barycenters by~\citet{altschuler2021barycenterComplexity}. For simplicity of exposition, we assume throughout this discussion that $N_k=N$ for all~$k \in [K]$. We also assume that $\log U$ represents an upper bound on the number of bits needed to encode any radius~$\varepsilon_k$, any probability~$p_{k,j}$ and any component of~$\widehat{\xi}_{k,j}$ for all $j\in[N]$ and~$k\in[K]$. Our tractability results rely on three minimal assumptions. All of them stipulate that certain elementary computations involving the problem data can be carried out in polynomial time. We tacitly assume that the polynomials bounding the worst-case runtimes of these computations are {\em known}. The first assumption concerns the generalized Moreau envelope problem shown below.
\begin{equation}
    \label{eq:Moreau-envelope}
    \ell^\star(\lambda)~ = ~\sup_{\xi\in\Xi} ~ \ell(\xi)-\sum_{k=1}^K\lambda_k c(\xi,\widehat \xi_k)
\end{equation}

\begin{assumption}[Generalized Moreau Envelope]
\label{ass:compl-subproblem}
There exists an oracle for solving~\eqref{eq:Moreau-envelope}. If $\ell^\star(\lambda)<\infty$, then the oracle outputs an optimal solution of~\eqref{eq:Moreau-envelope}. If~$\ell^\star(\lambda)=\infty$, then the oracle outputs a closed halfspace $\Lambda\subseteq\mathbb R^K$ such that~$\lambda\notin \Lambda$ and $\lambda'\in \Lambda$ for every~$\lambda'\in\mathbb R^K$ with $\ell^\star(\lambda')<\infty$. The oracle runs in time polynomial in~$d$, $K$, $N$, $\log U$ and the bit sizes of~$\lambda_1,\ldots,\lambda_K\geq 0$ and~$\widehat \xi_1,\ldots, \widehat \xi_K\in\mathbb R^d$.
\end{assumption}

Assumption~\ref{ass:compl-subproblem} will allow us to construct an efficient separation oracle for the feasible set of problem~\eqref{eq:dro_dual}. It is satisfied, for example, if~$c$ is a convex function with efficiently computable subgradients, $\ell$ is the pointwise maximum of finitely many concave functions with efficiently computable supergradients and~$\Xi$ is a compact convex set that admits an efficient separation oracle. In this case, the generalized Moreau envelope problem~\eqref{eq:Moreau-envelope} is susceptible to the ellipsoid algorithm.

We emphasize that Assumption~\ref{ass:compl-subproblem} is restrictive in that it requires the supremum of problem~\eqref{eq:Moreau-envelope} to be attained and a corresponding maximizer to be computable {\em exactly} and in polynomial time. As we will see in Remark~\ref{rem:weak-oracles} below, the tractability results of this section remain valid even if we have only access to an {\em approximate} Moreau envelope oracle. For ease of exposition, however, we present the key ideas of this section under the simplifying assumption that an exact oracle is available.

The next assumption requires the loss function to admit an efficient zeroth-order oracle.

\begin{assumption}[Computational Properties of~$\ell$]\label{ass:loss-complexity}
    The loss $\ell(\xi)$ can be computed in time polynomial in the bit length of~$\xi\in\mathbb R^d$.
\end{assumption}

The last assumption concerns the transportation cost function~$c$. Before we can state this assumption, we first need to introduce the notion of a $c$-diagram. 

\begin{definition}[$c\,$-Diagram]
    \label{def:c-diagram}
    For any transportation cost function~$c$ satisfying Assumption~\ref{ass:assumptions-c}, the $c$-diagram associated with points $\xi_1, \ldots, \xi_{N}\in\Xi$ and weights $\gamma_1, \ldots, \gamma_{N}\in\mathbb R$ is the family of cells
    \[
        \Xi_i=\left\{\xi\in\Xi \;:\;c(\xi, \xi_i)-\gamma_i <c(\xi, \xi_{j})-\gamma_j ~~ \forall j\in[N],\; j\neq i \right\}\quad \forall i\in[N].
    \]
\end{definition}

Note that if~$c(\xi,\xi')=\|\xi-\xi'\|_2^2$ is the quadratic transportation cost function and all weights~$\gamma_i$, $i\in[N]$, are equal, then the $c$-diagram coincides with the Voronoi diagram induced by $\xi_i$, $i\in[N]$.

\begin{assumption}[Computational Properties of~$c$]
\label{ass:compl-costf}
The following hold.
\begin{enumerate}
    \item[(i)] The cost $c(\xi, \xi')$ can be computed in time polynomial in the bit length of~$\xi,\xi'\in\mathbb R^d$.
    \item[(ii)] If $\{\Xi_i\}_{i=1}^N$ is any $c$-diagram, then $\cup_{i=1}^N \text{\rm cl}(\Xi_i)=\Xi$.
    \item[(iii)] If $\{\Xi_{k,i}\}_{i=1}^N$, $k\in[K]$, are $K$ different $c$-diagrams, then the set $\MI'=\{\mi\in\MI: \mathcal C_\mi\neq\emptyset\}$ defined in terms of the cells $\mathcal C_\mi = \cap_{k=1}^K \Xi_{k,\mi_k}$, $\mi\in\MI$, can be computed in time $\poly(N,K,\log U)$.
\end{enumerate}
\end{assumption}
Assumption~\ref{ass:compl-costf}(\textit{i}) is unrestrictive and usually trivial to check. Assumptions~\ref{ass:compl-costf}(\textit{ii}) and~(\textit{iii}) will be instrumental for constructing a 
separation oracle for the feasible set of~\eqref{eq:dro_dual} that runs in time $\poly(N, K, \log U)$. One readily verifies that Assumption~\ref{ass:compl-costf}(\textit{ii}) holds whenever~$c$ is convex and satisfies Assumption~\ref{ass:assumptions-c}. By \citep[Lemma~18]{altschuler2021barycenterComplexity}, Assumption~\ref{ass:compl-costf}(\textit{iii}) is satisfied if $c(\xi,\xi')=\|\xi-\xi'\|_2^2$. In addition, it is also satisfied, for example, if~$c(\xi,\xi')$ is given by $\|\xi-\xi'\|_1$, $\|\xi-\xi'\|_\infty$ or $\|\xi-\xi'\|_2$; see \citep[Section~5.1]{altschuler2021barycenterComplexity}.

We are now ready to state our main tractability result for problem~\eqref{eq:dro_dual}.

\begin{theorem}[Tractability of~\eqref{eq:dro_dual}]
\label{thm:dual-complexity}
Suppose that Assumptions~\ref{ass:compl-subproblem}, \ref{ass:loss-complexity} and~\ref{ass:compl-costf}(i) hold, and there is~$\delta>0$ such that all $\delta$-optimal solutions of~\eqref{eq:dro_dual} belong to the ball of radius~$R$ around~$0$, where~$R$ is computable in time $\poly(N,d,K,\log U)$. Then, the optimal value of~\eqref{eq:dro_dual} can be computed to any accuracy~$\delta>0$ in time $\poly(N,d, \log U,\log 1/\delta)$ (with exponential dependence on~$K$). If additionally Assumptions~\ref{ass:compl-costf}(ii) and~\ref{ass:compl-costf}(iii) hold, then the optimal value of~\eqref{eq:dro_dual} can be computed to any accuracy~$\delta>0$ in time $\poly(N,K, \log U,\log 1/\delta)$ (with exponential dependence on~$d$).
\end{theorem}

Theorem~\ref{thm:dual-complexity} is proved by showing that---under appropriate conditions---the feasible set of problem~\eqref{eq:dro_dual} admits two different separation oracles with exponential dependence on~$d$ or~$K$, respectively. Thus, given either of the two separation oracles, \eqref{eq:dro_dual} can be solved with the ellipsoid algorithm \citep[Theorem~5.2.1]{ben2001lectures}. The polynomial time solvability of problem~\eqref{eq:dro_dual} for fixed~$K$ is expected in view of Theorem~\ref{th_piecew_concave}, which reformulates~\eqref{eq:dro_dual} as a convex program with $\mathcal{O}(N^K\cdot d\cdot L)$ decision variables and constraints. However, the insight that~\eqref{eq:dro_dual} can sometimes even be solved in time polynomial in~$K$ for fixed~$d$ is perhaps more surprising and mirrors recent results by~\cite{altschuler2021barycenterComplexity} for Wasserstein barycenters. The assumption that~$R$ can be computed in time $\poly(d,K,N,\log U)$ seems difficult to check. One can show, however, that this assumption holds if there exists a discrete distribution $\mathbb{P}_{\rm S}\in\mathcal P(\Xi)$ with $C(\mathbb{P}_{\rm S},\widehat{\mathbb P}_k)<\varepsilon_k$ for all~$k\in[K]$, which can be encoded in $\poly(d,K,N,\log U)$ bits, and if Assumption~\ref{ass:growth-condition} holds with a known growth constant~$g'$, which can be encoded in~$\log U$ bits. Details are omitted for brevity.

\begin{remark}[Weak Separation Oracles]
    \label{rem:weak-oracles}
    Assumption~\ref{ass:compl-subproblem} is too restrictive for many problems of practical interest. Indeed, it is easy to construct instances of~\eqref{eq:dro_dual} for which the optimal values of the embedded maximization problems are irrational or not even attained. In these cases, the Moreau envelope oracle cannot run in polynomial time or does not even exist, respectively. However, Theorem~\ref{thm:dual-complexity} remains valid if the Moreau envelope oracle of Assumption~\ref{ass:compl-subproblem} outputs only a $\delta$-optimal solution of the underlying optimization problem for any prescribed tolerance $\delta>0$~\cite[Corollary~4.2.7]{grotschel2012geometric}. This relaxation complicates the proof of Theorem~\ref{thm:dual-complexity} but is standard. \HOUSE
\end{remark}

 To close this section, we note that the computational complexity of problem~\eqref{eq:oringal-problem:inner-max} can be significantly reduced by constructing ambiguity sets from clustered data rather than raw samples from the source distributions. This could be achieved, for instance, by extending the clustering and ambiguity set construction methods of~\cite{wang2025mean} to multi-source DRO.


\section{Statistical Guarantees}
\label{sec:stat_guarantees}

We now investigate the relation between the multi-source DRO problem~\eqref{eq:original-problem} and the original stochastic optimization problem~\eqref{eq:sp}. To this end, we distinguish five scenarios characterized by increasingly realistic informational assumptions about the unknown target distribution~$\PP$. 

\paragraph{Scenario~1.} Assume that we know~$K$ source distributions~$\PP_k$, $k\in[K]$, which are all different from~$\PP$. In addition, even though~$\PP$ is unknown, the $p$-Wasserstein distance between~$\PP$ and~$\PP_k$ has a {\em known} upper bound~$r_k\geq 0$ for all~$k\in[K]$. As an example, imagine that a retailer expands into a new country, where the demands of the offered products follow an unknown distribution~$\PP$. However, the retailer is already active in~$K$ other countries and knows the respective demand distributions~$\PP_k$, $k\in[K]$. While demand distributions invariably differ across countries, it is reasonable to assume that each~$\PP_k$ provides {\em some} useful information about~$\PP$. Specifically, the retailer might know that~$W_p(\PP,\PP_k)\leq r_k$, where $W_p$ denotes the $p$-Wasserstein distance on~$\mathcal P_p(\Xi)$ induced by a metric~$d$ on~$\Xi$. Hence, $\PP$ is known to belong to~$\cap_{k=1}^K\mathbb{B}^p_{r_k}(\PP_k)$, where~$\mathbb{B}^p_{r_k}(\PP_k)=\{{ \QQ}\in\mathcal{P}(\Xi):W_p({ \QQ},\PP_k)\leq r_k\}$ is the $p$-Wasserstein ball of radius~$r_k$ around~$\PP_k$. The multi-source DRO problem with ambiguity set $\cap_{k=1}^K\mathbb{B}^p_{r_k}(\PP_k)$ thus provides a deterministic upper bound on the corresponding stochastic program.

\paragraph{Scenario~2.} As in scenario~1, we continue to assume that $W_p(\PP,\PP_k)$ has a {\em known} upper bound~$r_k\geq 0$. In contrast to scenario~1, however, we no longer assume that~$\PP_k$ is known. Instead, we only assume to have access to~$N_k$ independent samples~$\widehat\xi_{k,j}$, $j\in[N_k]$, from~$\PP_k$. In the remainder of this section, we denote the corresponding empirical distribution by $\widehat\PP_k=\frac{1}{N_k}\sum_{j=1}^{N_k} \delta_{\widehat \xi_{k,j}}$. Note that~$\widehat\PP_k$ constitutes a random object governed by the joint distribution~$\PP_k^{N_k}$ of the given independent samples from~$\PP_k$. It is well known that $\widehat \PP_k$ approximates~$\PP_k$ in $p$-Wasserstein distance provided that~$\PP_k$ is light-tailed.

\begin{definition}[Light-Tailed Distributions]
\label{def:light_tails}
A distribution $\PP\in\mathcal P(\Xi)$ is $(a,A)$-light tailed for some positive constants $a,A>0$ if $\mathbb E_{\PP}[\exp(\|\xi \|^{a})]\leq A$.
\end{definition}

Most of the subsequent results will depend on the following assumption.

\begin{assumption}[Light-Tailed Reference Distributions]
\label{light_tailed_assumption}
There exist $a>p$ and $A>0$ such that all reference distributions $\PP_k$, $k\in[K]$, are $(a,A)$-light tailed.
\end{assumption}

\begin{lemma}[Concentration Inequalities {\cite{Fournier2015}}]
\label{lem:concentration-empirical-FG}
Suppose that $p \neq d/2$ and that Assumption~\ref{light_tailed_assumption} holds. Then, there exist constants $c_1, c_2>0$ that depend only on~$a$, $A$ and~$d$ such that $\PP_k^{N_k}[\PP_k \in \mathbb{B}^p_{\eps_k}(\widehat \PP_k)] \geq 1-\beta(\eps_k,N_k)$ for all $\eps_k\in\mathbb R_+$ and~$N_k\in\mathbb N$, where 
\begin{equation*}
	\beta(\eps_k,N_k) = \left\{ \begin{array}{ll} \displaystyle c_1\exp\left( -c_2 N_k \eps_k^{\max\{d/p,2\}}\right) & \displaystyle \text{if } \eps_k \le 1, \\[0ex]
	\displaystyle c_1\exp\left( -c_2 N_k \eps_k^{a/p}\right) & \displaystyle \text{if } \eps_k >1. \end{array}\right.    
\end{equation*}
\end{lemma}
Lemma~\ref{lem:concentration-empirical-FG} readily extends to~$p = d/2$ at the expense of a more intricate formula for $\beta(\eps_k,N_k)$. Given a fixed significance level $\beta_k\in(0,1]$, one can solve the equation $\beta(\eps_k,N_k)=\beta_k$ for~$\eps_k$ to obtain 
\begin{equation*}
\eps(\beta_k,N_k) = \left\{ \begin{array}{ll}
\displaystyle \left( \log(c_1/\beta_k) / (c_2 N_k) \right)^{\min\{p/d,\,1/2\}} 
& \displaystyle \text{if } N_k \ge \log(c_1/\beta_k)/c_2, \\[1ex]
\displaystyle \left( \log(c_1/\beta_k) / (c_2 N_k) \right)^{p/a} 
& \displaystyle \text{if } N_k < \log(c_1/\beta_k)/c_2.
\end{array} \right.
\end{equation*}
Lemma~\ref{lem:concentration-empirical-FG} then implies that the $p$-Wasserstein ball $\mathbb{B}^p_{\eps_k}(\widehat \PP_k)$ represents a $(1-\beta_k)$-confidence set for~$\PP_k$ whenever its radius satisfies~$\eps_k\geq \eps(\beta_k,N_k)$; see also \citep[Theorem~18]{kuhn2019wasserstein}.

The following proposition leverages Lemma~\ref{lem:concentration-empirical-FG}  to construct a confidence set for the unknown distribution~$\PP$. This is possible even though there is no data from~$\PP$.  In the following, we use $\PP_{[K]} = \otimes_{k=1}^K\mathbb{P}_k^{N_k}$ as a shorthand for the joint distribution of the samples $\widehat\xi_{k,j}$ for all~$j\in[N_k]$ and~$k\in[K]$.

\begin{proposition}[Known Distribution Shifts]
\label{prop:conc_ineq_balls}
Suppose that $p \neq d/2$ and that Assumption~\ref{light_tailed_assumption} holds. If $W_p(\mathbb{P}_k,\mathbb{P})\leq r_k$ and if $\eps_k\geq r_k+\eps(\beta_k,N_k)$ for some $\beta_k\in(0,1]$ and for all $k\in[K]$, then
\begin{equation*}
    \PP_{[K]} \left( \mathbb{P} \in \cap_{k=1}^{K} \mathbb{B}^p_{\eps_k}(\widehat{\mathbb{P}}_k) \right)\geq 1 - \sum_{k=1}^K \beta_k.
\end{equation*}
\end{proposition}
The ambiguity set $\cap_{k=1}^{K} \mathbb{B}^p_{\eps_k}(\widehat{\mathbb{P}}_k)$ inherits the randomness of the empirical distributions. In particular, it may be empty with a small positive probability. Depending on the realizations of~$\widehat \PP_k$, $k\in[K]$, one should thus select sufficiently large radii~$\eps_k\geq r_k+\eps(\beta_k,N_k)$ that lead to a nonempty ambiguity set.

Proposition~\ref{prop:conc_ineq_balls} implies that the optimal value of the multi-source DRO problem~\eqref{eq:original-problem} with ambiguity set $ \cap_{k=1}^{K} \mathbb{B}^p_{\eps_k}(\widehat{\mathbb{P}}_k)$ provides a $(1-\sum_{k=1}^K \beta_k)$-upper confidence bound on the optimal value of the original stochastic optimization problem~\eqref{eq:sp}. Note that different choices of the parameters~$\beta_k$, $k\in[K]$, yield different bounds with the same confidence level. Hence, these bounds can be tuned.

\paragraph{Scenario~3.} As in scenario~2, we continue to assume that $\PP$ as well as all reference distributions~$\PP_k$, $k\in[K]$, are unknown. In contrast to scenario~2, however, we now assume that~$\PP=\PP_1$ and that we have {\em no} information about the $p$-Wasserstein distance between~$\PP$ and~$\PP_k$ for~$k\geq 2$. Instead, we only have access to independent samples from all source distributions and can thus construct  empirical distributions~$\widehat \PP_k$, $k\in[K]$, as in scenario~2. It is now natural to set $\eps_k=W_p(\widehat \PP_1, \widehat \PP_k) +\eps(\beta_1,N_1)$, $k\in[K]$, for some~$\beta_1\in(0,1]$. In this case, one can show that the ambiguity set $\cap_{k=1}^{K} \mathbb{B}^p_{\eps_k}(\widehat{\mathbb{P}}_k)$ contains~$\PP$ with probability at least $1-\beta_1$. To see this, note first that the inequality $W_p({\mathbb{Q}}, \widehat\PP_1)\leq { \eps_1}$ implies
\[
    W_p({\mathbb{Q}},\widehat \PP_k)\leq W_p({\mathbb{Q}},\widehat \PP_1) + W_p(\widehat \PP_1,\widehat \PP_k) \leq W_p(\widehat \PP_1,\widehat \PP_k)+ { \eps_1}= { \eps_k},
\]
where the equality follows from the definitions  of~$\eps_1$ and~$\eps_k$. Therefore, we have $\mathbb{B}^p_{\eps_1}(\widehat{\mathbb{P}}_1) \subseteq \mathbb{B}^p_{\eps_k}(\widehat{\mathbb{P}}_k)$ for all $k\in[K]$, and thus the ambiguity set $\cap_{k=1}^{K} \mathbb{B}^p_{\eps_k}(\widehat{\mathbb{P}}_k)$ collapses to the crisp Wasserstein ball~$\mathbb{B}^p_{\eps_1}(\widehat{\mathbb{P}}_1)$. Lemma~\ref{lem:concentration-empirical-FG} then implies that $\PP\in \cap_{k=1}^{K} \mathbb{B}^p_{\eps_k}(\widehat{\mathbb{P}}_k)$ with probability at least $1-\beta_1$. Thus, \eqref{eq:original-problem} provides again an upper confidence bound on~\eqref{eq:sp}. However, the data from the reference distributions~$\PP_k$ with $k\neq 1$ is wasted. It is unclear how to make use of this data without additional assumptions.

\paragraph{Scenario~4.} As in scenario~3, we assume that $\PP=\PP_1$ and that the reference distributions~$\PP_k$, $k\in[K]$, are only accessible via data. In contrast to scenario~3, however, we now adopt a Bayesian perspective, that is, we assume that the decision-maker has a prior belief about the Wasserstein distances between the reference distributions. Arguably, the main motivation for using data from reference distributions is the belief that $\mathbb{P}$ is in some sense ``close'' to~$\mathbb{P}_k$ (and thus also to~$\widehat{\mathbb{P}}_k$).

In the following, we treat not only the samples $\{\widehat\xi_{k,j}\}_{j=1}^{N_k}$, $k\in[K]$, but also the reference distributions~$\PP_k$, $k\in[K]$, as random objects on an abstract probability space~$(\Omega,\mathcal F,\mu)$. We assume that the distributions~$\PP_k$, $k\in[K]$, satisfy Assumption~\ref{light_tailed_assumption} $\mu$-almost surely. The only property of~$\mu$ to be used below is that, conditional on any fixed realizations of~$\PP_k$, $k\in[K]$, all samples are mutually independent under~$\mu$, and the conditional distribution of~$\widehat\xi_{k,j}$ under~$\mu$ coincides with~$\PP_k$ for all~$j\in[N_k]$ and~$k\in[K]$. We also define~$F_k$ as the cumulative distribution function of~$W_p(\PP_1,\PP_k)$ under~$\mu$. This function captures the decision-maker's prior beliefs about the proximity between~$\PP_1$ and~$\PP_k$.

In Bayesian statistics, the prior beliefs are updated when new evidence emerges. In the example at hand, the decision-maker observes the $p$-Wasserstein distance between~$\widehat\PP_1$ and~$\widehat\PP_k$, which we denote as~$\widehat r_k\in\mathbb R_+$. The prior beliefs can then be updated by conditioning on the event~$W_p(\widehat\PP_1,\widehat\PP_k)\leq \widehat r_k$.

\begin{remark}
Note that conditioning on the event~$W_p(\widehat\PP_1,\widehat\PP_k) = \widehat r_k$ would provide a stronger update of the prior. We nevertheless prefer to condition on the less informative event~$W_p(\widehat\PP_1,\widehat\PP_k)\leq \widehat r_k$ for two reasons. First, it enhances the transparency of the subsequent derivations and obviates the need to mobilize advanced measure-theoretic tools. In addition, we aim to derive an upper confidence bound on $W_p(\PP_1,\widehat\PP_k)$. As $\PP_1$ is close to $\widehat\PP_1$ by virtue of Lemma~\ref{lem:concentration-empirical-FG}, such a confidence bound is within reach if $\widehat\PP_1$ and $\widehat\PP_k$ are known to be close. The extra information that $\widehat\PP_1$ and $\widehat\PP_k$ are ``not \emph{too} close,'' which is available in the event~$W_p(\widehat\PP_1,\widehat\PP_k) = \widehat r_k$, cannot meaningfully strengthen this bound. \HOUSE
\end{remark}

\begin{proposition}[Bayesian Measure Concentration]
\label{lemma:bayesian-fixed_N}
If~$\mu$ and~$F_k$, $k\in[K]$, are constructed as above, then we have $\mu(\PP_1 \in \mathbb{B}^p_{\eps_k}(\widehat \PP_k) | W_p( \widehat\PP_1,\widehat\PP_k)\leq\widehat r_k) \geq 1-\beta_{F_k}(\eps_k,\widehat r_k,N_1,N_k)$ for all $\eps_k, \widehat r_k \in\mathbb R_+$ with $\eps_k\geq \widehat r_k$ and for all $N_1,N_k\in\mathbb N$, $k\in[K]$, where 
\begin{equation}
\label{eq:beta_F_k}
    \beta_{F_k}(\eps_k,\widehat r_k,N_1,N_k) = \frac{\beta(\eps_k- \widehat r_k,N_1) \cdot \left(\int_0^{\eps_k} \beta(\eps_k-r_k,N_k) \, {\rm d}F_k(r_k)+1-F_k(\eps_k) \right) }{ \mu \left(W_p( \widehat\PP_1,\widehat\PP_k)\leq\widehat r_k\right) }.
\end{equation}
\end{proposition}

We highlight that $\beta_{F_k}(\eps_k,\widehat r_k,N_1,N_k)$ is non-increasing in~$\eps_k$. Indeed, the denominator in~\eqref{eq:beta_F_k} is positive and independent of~$\eps_k$. In addition, both terms in the numerator are nonnegative and non-increasing in~$\eps_k$ whenever $\eps_k\geq \widehat{r}_k$. Note that the derivative of the second term can be expressed as $\int_0^{\eps_k} \partial_\eps\beta(\eps_k-r_k,N_k) \, {\rm d}F_k(r_k)$ by virtue of the Reynolds theorem. Hence, it is strictly negative. Proposition~\ref{lemma:bayesian-fixed_N} and the union bound imply that, given any fixed $\beta_k\in(0,1]$, $k\in[K]$, we have
\begin{equation*}
    \mu \left( \left. \mathbb{P} \in \cap_{k=1}^{K} \mathbb{B}^p_{\eps_k}(\widehat{\mathbb{P}}_k) \right| W_p( \widehat\PP_1,\widehat\PP_k)\leq\widehat r_k\right)\geq 1 - \sum_{k=1}^K \beta_k
\end{equation*}
whenever $\eps_k\geq \eps_{F_k}(\beta_k, \widehat r_k, N_1, N_k)$ for all $k\in[K]$, where
\[
    \eps_{F_k}(\beta_k, \widehat r_k, N_1, N_k)=\inf \left\{\eps_k\geq \widehat r_k\; :\; \beta_{F_k}(\eps_k,\widehat r_k,N_1,N_k)\leq \beta_k\right\}.
\]
Thus, \eqref{eq:original-problem} provides a meaningful upper confidence bound on~\eqref{eq:sp} that does not waste data. 

\begin{figure}
    \centering
    \begin{tikzpicture}
        \node[anchor=south west,inner sep=0] (image) at (0,0) {\includegraphics[width=0.3\linewidth]{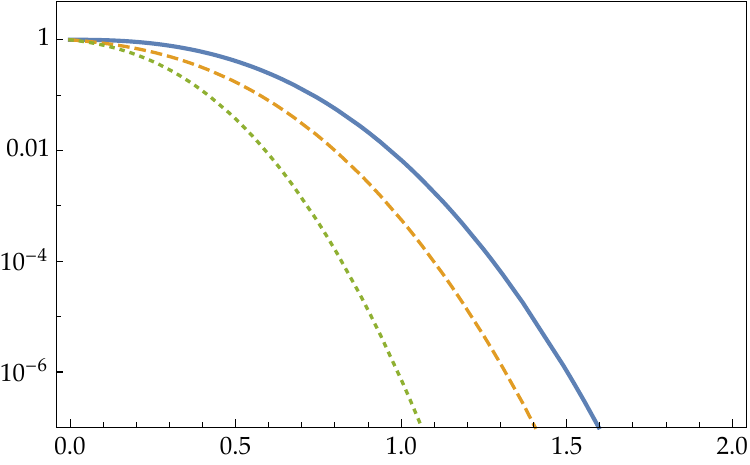}};

        \node at (2.6,0.1)[below] {\scriptsize{$c_2^{p/a}(\eps_k-\widehat r_k)$}};

        \node at (-0.4,2.8) [left, rotate=90] {\scriptsize{$\frac{\beta_{F_k}(\eps_k,\widehat r_k,N_1,N_k)}{\beta_{F_k}(\widehat r_k,\widehat r_k,N_1,N_k)}$}};
    \end{tikzpicture}
    \caption{Normalized Bayesian significance level as a function of the normalized Wasserstein radius $c_2^{p/a}(\eps_k-\widehat r_k)$ when there is no prior (solid blue line), a weak prior (dashed orange line) and a strong prior (dotted green line).}
    \label{fig:bayesian_effect}
\end{figure}

Figure~\ref{fig:bayesian_effect} shows the normalized Bayesian significance level $\beta_{F_k}(\eps_k,\widehat r_k,N_1,N_k)/\beta_{F_k}(\widehat r_k,\widehat r_k,N_1,N_k)$ as a function of the normalized Wasserstein radius $c_2^{p/a}(\eps_k-\widehat r_k)$ for the $2$-Wasserstein distance on~$\mathcal P_2(\mathbb R^5)$ when there are $N_1=5$ target samples and $N_k=50$ source samples. Here, $\PP_1$ and~$\PP_k$ are light-tailed with critical exponent~$a=5$, and the $2$-Wasserstein distance between~$\PP_1$ and~$\PP_k$ is believed to follow a normal distribution with mean~$\widehat r_k$ and variance $\sigma^2=0.2c_2^{-2p/a}$ (strong prior), $\sigma^2=0.5c_2^{-2p/a}$ (weak prior) or $+\infty$ (no prior). The chosen normalizations ensure that the curves in Figure~\ref{fig:bayesian_effect} are independent of the unknown constants~$c_1$, $c_2$ and $\mu (W_2( \widehat\PP_1,\widehat\PP_k)\leq\widehat r_k)$. As expected, we observe that a stronger prior leads to stronger statistical guarantees. 

\paragraph{Scenario~5.} Consider now a special case of scenario~4, which arises when the decision-maker observes {\em no} data from the target distribution~$\mathbb{P}_1=\PP$ at all (that is, when~$N_1=0$). From the proof of Proposition~\ref{lemma:bayesian-fixed_N} we already know that the prior probability obeys the estimate
\begin{equation}
    \label{eq:scenario5-guarantee}
    \mu\left(W_p(\mathbb{P}_1,\widehat{\PP}_k)>\varepsilon_k\right) \leq \int_0^{\eps_k} \beta(\eps_k-r_k,N_k) \, {\rm d}F_k(r_k) + 1-F_k(\eps_k).
\end{equation}
As there is no data from~$\PP_1$, this prior probability cannot be updated. However, given any significance level~$\beta_k\in(0,1]$, one can increase~$\eps_k$ until the right hand side of~\eqref{eq:scenario5-guarantee} drops below~$\beta_k$. Constructing~$\eps_k$ in this manner for every~$k\in[K]$ enables us again to show that~\eqref{eq:original-problem} provides a $(1-\sum_{k=1}^K \beta_k)$-upper confidence bound on the optimal value of the original stochastic program.

Further information about the relationships between the five scenarios is given in Appendix~\ref{sec:scenario-relationships}.

\section{Experiments}\label{sec:experiments}

Consider the problem of investing a given amount of capital into
$d$~assets with uncertain rates of return~$\xi\in\mathbb R^d$. The portfolio weight vector~$\theta\in\mathbb R^d$ collects the percentage weights of the capital allocated to the available assets. In the absence of short-sales, it ranges over the $d$-dimensional unit simplex~$\Theta$. The goal is to minimize the expected value of the portfolio loss~$-\langle \theta, \xi\rangle$ adjusted by the portfolio risk, measured by the conditional value-at-risk (CVaR) of $-\langle \theta, \xi\rangle$ at level $\cvar\in (0,1]$. Denoting the investor's risk aversion as $\rho\geq 0$, the portfolio selection problem can thus be expressed~as
$$
    \inf_{\theta\in\Theta}~ \mathbb E_\PP\left[-\langle \theta, \xi\rangle\right] + \rho \cdot \PP\text{-CVaR}_\cvar (-\langle \theta, \xi\rangle),
$$
where $\mathbb{P}$ denotes the target distribution of the asset returns. Replacing the CVaR by its definition due to \citet{rockafellar2000cvar}, this problem simplifies to
\begin{equation}\label{sec:Exp:portfolio_problem}
    \inf_{\theta\in\Theta, \, \tau\in\mathbb{R}} ~ \int_\Xi \max_{t\in\{1,2\}} \left\{ a_t\langle \theta, \xi\rangle +b_t\tau\right\}\,\rm d\mathbb{P}(\xi),
\end{equation}
where $a_1=-1$, $a_2=-1-\frac{\rho}{\cvar}$, $b_1=\rho$, and $b_2=\rho(1-\frac{1}{\cvar})$; see, {\em e.g.}, \citep[\S~7.1]{mohajerin2018data}. We henceforth assume that the target distribution~$\mathbb{P}$ is unknown but close to $K=2$ source distributions~$\PP_1$ and~$\PP_2$. While~$\PP_1$ and~$\PP_2$ are unknown, too, the investor has access to respective empirical distributions~$\widehat\PP_1$ and~$\widehat\PP_2$ and can thus construct the following instance of~\eqref{eq:original-problem}.
\begin{equation}
\label{eq:portfolio-problem}
    \inf_{\theta\in\Theta, \,\tau\in\mathbb{R}} ~ \sup_{{ \QQ}\in\cap_{k=1}^2\mathbb{B}_{\varepsilon_k}(\widehat{\mathbb{P}}_k)} ~ \int_\Xi \max_{t\in\{1,2\}} \left\{ a_t\langle \theta, \xi\rangle +b_t\tau\right\}\,\rm d{ \QQ}(\xi).
\end{equation}
The ambiguity set in~\eqref{eq:portfolio-problem} is defined as the intersection of two 1-Wasserstein balls induced by the transportation cost function~$c(\xi, \xi')=\|\xi-\xi'\|_1$. By Corollary~\ref{cor:pw-affine}, problem~\eqref{eq:portfolio-problem} can be reformulated as a tractable linear program and is thus susceptible to highly efficient optimization algorithms.

\subsection{Sensitivity Analysis}
\label{sec:sensitivity}
\begin{figure}
     \centering
     \begin{subfigure}[b]{0.45\textwidth}
         \centering
         \includegraphics[width=\textwidth, trim={0 0.3cm 0 5cm},clip]{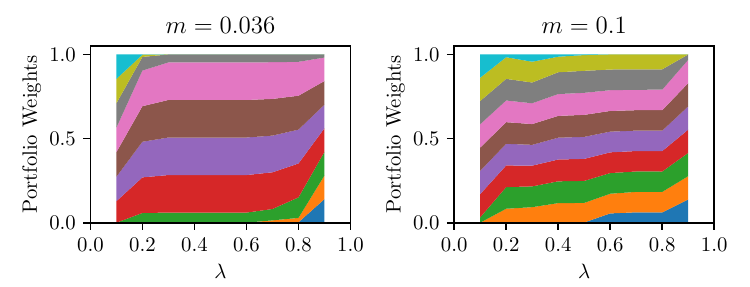}
         \caption{Sensitivity with respect to $\lambda$}
         \label{fig:portfolio-results:synthetic:lambda}
     \end{subfigure}
     \hfill
     \begin{subfigure}[b]{0.45\textwidth}
         \centering
         \includegraphics[width=\textwidth,trim={0 0.3cm 0 5cm},clip]{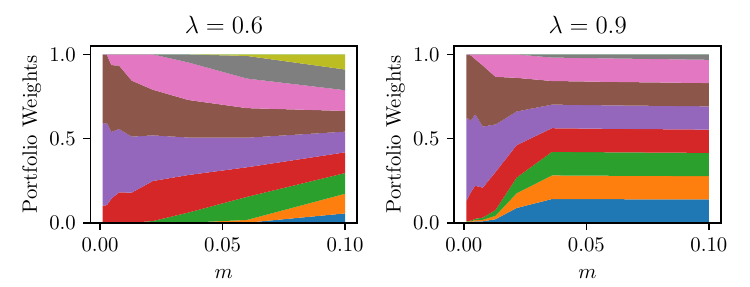}
         \caption{Sensitivity with respect to $m$}
         \label{fig:portfolio-results:synthetic:m}
     \end{subfigure}
    \caption{Impact of~$m$ and~$\lambda$ on the portfolio weights. The assets are ordered from bottom (asset~1) to top (asset~10).}
    \label{fig:portfolio-results:synthetic}
\end{figure}

In the first experiment we use synthetic data to analyze the impact of the radii~$\eps_1$ and~$\eps_2$ on the optimal portfolio. Throughout this experiment we assume that there are $d=10$ assets. We further assume that the asset returns are mutually independent under both source distributions and that $\xi_{i}\sim\mathcal{N}(i\%, 1\%)$ under~$\PP_1$, whereas $\xi_{i}\sim\mathcal{N}((11-i)\%, 1\%)$ under~$\PP_2$ for all~$i\in[d]$.
We also measure risk by the CVaR at level~$\cvar=20\%$ and set~$\rho=10$. In addition, we construct two empirical distributions~$\widehat{\mathbb{P}}_1$ and~$\widehat{\mathbb{P}}_2$ from $N_1=N_2=30$ independent training samples from~$\PP_1$ and~$\PP_2$, respectively. Finally, we set the radii of the two 1-Wasserstein balls to $\varepsilon_1=\lambda \cdot (1+m) \cdot W_1(\widehat{\mathbb{P}}_1, \widehat{\mathbb{P}}_2)$ and $\varepsilon_2=(1-\lambda) \cdot (1+m)\cdot W_1(\widehat{\mathbb{P}}_1, \widehat{\mathbb{P}}_2)$, where~$m\in\mathbb{R}_+$ captures the absolute and~$\lambda\in [0,1]$ the relative degree of ambiguity of the two source distributions. This construction ensures that $\cap_{k=1}^2 \mathbb B_{\eps_k}(\widehat \PP_k)$ is nonempty. Figure~\ref{fig:portfolio-results:synthetic} visualizes the dependence of the optimal portfolio weights obtained from problem~\eqref{eq:portfolio-problem} on~$\lambda$ and~$m$. Specifically, Figure~\ref{fig:portfolio-results:synthetic:lambda} shows that, as~$\lambda$ increases (decreases), assets with low (high) indices are assigned more weight. Indeed, the ambiguity set of problem~\eqref{eq:portfolio-problem} degenerates to~$\{\widehat\PP_2\}$ as~$\lambda$ tends to~$1$ and to~$\{\widehat\PP_1\}$ as~$\lambda$ tends to~$0$, and assets with high (low) indices perform well under~$\PP_1\approx\widehat\PP_1$ ($\PP_2\approx\widehat\PP_2$). Figure~\ref{fig:portfolio-results:synthetic:m} shows that, as~$m$ increases, the portfolio weights all become approximately equal, that is, the equally-weighted portfolio appears to be optimal under high ambiguity. This observation is in agreement with theoretical findings for distributionally robust single-source portfolio selection problems due to \cite{pflug2012}; see also \cite[Proposition~7.2]{mohajerin2018data}.

\subsection{Backtest on Synthetic Data}
\label{sec:synthetic-data}
In the second experiment we compare the optimal portfolios resulting from the multi-source DRO problem~\eqref{eq:portfolio-problem} against several other distributionally robust portfolios in terms of their out-of-sample means, standard deviations and Sharpe ratios. Similar to Section~\ref{sec:sensitivity}, we assume that the asset returns~$\xi_i$ for $i \in [10]$ follow independent normal distributions $\xi_i \sim \mathcal{N}(r_i, 1\%)$, and we use two different models for the expected returns. In model~1, we set~$r_i=+0.4\%$ if~$i\leq 2$ and~$r_i=-0.4\%$ if~$i> 2$. In model~2, on the other hand, we set~$r_i=+0.4\%$ if~$i\leq C$ and~$r_i=-0.4\%$ if~$i> C$ for some \textit{cutoff index}~$C$. The distributions~$\PP_1$ and~$\PP_2$ are then defined as the distribution of~$\xi$ in model~1 and model~2, respectively. We assume that~$\PP_1$ coincides with the target distribution~$\PP$ and that we have access to~$N_1=5$ samples from~$\PP_1$ and to~$N_2=30$ samples from~$\PP_2$. All other problem parameters are chosen as in Section~\ref{sec:sensitivity}. 

Note that~$C$ separates low-performing from high-performing assets in the source distribution. If~$C = 2$, then~$\PP_2 = \PP_1$. As~$C$ deviates further from~$2$, the distribution shift from~$\PP_2$ to~$\PP_1$ becomes more pronounced. The purpose of this stylized example is to show that if data from the target distribution~$\PP_1$ is scarce, then using data from the source distribution~$\PP_2$ can improve the performance of the optimal portfolio. Intuitively, the (abundant) source data provides an initial guess for the unknown target distribution, and the (scarce) target data allows us to improve this initial guess.

We compare the optimal portfolios of the multi-source DRO problem~\eqref{eq:portfolio-problem} against several baselines. Every baseline portfolio is obtained by solving a single-source version of problem~\eqref{eq:portfolio-problem}, which involves only a single 1-Wasserstein ball. The center of this ball is set to the empirical distribution on the \textit{target data}, the \textit{source data}, the \textit{pooled data} (that is, the union of the target and source data), or a \textit{$1$-Wasserstein barycenter} of the empirical distributions corresponding to the target and the source data. The latter approach is inspired by~\citep{lau2022wassersteinbarydro}. We set the Wasserstein radii in problem~\eqref{eq:portfolio-problem} to $\varepsilon_1=m+\lambda\cdot W_1(\widehat \PP_1, \widehat\PP_2)$ and $\varepsilon_2= m+(1-\lambda)\cdot W_1(\widehat \PP_1, \widehat\PP_2)$, where~$\lambda$ ranges from~0 to~1 in steps of 0.1, and~$m \in M\triangleq \{0.001, 0.005, 0.01, 0.05, 0.1, 0.2, 0.5, 1\}$. The Wasserstein radius in the single-source DRO problems that are used to construct the baseline portfolios is also chosen from~$M$. In all cases, we select hyperparameters from the underlying search grids that maximize the average portfolio return across five (additional) independent validation samples from~$\PP_1$.

\begin{figure}
     \centering
     \begin{subfigure}[b]{0.3\textwidth}
         \centering
         \includegraphics[width=\textwidth]{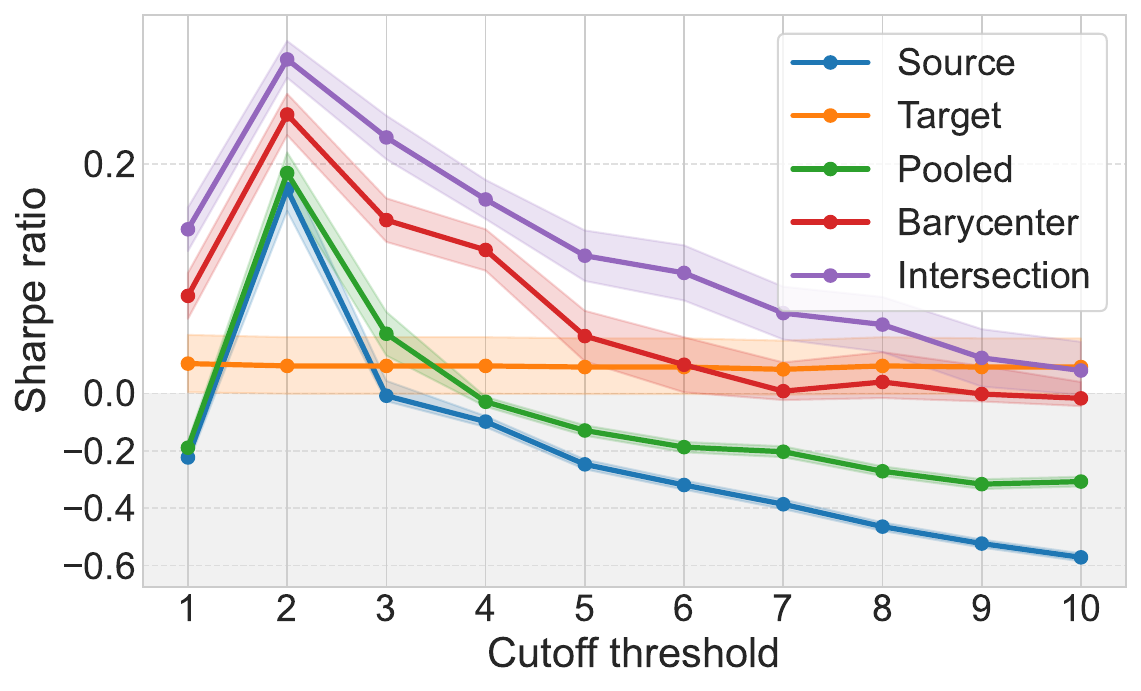}
         \caption{Sharpe Ratio}
     \end{subfigure}
     \begin{subfigure}[b]{0.3\textwidth}
         \centering
         \includegraphics[width=\textwidth]{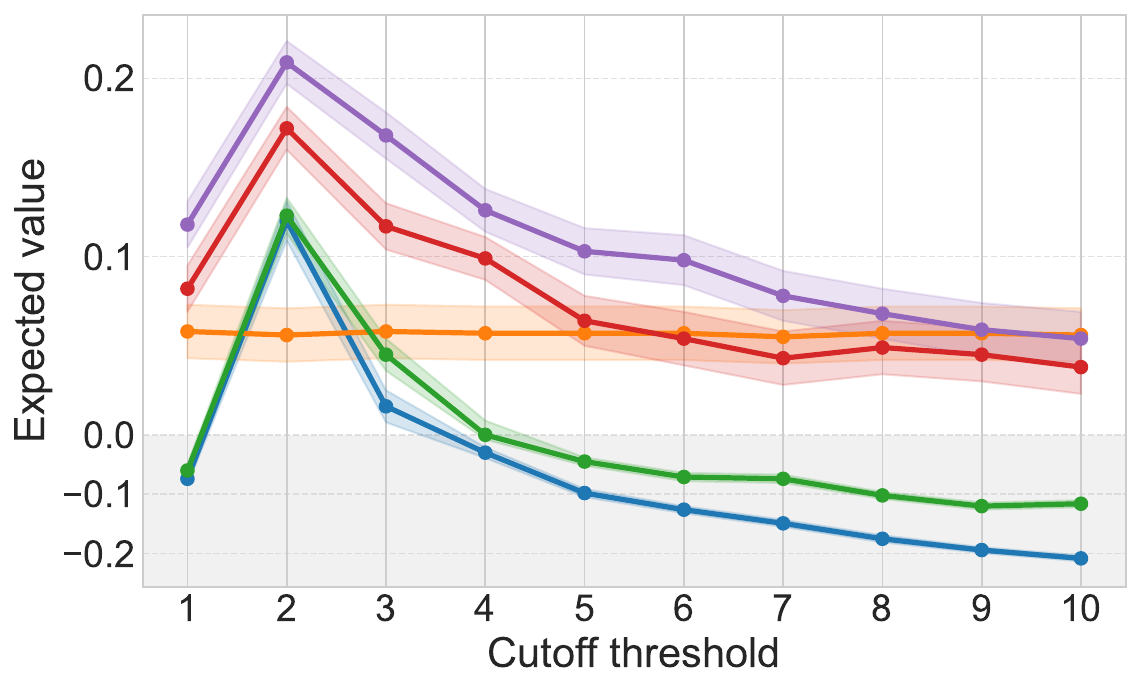}
         \caption{Expected Return}
     \end{subfigure}
     \begin{subfigure}[b]{0.3\textwidth}
         \centering
         \includegraphics[width=\textwidth]{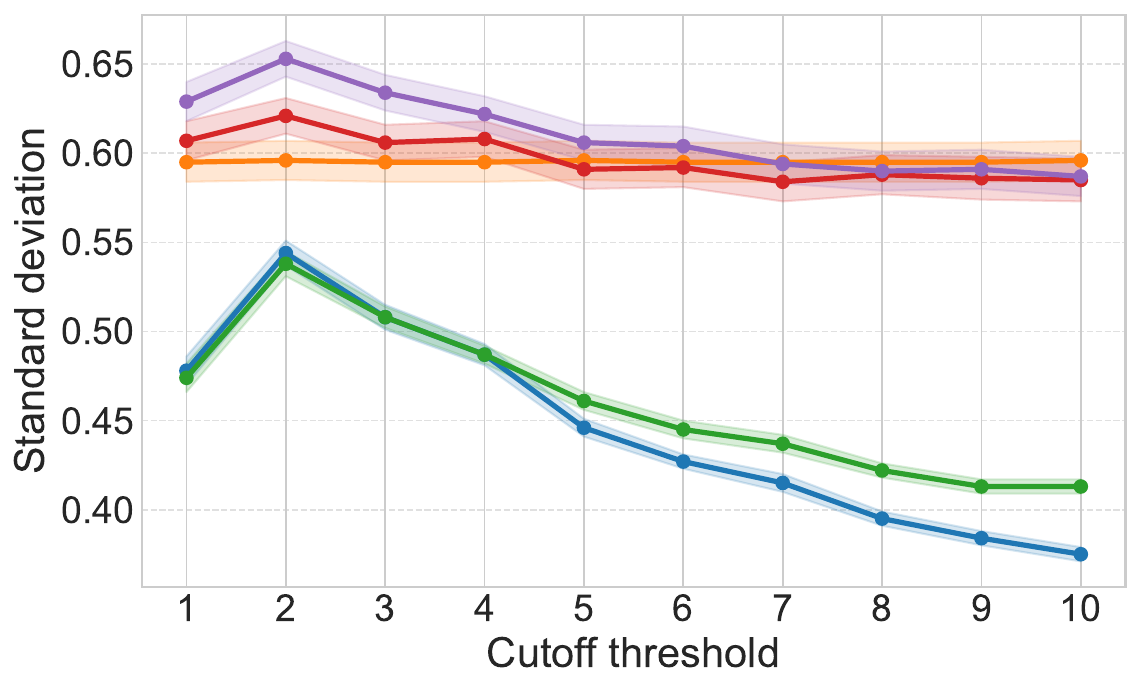}
         \caption{Return Standard Deviation}
     \end{subfigure}
     \begin{subfigure}[b]{0.3\textwidth}
         \centering
     \end{subfigure}
    \caption{Out-of-sample performance of optimal portfolios on synthetic data (mean (std.\ error) over 200~replications). For readability, the negative range is shown on a compressed linear scale in panels (a) and (b).}
    \label{fig:portfolio-results:synthetic_performance}
\end{figure}

For all portfolios under consideration, we compute the expected value, the standard deviation and the Sharpe ratio of the corresponding portfolio return under the target distribution~$\PP_1$. Figure~\ref{fig:portfolio-results:synthetic_performance} reports the means and standard errors of these performance metrics across 200~independent simulation runs. We observe that the multi-source DRO portfolios achieve the highest expected returns and Sharpe ratios. Optimizing only in view of the scarce target data results in noisy portfolios that are susceptible to the optimizer's curse and thus perform poorly. Similarly, optimizing only in view of the source data results in biased portfolios because of the distribution shift from~$\PP_2$ to~$\PP_1$. Indeed, under the source distribution~$\mathbb{P}_2$, it would be optimal to distribute the available capital evenly across the assets from~1 to~$C$. 
However, if $C = 1$, then this means only investing in one asset and missing out on diversification. If $C>2$, then this means investing a fraction $(C-2)/C$ of the available capital in assets that destroy wealth under the target distribution~$\PP_1$.
A notable exception is $C=2$ in which case $\PP_2=\PP_1$.
Optimizing in view of the pooled data does not significantly improve performance because the target samples are outnumbered by the source samples and have thus only a minor impact on the optimal portfolio. These shortcomings are mitigated by optimizing over a ball centered at a 1-Wasserstein barycenter between the empirical distributions~$\widehat\PP_1$ and~$\widehat\PP_2$. By Corollary~\ref{cor:barycenterK2}, however, if~$\widehat{\PP}_1$ and~$\widehat{\PP}_2$ are assigned the same weights, then either of them constitutes a 1-Wasserstein barycenter. In the experiment we choose the one that generates the smallest worst-case risk on the validation dataset. The portfolios obtained via multi-source DRO outperform all baselines in terms of Sharpe ratio and expected return uniformly across all distribution shifts~$C<10$ and match the performance of single-source DRO on the target distribution for $C=10$. However, their returns display high standard deviations.
This is neither surprising nor troubling. The optimal portfolio under complete knowledge of the target distribution~$\PP_1$ allocates the available capital evenly to assets~1 and~2. While this portfolio maximizes expected return, it suffers from poor diversification. It is easy to show that the lowest possible portfolio standard deviation is obtained by investing equal amounts of money in all assets (even those with negative expected returns).

\subsection{Backtest on Real Data}
In the third experiment the asset universe consists of $d=3$ regional sector indices (\textit{Agriculture and Food Chain}, \textit{Commodity Producers} and \textit{Infrastructure}) provided by MSCI, which follow the Global Industry Classification Standard (GICS). Each index includes large and mid cap companies from a particular geographical region (Europe, USA or Pacific) denominated in US dollars. We believe that the European and US economies are sufficiently similar so that the European index return distributions provide a strong prior for those of the US and vice versa. However, we believe that the index return distributions of the Pacific area are not informative for those of Europe and the US. Note that the three sector indices were chosen because they cover basic industries and thus exist in all  considered geographical regions. We then define~$\PP_1$ and~$\PP_2$ as the index return distributions corresponding to (any) two different regions. These distributions are unknown, but we can construct respective empirical distributions~$\widehat \PP_1$ and~$\widehat\PP_2$ from historical return data available from the MSCI online database (\url{https://www.msci.com/end-of-day-data-regional}). As in Section~\ref{sec:synthetic-data}, we compare the optimal portfolios of the multi-source DRO problem~\eqref{eq:portfolio-problem} against those of different single-source DRO baselines. We use the same parametrization of $\varepsilon_1$ and $\varepsilon_2$ as in Section~\ref{sec:sensitivity}, vary $\lambda$ from 0 to 1 in steps of 0.1, and select $m\in\{0.002,0.005,0.01,0.02\}$. The Wasserstein radius
in the single-source DRO problems is chosen from the set $\{a \cdot 10^b : a \in \{1,2,5\},\; b \in \{-3,-2,-1,0,1,2\}\}$.

The training, test and validation datasets cover monthly return data over the periods from January 2008 to December 2014, from January 2015 to December 2015, and from January 2016 to December 2021, respectively. In each simulation run, we randomly sample $N_1\in\{5,10,20,50\}$ return vectors from the training dataset of the target distribution~$\PP_1$ and $N_2= 5$ return vectors from the training dataset of the source distribution~$\PP_2$. These samples are used to construct~$\widehat\PP_1$ and~$\widehat\PP_2$, respectively. In addition, we sample 5~return vectors from the validation dataset of the target distribution, which are used to tune the hyperparameters. For each portfolio under consideration, we then compute the empirical Sharpe ratio on the test samples of the target distribution. Table~\ref{tab:msci-portfolio-results} reports the means and standard errors of these empirical Sharpe ratios across 10~independent simulation runs for different choices of source and target regions and for different choices of the target sample size~$N_1$. We observe that multi-source DRO outperforms single-source DRO when the source data originates from a similar market ({\em e.g.}, European data can help to improve portfolios in the US but {\em not} in the Pacific area). When the source and target distributions differ too much, it is better to use single-source DRO based solely on the target data.

\begin{table}\caption{Out-of-sample Sharpe ratios of optimal portfolios on real data (mean (std.\ error) over 10~replications)}
    \centering
    {\scriptsize
\begin{tabular}{llr|lllll}
\toprule
Source & Target & $N_1$ &        Single-s.\ DRO on &        Single-s.\ DRO on&        Single-s.\ DRO on&    Single-s.\ DRO on&          Multi-s.\ DRO \\
 &  &  & target distribution & source distribution & pooled data &  barycenter &           \\
\hline\hline
\multirow{8}{*}{Europe} & \multirow{4}{*}{Pacific} & 5  &  \hspace{0.5cm}$\boldsamewidth{0.059\ (0.004)}$ &  \hspace{0.5cm}0.054\ (0.006) &  \hspace{0.5cm}0.052\ (0.005) &  0.052\ (0.003) &  0.039\ (0.006) \\
    &         & 10 &   \hspace{0.5cm}$\boldsamewidth{0.047\ (0.005)}$ &  \hspace{0.5cm}0.04\phantom{0}\ (0.004) &    \hspace{0.5cm}0.043\ (0.0) &  0.036\ (0.004) &  0.035\ (0.006) \\
    &         & 20 &  \hspace{0.5cm}$\boldsamewidth{0.052\ (0.004)}$ &  \hspace{0.5cm}0.049\ (0.008) &  \hspace{0.5cm}0.049\ (0.006) &  0.042\ (0.006) &  0.035\ (0.006) \\
    &         & 50 &  \hspace{0.5cm}$\boldsamewidth{0.051\ (0.004)}$ &  \hspace{0.5cm}0.047\ (0.005) &  \hspace{0.5cm}0.047\ (0.003) &  0.045\ (0.004) &  0.028\ (0.004) \\
\cline{2-8}
    & \multirow{4}{*}{USA} & 5  &   \hspace{0.5cm}0.101\ (0.01) &  \hspace{0.5cm}0.107\ (0.006) &  \hspace{0.5cm}0.103\ (0.007) &  0.112\ (0.008) &  $\boldsamewidth{0.113\ (0.003)}$ \\
    &         & 10 &  \hspace{0.5cm}0.095\ (0.007) &  \hspace{0.5cm}0.103\ (0.007) &  \hspace{0.5cm}0.104\ (0.007) &  0.108\ (0.006) &  $\boldsamewidth{0.112\ (0.004)}$ \\
    &         & 20 &  \hspace{0.5cm}0.094\ (0.008) &  \hspace{0.5cm}0.092\ (0.009) &  \hspace{0.5cm}0.103\ (0.007) &  0.107\ (0.008) &  $\boldsamewidth{0.116\ (0.004)}$ \\
    &         & 50 &  \hspace{0.5cm}0.105\ (0.006) &  \hspace{0.5cm}0.095\ (0.008) &  \hspace{0.5cm}0.105\ (0.006) &  0.113\ (0.005) &  $\boldsamewidth{0.115\ (0.004)}$ \\
\hline
\multirow{8}{*}{Pacific} & \multirow{4}{*}{Europe} & 5  &  \hspace{0.5cm}0.111\ (0.009) &    \hspace{0.5cm}0.088\ (0.0) &    \hspace{0.5cm}0.088\ (0.0) &  0.096\ (0.008) &  $\boldsamewidth{0.123\ (0.014)}$ \\
    &         & 10 &  \hspace{0.5cm}0.112\ (0.012) &    \hspace{0.5cm}0.088\ (0.0) &    \hspace{0.5cm}0.088\ (0.0) &  0.111\ (0.012) &  $\boldsamewidth{0.126\ (0.014)}$ \\
    &         & 20 &   \hspace{0.5cm}0.107\ (0.01) &    \hspace{0.5cm}0.088\ (0.0) &    \hspace{0.5cm}0.088\ (0.0) &   0.106\ (0.01) &  $\boldsamewidth{0.111\ (0.012)}$ \\
    &         & 50 &  \hspace{0.5cm}0.113\ (0.01) &    \hspace{0.5cm}0.089\ (0.0) &    \hspace{0.5cm}0.088\ (0.0) &  0.101\ (0.008) &  $\boldsamewidth{0.145\ (0.011)}$ \\
\cline{2-8}
    & \multirow{4}{*}{USA} & 5  &  \hspace{0.5cm}0.109\ (0.006) &  \hspace{0.5cm}0.115\ (0.006) &   \hspace{0.5cm}0.115\ (0.006) &  0.113\ (0.006) &  $\boldsamewidth{0.116\ (0.006)}$ \\
    &         & 10 &  \hspace{0.5cm}0.088\ (0.006) &   \hspace{0.5cm}0.11\phantom{0}\ (0.007) &   \hspace{0.5cm}0.11\phantom{0}\ (0.007) &  0.109\ (0.008) & $\boldsamewidth{0.113\ (0.006)}$ \\
    &         & 20 &  \hspace{0.5cm}0.099\ (0.007) &  \hspace{0.5cm}0.107\ (0.007) &  \hspace{0.5cm}0.106\ (0.006) &  0.114\ (0.006) &  $\boldsamewidth{0.114\ (0.006)}$ \\
    &         & 50 &  \hspace{0.5cm}0.111\ (0.005) &  \hspace{0.5cm}$\boldsamewidth{0.119\ (0.004)}$ &  \hspace{0.5cm}$\boldsamewidth{0.119\ (0.004)}$ & $ \boldsamewidth{0.119\ (0.004)}$ & $ \boldsamewidth{0.119\ (0.004)}$ \\
\hline
\multirow{8}{*}{USA} & \multirow{4}{*}{Europe} & 5  &  \hspace{0.5cm}0.117\ (0.014) &  \hspace{0.5cm}0.091\ (0.001) &   \hspace{0.5cm}0.09\phantom{0}\ (0.001) &  0.116\ (0.014) &  $\boldsamewidth{0.144\ (0.013)}$ \\
    &         & 10 &  \hspace{0.5cm}0.117\ (0.014) &  \hspace{0.5cm}0.091\ (0.001) &  \hspace{0.5cm}0.091\ (0.001) &  0.117\ (0.014) &  $\boldsamewidth{0.151\ (0.014)} $\\
    &         & 20 &  \hspace{0.5cm}0.119\ (0.011) &   \hspace{0.5cm}0.09\phantom{0}\ (0.001) &  \hspace{0.5cm}0.095\ (0.005) &  0.108\ (0.009) &   $\boldsamewidth{0.15\phantom{0}\ (0.011)}$ \\
    &         & 50 &    \hspace{0.5cm}0.11\phantom{0}\ (0.01) &  \hspace{0.5cm}0.091\ (0.001) &  \hspace{0.5cm}0.091\ (0.002) &    0.11\phantom{0}\ (0.01) &  $\boldsamewidth{0.141\ (0.013)} $\\
\cline{2-8}
    & \multirow{4}{*}{Pacific} & 5  &  \hspace{0.5cm}0.046\ (0.007) &  \hspace{0.5cm}0.049\ (0.002) & \hspace{0.5cm}$\boldsamewidth{0.051\ (0.004)}$ &  0.042\ (0.004) &  0.046\ (0.006) \\
    &         & 10 &   \hspace{0.5cm}$\boldsamewidth{0.07\phantom{0}\ (0.004)}$ &   \hspace{0.5cm}0.06\phantom{0}\ (0.004) &  \hspace{0.5cm}0.058\ (0.004) &  0.058\ (0.005) &  0.046\ (0.006) \\
    &         & 20 &  \hspace{0.5cm}0.048\ (0.004) & \hspace{0.5cm}$\boldsamewidth{0.053\ (0.003)}$ &  \hspace{0.5cm}0.049\ (0.003) &  0.047\ (0.003) &  0.041\ (0.004) \\
    &         & 50 &  \hspace{0.5cm}0.047\ (0.003) & \hspace{0.5cm}$\boldsamewidth{0.052\ (0.003)}$ &  \hspace{0.5cm}0.048\ (0.003) &  0.047\ (0.003) &  0.046\ (0.002) \\
\bottomrule
\end{tabular}
}
    \label{tab:msci-portfolio-results}
\end{table}

\paragraph{Acknowledgments} This research was supported by the Swiss National Science Foundation under the NCCR Automation, grant agreement 51NF40\_180545, and by the Spanish Ministry of Science and Innovation (AEI/10.13039/501100011033) through project PID2023-148291NB-I00.  Finally, the authors thankfully acknowledge the computer resources, technical expertise, and assistance provided by the SCBI (Supercomputing
and Bioinformatics) center of the University of M\'alaga.

\bibliographystyle{informs2014}
\bibliography{References}

\newpage

\begin{APPENDICES}

\section{Stability of Multi-Source DRO}
\label{app:stability}
Here we establish conditions under which the optimal value as well as the minimizer of problem~\eqref{eq:original-problem} change continuously under $p$-Wasserstein perturbations of the estimated source distributions~$\widehat{\mathbb P}_k$, $k\in[K]$. All Wasserstein distances are defined relative to a metric~$d$ on~$\Xi$ satisfying the following mild assumption.

\begin{assumption}[Metric Spaces]
\label{ass:metric-spaces}
The metric~$d$ on~$\Xi$ is such that $(\Xi,d)$ constitutes a complete, separable, geodesic metric space. Also, $p \geq 1$ is a fixed exponent, and~$\mathcal{P}_p(\Xi)$ is equipped with the topology induced by the $p$-Wasserstein distance corresponding to the metric~$d$.
\end{assumption}

The properties of $(\Xi, d)$ put forth in Assumption~\ref{ass:metric-spaces} are satisfied, for example, when~$d$ is induced by a norm. Since $(\mathcal{P}_p(\Xi), W_p)$ is a metric space, subsets of~$\mathcal{P}_p(\Xi)$ are compact if and only if they are sequentially compact. While $p$-Wasserstein balls in $\mathcal{P}_p(\Xi)$ may fail to be compact, they are always weakly compact \cite[Theorem~1]{yue2021linearWasserstein}. If~$\Xi$ is compact, however, then convergence in $p$-Wasserstein distance is equivalent to weak convergence \cite[Theorem~6.9]{villani2009OT}. In this case, $p$-Wasserstein balls in~$\mathcal{P}_p(\Xi)$ are compact. 

The next lemma shows that the intersection of $K$~Wasserstein balls constitutes a continuous multifunction of the estimated source distributions $\widehat{\mathbb P}_k$, $k\in[K]$, at the centers of these balls.

\begin{lemma}[Continuity of Intersections of Wasserstein Balls]
\label{lem:continuity-multifunction}
Suppose that $\Xi$ is compact and that Assumption~\ref{ass:metric-spaces} holds. Define the multifunction
\[
	F : \mathcal P_p(\Xi)^K \rightrightarrows \mathcal P_p(\Xi) \quad \text{via} \quad F(\widehat{\mathbb P}_1,\dots,\widehat{\mathbb P}_K) = \bigcap_{k=1}^K \mathbb B_{\varepsilon_k}(\widehat{\mathbb P}_k),
\]
where $\mathbb B_{\varepsilon_k}(\widehat{\mathbb P}_k)$ denotes the $p$-Wasserstein ball of radius $\varepsilon_k>0$ around $\widehat{\mathbb P}_k$. Let $U \subseteq \mathcal P_p(\Xi)^K$  be the set on which $F(\widehat{\mathbb P}_1,\dots,\widehat{\mathbb P}_K)$ has nonempty interior. Then, $F$ is continuous on~$U$.
\end{lemma}

\proof{Proof of Lemma~\ref{lem:continuity-multifunction}}
For each $k\in[K]$, define $F_k : \mathcal P_p(\Xi) \rightrightarrows \mathcal P_p(\Xi)$ via $F_k(\widehat{\mathbb P}_k) = \mathbb B_{\varepsilon_k}(\widehat{\mathbb P}_k)$. We first prove that~$F_k$ is $1$-Lipschitz continuous in the Hausdorff distance $d_H$ on $\mathcal P_p(\Xi)$, that is,
\begin{align*}
	d_H\!\left(F_k(\widehat{\mathbb P}), F_k(\widehat{\mathbb P}')\right) 
	& \triangleq  \max\left\{ \sup_{{ \QQ}\in F_k(\widehat{\mathbb P})} \inf_{{ \QQ'}\in F_k(\widehat{\mathbb P}')} W_p({ \QQ}, { \QQ'}), \sup_{{ \QQ'}\in F_k(\widehat{\mathbb P}')} \inf_{{ \QQ}\in F_k (\widehat{\mathbb P})} W_p({ \QQ}, { \QQ'}) \right\} \\
	& \le W_p(\widehat{\mathbb P}, \widehat{\mathbb P}')  \quad \forall \,\widehat{\mathbb P},\widehat{\mathbb P}' \in \mathcal P_p(\Xi).
\end{align*}
To see this, select any ${ \QQ}\in \mathbb B_{\varepsilon_k}(\widehat{\mathbb P})$, and assume without loss of generality that ${ \QQ}\not\in \mathbb B_{\varepsilon_k}(\widehat{\mathbb P}')$. By Assumption~\ref{ass:metric-spaces}, there exists a geodesic from~${ \QQ}$ to~$\widehat{\mathbb P}'$ as well as a point~$\mathbb P^\star$ on this geodesic with $W_p({ \QQ}, \mathbb P^\star)= W_p({ \QQ}, \widehat{\mathbb P}')-\varepsilon_k$ and $W_p(\mathbb P^\star, \widehat{\mathbb P}')= \varepsilon_k$; see, {\em e.g.}, \citep[Corollary~7.22]{villani2009OT}. This implies that
\begin{align*}
	 \inf_{{ \QQ'}\in \mathbb B_{\varepsilon_k}(\widehat{\mathbb P}')} W_p({ \QQ}, { \QQ'}) = W_p({ \QQ}, \widehat{\mathbb P}')-\varepsilon_k \leq W_p(\widehat{\mathbb P}, \widehat{\mathbb P}'),
\end{align*}
where the equality holds because one can use the triangle inequality to prove that~$\mathbb P^\star$ attains the minimum. The inequality follows from the triangle inequality and the assumption that $W_p({ \QQ}, \widehat{\mathbb P})\leq\varepsilon_k$. Maximizing the left hand side of the above inequality over all ${ \QQ}\in \mathbb B_{\varepsilon_k}(\widehat{\mathbb P})$ yields
\[
	\sup_{{ \QQ}\in \mathbb B_{\varepsilon_k}(\widehat{\mathbb P})} \inf_{{ \QQ'}\in \mathbb B_{\varepsilon_k}(\widehat{\mathbb P}')} W_p({ \QQ}, { \QQ'}) \leq W_p(\widehat{\mathbb P},\widehat{ \mathbb P}').
\]
By symmetry, the same holds with~$\widehat{\mathbb P}$ and~$\widehat{\mathbb P}'$ exchanged. Therefore, the Hausdorff distance satisfies
\[
	d_H\!\left(F_k(\widehat{\mathbb P}), F_k(\widehat{\mathbb P}')\right) \le W_p(\widehat{\mathbb P}, \widehat{\mathbb P}')  \quad \forall \,\widehat{\mathbb P},\widehat{\mathbb P}' \in \mathcal P_p(\Xi).
\]
Hence, $F_k$ is indeed $1$-Lipschitz continuous in the Hausdorff distance on $\mathcal P_p(\Xi)$.

Next, we prove that the multifunction~$F$ is upper semicontinuous. To see this, fix $(\widehat{\mathbb P}_1,\dots,\widehat{\mathbb P}_K)\in U$ as well as an open set $V \subseteq \mathcal P_p(\Xi)$ with $F(\widehat{\mathbb P}_1,\dots,\widehat{\mathbb P}_K) \subseteq V$. As each~$F_k$ has compact values (owing to the compactness of~$\Xi$) and as~$F$ has nonempty values on~$U$, standard topological arguments ensure that there is~$\delta>0$ with
\[
	\bigcap_{k=1}^K \mathbb B_\delta(F_k(\widehat{\mathbb P}_k)) \subseteq V,
\]
where $\mathbb B_\delta(F_k(\widehat{\mathbb P}_k))$ denotes the $\delta$-neighborhood of~$F_k(\widehat{\mathbb P}_k)$. For any $\widehat{\mathbb Q}_k\in \mathbb B_\delta(\widehat{\mathbb P}_k)$ we further have
\[
	d_H\!\left(F_k(\widehat{\mathbb P}_k),F_k(\widehat{\mathbb Q}_k)\right) \leq W_p(\widehat{\mathbb P}_k,\widehat{\mathbb Q}_k) \leq \delta\quad \implies \quad F_k(\widehat{\mathbb Q}_k) \subseteq \mathbb B_\delta(F_k(\widehat{\mathbb P}_k))
\]
thanks to the Hausdorff-Lipschitz property of~$F_k$. In summary, these insights imply that
\[
	F(\widehat{\mathbb Q}_1,\dots,\widehat{\mathbb Q}_K) = \bigcap_{k=1}^K F_k(\widehat{\mathbb Q}_k) \subseteq \bigcap_{k=1}^K \mathbb B_\delta(F_k(\widehat{\mathbb P}_k)) \subseteq V \quad \forall (\widehat{\mathbb Q}_1,\dots,\widehat{\mathbb Q}_K) \in \bigtimes_{k=1}^K \mathbb B_\delta(\widehat{\mathbb P}_k).
\]
This confirms that if $F(\widehat{\mathbb P}_1,\dots,\widehat{\mathbb P}_K) \subseteq V$, then the values of~$F$ remain inside the open set~$V$ on a neighborhood of $(\widehat{\mathbb P}_1,\dots,\widehat{\mathbb P}_K)$. Hence, $F$ is indeed an upper semicontinuous multifunction.

Finally, we prove that~$F$ is lower semicontinuous. To this end, fix $(\widehat{\mathbb P}_1,\dots,\widehat{\mathbb P}_K)\in U$ as well as an open set $V \subseteq \mathcal P_p(\Xi)$ with $F(\widehat{\mathbb P}_1,\dots,\widehat{\mathbb P}_K) \cap V\neq \emptyset$. As the values of~$F$ have nonempty interior on~$U$ and as~$V$ is open, there exist~${ \QQ}\in V$ and~$\delta>0$ with $\mathbb B_\delta({ \QQ})\subseteq F(\widehat{\mathbb P}_1,\dots,\widehat{\mathbb P}_K)$.  Thus, $\mathbb B_\delta({ \QQ})\subseteq F_k(\widehat{\mathbb P}_k)$, and the Hausdorff-Lipschitz property of~$F_k$ implies that any $\widehat{\mathbb Q}_k\in \mathbb B_{\delta/2} (\widehat{\mathbb P}_k)$ satisfies
\[
	d_H\! \left(F_k(\widehat{\mathbb Q}_k), F_k(\widehat{\mathbb P}_k) \right) \leq \frac{\delta}{2} \quad \implies \quad { \QQ}\in F_k(\widehat{\mathbb Q}_k) \quad \forall k\in[K].
\]
In summary, we may therefore conclude that
\[
	{ \QQ}\in F(\widehat{\mathbb Q}_1,\dots,\widehat{\mathbb Q}_K) = \bigcap_{k=1}^K F_k(\widehat{\mathbb Q}_k) \quad \forall (\widehat{\mathbb Q}_1,\dots,\widehat{\mathbb Q}_K) \in \bigtimes_{k=1}^K \mathbb B_{\delta/2} (\widehat{\mathbb P}_k).
\]
This confirms that if $F(\widehat{\mathbb P}_1,\dots,\widehat{\mathbb P}_K) \cap V\neq\emptyset$, then the values of~$F$ have nonempty intersection with the open set~$V$ on a neighborhood of $(\widehat{\mathbb P}_1,\dots,\widehat{\mathbb P}_K)$. Hence, $F$ is indeed lower semicontinuous.

Since $F$ is both upper and lower semicontinuous and has compact values, it constitutes a continuous multifunction. This observation completes the proof. \QEDA
\endproof

The next lemma establishes sufficient conditions under which the expected loss is jointly continuous in the decisions and the probability distribution used to evaluate the expectation. 

\begin{lemma}[Continuity of Expected Loss]
\label{lem:continuity-loss}
If Assumption~\ref{ass:metric-spaces} holds, $\Xi$ is compact and~$\ell(\theta,\xi)$ is jointly continuous in~$\theta$ and~$\xi$, then the expected loss function $L: \Theta\times \mathcal P_p(\Xi) \to \mathbb R$ defined through $L(\theta,{ \QQ}) = \int_{\Xi} \ell(\theta,\xi)\,{\rm d}{ \QQ}(\xi)$ is jointly continuous in~$\theta$ and~${ \QQ}$.
\end{lemma}

\proof{Proof of Lemma~\ref{lem:continuity-loss}}
Let $\{\theta_t \}_{t\in\mathbb N}$ be a sequence in~$\Theta$ converging to~$\theta$, and let $\{{ \QQ_t}\}_{t\in\mathbb N}$ be a sequence in~$\mathcal P_p(\Xi)$ converging to~${ \QQ}$. The triangle inequality readily implies that
\begin{align}
	\label{eq:triangle-bound-L}
	\big|L(\theta_t,{ \QQ}_t)-L(\theta, { \QQ}) \big| \leq \big|L(\theta_t,{ \QQ_t})-L(\theta, { \QQ_t})\big| + \big|L(\theta,{ \QQ_t})-L(\theta, { \QQ})\big|
\end{align}
for every $t\in \mathbb N$. Recall now that~$\ell(\theta,\xi)$ is continuous in~$\xi$. Recall also that~$\Xi$ is compact, which implies that~$\ell(\theta,\xi)$ is bounded in~$\xi$ and that convergence in $\mathcal P_p(\Xi)$ is equivalent to weak convergence \cite[Theorem~6.9]{villani2009OT}. The second term on the right hand side of~\eqref{eq:triangle-bound-L} thus satisfies
\[
	\lim_{t\to\infty} \big|L(\theta,{ \QQ_t})-L(\theta, { \QQ})\big| = \lim_{t\to\infty} \left| \int_{\Xi} \ell(\theta,\xi)\,{\rm d}{ \QQ_t}(\xi) - \int_{\Xi} \ell(\theta,\xi)\,{\rm d}{ \QQ}(\xi)\right| = 0
\]
by the definition of weak convergence. As the continuous loss function~$\ell(\theta,\xi)$ is uniformly continuous on every compact subset of its domain~$\Theta\times\Xi$, we further have
\begin{equation*}
	\lim_{t\to\infty} \sup_{\xi\in\Xi} |\ell(\theta_t,\xi)-\ell(\theta,\xi)| = 0.
\end{equation*}
Consequently, the first term on the right hand side of~\eqref{eq:triangle-bound-L} satisfies
\begin{align*}
	\lim_{t\to\infty} \big|L(\theta_t,{ \QQ_t})-L(\theta, { \QQ_t})\big| \leq \lim_{t\to\infty} \int_{\Xi} |\ell(\theta_t,\xi)-\ell(\theta,\xi)|\,{\rm d}{ \QQ_t} (\xi) =0.
\end{align*}
Thus, the right hand side of~\eqref{eq:triangle-bound-L} tends to zero, and $L(\theta, { \QQ})$ is jointly continuous in~$\theta$ and~${ \QQ}$. \QEDA
\endproof

We can now prove that multi-source DRO is stable under perturbations of the source distributions. 

\begin{theorem}[Stability of Multi-Source DRO]
\label{thm:stability}
Suppose that Assumption~\ref{ass:metric-spaces} holds, $\Theta$ is compact, $\Xi$ is compact and~$\ell(\theta,\xi)$ is jointly continuous in~$\theta$ and~$\xi$. For each $k\in[K]$, let $\mathbb B_{\varepsilon_k}(\widehat{\mathbb P}_k)$ be the $p$-Wasserstein ball of radius $\varepsilon_k>0$ around $\widehat{\mathbb P}_k$, and define~$U\subseteq \mathcal P_p(\Xi)^K$ as the set of distributions~$\widehat{\mathbb P}_1,\ldots,\widehat{\mathbb P}_K$ for which the intersection of these Wasserstein balls has nonempty interior. Then, the optimal value of  problem~\eqref{eq:original-problem} is continuous in $\widehat{\mathbb P}_1,\dots,\widehat{\mathbb P}_K$ on~$U$.
\end{theorem}

\proof{Proof of Theorem~\ref{thm:stability}}
Define the multifunction~$F$ and the expected loss function~$L$ as in Lemmas~\ref{lem:continuity-multifunction} and~\ref{lem:continuity-loss}, respectively. Using this notation, problem~\eqref{eq:original-problem} can be represented compactly as $\inf_{\theta\in\Theta} \bar L(\theta;\widehat{\mathbb P}_1,\dots,\widehat{\mathbb P}_K)$, where the adversary's value function $\bar L: \Theta\times \mathcal P_p(\Xi)^K \to [-\infty,\infty]$ is given~by
\[
	\bar L(\theta;\widehat{\mathbb P}_1,\dots,\widehat{\mathbb P}_K) = \sup_{{ \QQ}\in F(\widehat{\mathbb P}_1,\dots,\widehat{\mathbb P}_K)} L(\theta, { \QQ}).
\]
Recall now from Lemmas~\ref{lem:continuity-multifunction} and~\ref{lem:continuity-loss} that~$F$ is continuous on~$U$ while~$L$ is continuous on~$\Theta\times\mathcal P_p(\Xi)$. Therefore, Berge's maximum theorem \cite[pp.~115--116]{berge1997topological} implies that nature's optimal value function~$\bar L$ is continuous on~$\Theta\times U$. As~$\Theta$ is compact and independent of the source distributions, a second application of Berge's maximum theorem finally shows that the optimal value of problem~\eqref{eq:original-problem} is continuous in $\widehat{\mathbb P}_1,\dots,\widehat{\mathbb P}_K$ on~$U$. Thus, the claim follows.
\QEDA
\endproof

\begin{corollary}
\label{cor:stability-optimizers}
If all conditions of Theorem~\ref{thm:stability} hold and problem~\eqref{eq:original-problem} has a unique minimizer $\theta^\star (\widehat{\mathbb P}_1,\dots,\widehat{\mathbb P}_K)$ for every $(\widehat{\mathbb P}_1,\dots,\widehat{\mathbb P}_K)\in U$, then $\theta^\star$ represents a continuous function on~$U$.
\end{corollary}

\proof{Proof of Corollary~\ref{cor:stability-optimizers}}
By Berge's maximum theorem \cite[pp.~115--116]{berge1997topological}, the argmin multifunction of problem~\eqref{eq:original-problem} is upper semicontinuous in $\widehat{\mathbb P}_1,\dots,\widehat{\mathbb P}_K$ on~$U$. By assumption, for any given $(\widehat{\mathbb P}_1,\dots,\widehat{\mathbb P}_K)\in U$ this multifunction returns the singleton $\{ \theta^\star (\widehat{\mathbb P}_1,\dots,\widehat{\mathbb P}_K) \}$. The definition of upper semicontinuity for multifunctions thus readily implies that $\theta^\star$ must be continuous.
\QEDA
\endproof

\vspace{3ex}

\section{Proofs}
\label{sec:proofs}
This appendix contains the proofs of all formal results in the main paper in order of appearance.

\proof{Proof of Lemma~\ref{lemma:barycenter:existence}}
If $\lambda_k=0$ for every $k\in[K]$, then $\BARY=\delta_{\xi_0}$ solves~\eqref{eq:barycenter}, and thus the claim trivially holds. From now on we may thus assume without loss of generality that $\lambda_1>0$. In this case, the objective function of~\eqref{eq:barycenter} is bounded below by $\lambda_1 C({ \QQ},\mathbb{P}_1)$. Assumption~\ref{ass:barycenter-existence} further implies that ${ \QQ}=\delta_{\xi_0}$ is feasible in~\eqref{eq:barycenter}, and thus $\Bar{C}=\sum_{k=1}^K\lambda_k C(\delta_{\xi_0},\mathbb{P}_k)<\infty$ provides an upper bound on the optimal value of~\eqref{eq:barycenter}. Taken together, these insights reveal that we can restrict the feasible set of problem~\eqref{eq:barycenter} to $B=\{{ \QQ}\in\mathcal{P}(\Xi): C({ \QQ},\mathbb{P}_1)\leq \Bar{C}/\lambda_1\}$ without affecting its optimal value. The set $B$ is weakly compact by~\citep[Proposition~2.5]{sorooshOTDRO}. Furthermore, $C({ \QQ}, \mathbb{P}_k)$ is weakly lower semicontinuous in ${ \QQ}$ for every $k\in[K]$~\cite[Lemma~5.2]{clement2008wasserstein}, and thus  the objective function of problem~\eqref{eq:barycenter} is weakly lower semicontinuous. Problem~\eqref{eq:barycenter} is thus solvable by Weierstrass' theorem.~\QEDA
\endproof

\proof{Proof of Theorem~\ref{thm:barycenter:pushforward-reformulation}}
As in Lemma~\ref{lemma:barycenter:existence}, we may assume without loss of generality that $\lambda_1>0$. By Assumption~\ref{ass:assumptions-c}, the transportation cost function~$c$ is lower semicontinuous and has bounded sublevel sets. The minimum of~\eqref{eq:barycenter:inner-min} is thus attained by Weierstrass' Theorem. In addition, the $\arg\min$ multifunction of problem~\eqref{eq:barycenter:inner-min} admits a measurable selector thanks to~\cite[Corollary~14.6 and Theorem~14.37]{rockafellar2009variational}. Thus, the functions~$\phi$ and~$\Phi$ are well-defined.

Assumptions~\ref{ass:assumptions-c} and~\ref{ass:barycenter-existence} imply via Lemma~\ref{lemma:barycenter:existence} that~\eqref{eq:barycenter} is solvable. In the following we will construct a sequence of optimization problems that are all equivalent to~\eqref{eq:barycenter}. As usual, we say that two optimization problems are equivalent if any feasible solution of the first problem can be used to construct a feasible solution to the second problem with the same (or a smaller) objective function value and vice versa. This correspondence holds in particular for all optimal solutions (if they exist). Thus, if one of the two equivalent problems is solvable, then so is the other one.

In a first step we note that problem~\eqref{eq:barycenter} is equivalent to
\begin{align}\label{eq:otbary-1}\tag{BC-1}
\begin{split}
    \min\quad & \sum_{k=1}^K\lambda_k \int_{\Xi^2} c(\xi,\xi_k) \,{\rm d}\pi_k(\xi,\xi_k)\\
    {\rm s.t.}\quad & { \QQ}\in\mathcal{P}(\Xi),\ \pi_k\in\Pi({ \QQ}, \mathbb{P}_k)\quad\forall k\in[K].
\end{split}
\end{align}
To see this, assume first that ${ \QQ}$ is feasible in~\eqref{eq:barycenter}. By \citep[Theorem~4.1]{villani2009OT}, there exists an optimal transportation plan $\pi_k\in\Pi({ \QQ},\mathbb P_k)$ with $C({ \QQ}, \mathbb{P}_k)=\int_{\Xi^2} c(\xi, \xi_k) {\rm d}\pi_k (\xi, \xi_k)$ for every $k=1,\ldots, K$. By the definition of the OT cost, it is then clear that $({ \QQ},\pi_1,\ldots,\pi_K)$ is feasible in~\eqref{eq:otbary-1} and attains the same objective function value as {$\QQ$} in~\eqref{eq:barycenter}. Conversely, if $({ \QQ},\pi_1,\ldots,\pi_K)$ is feasible in~\eqref{eq:otbary-1}, then {$\QQ$} is feasible in~\eqref{eq:barycenter} with the same objective function value.

In a second step we prove that problem~\eqref{eq:otbary-1} is equivalent to
\begin{align}\label{eq:otbary-2}\tag{BC-2}
\begin{split}
    \min\quad & \int_{\Xi^{K+1}} \sum_{k=1}^K\lambda_k c(\xi,\xi_k)\, {\rm d}\pibar(\xi,\xi_1,\ldots, \xi_K)\\
    {\rm s.t.}\quad & { \QQ}\in\mathcal{P}(\Xi),\ \pibar\in\Pi({ \QQ}, \mathbb{P}_1, \ldots,\mathbb{P}_K).
    \end{split}
\end{align}
To see this, assume first that $({\QQ},\pi_1,\ldots,\pi_K)$ is feasible in~\eqref{eq:otbary-1}. An iterative application of the gluing lemma~\citep[Chapter~1]{villani2009OT} implies that the $K$ couplings $\pi_k\in\Pi({ \QQ},\mathbb{P}_k), k\in[K]$, can be merged to a multi-margin transportation plan $\pibar \in\Pi({ \QQ}, \mathbb{P}_1, \ldots,\mathbb{P}_K)$ satisfying $P_{0,k\#}\pibar=\pi_k$, where the truncation operator $P_{0,k}:\Xi\times\Xi^k\rightarrow\Xi^2$ is defined through $P_{0,k}(\xi,\xi_1,\ldots,\xi_K)=(\xi,\xi_k)$ for every $k\in[K]$. By the measure-theoretic change-of-variables formula, we thus have 
\[
    \sum_{k=1}^K \lambda_k \int_{\Xi^2} c(\xi,\xi_k) \,{\rm d}\pi_k(\xi,\xi_k) = \int_{\Xi^{K+1}}\sum_{k=1}^K\lambda_k c(\xi,\xi_k) \, {\rm d}\pibar (\xi,\xi_1,\ldots, \xi_K).
\]
Hence, $({\QQ}, \pibar)$ is feasible in~\eqref{eq:otbary-2} and attains the same objective function value as $({ \QQ},\pi_1,\ldots,\pi_K)$. Conversely, if $({ \QQ}, \pibar)$ is feasible in~\eqref{eq:otbary-2}, then $({ \QQ},\pi_1,\ldots,\pi_K)$ with $\pi_k=P_{0,k\#}\pibar$ for every $k\in[K]$ is feasible in~\eqref{eq:otbary-1} and attains the same objective function value as $({ \QQ}, \pibar)$.

Next, we show that problem~\eqref{eq:otbary-2} is equivalent to
\begin{align}\label{eq:otbary-3}\tag{BC-3}
\begin{split}
    \min\quad &\int_{\Xi^K}\int_\Xi\sum_{k=1}^K\lambda_k c(\xi,\xi_k)\, {\rm d}\pi_{\xi|\xi_1,\ldots,\xi_K}(\xi|\xi_1, \ldots,\xi_K) \, {\rm d}\pi(\xi_1, \ldots,\xi_K)\\
    {\rm s.t.}\quad & \pi\in\Pi(\mathbb{P}_1, \ldots,\mathbb{P}_K),\ \pi_{\xi|\xi_1, \ldots,\xi_K}\in\mathcal{P}_{\xi|\xi_1,\ldots,\xi_K},
\end{split}
\end{align}
where~$\mathcal{P}_{\xi|\xi_1,\ldots,\xi_K}$ denotes the family of all transition kernels $\pi_{\xi|\xi_1, \ldots,\xi_K}$ satisfying the following two conditions. First, $\pi_{\xi|\xi_1,\ldots, \xi_K}(\cdot|\xi_1,\ldots,\xi_K)$ is a probability distribution on $\Xi$ for all $(\xi_1,\ldots,\xi_K)\in\Xi^K$. In addition, $\pi_{\xi|\xi_1,\ldots, \xi_K}(B|\cdot)$ is a Borel-measurable function on $\Xi^K$ for all Borel sets $B\subseteq\Xi$.

To prove the equivalence of~\eqref{eq:otbary-2} and~\eqref{eq:otbary-3}, assume first that $({ \QQ}, \pibar)$ is feasible in~\eqref{eq:otbary-2}. Next, define $\pi\in \Pi(\mathbb{P}_1, \ldots,\mathbb{P}_K)$ as the marginal distribution of $(\xi_1,\ldots, \xi_K)$ under~$\pibar$, and define $\pi_{\xi|\xi_1, \ldots,\xi_K}\in\mathcal{P}_{\xi|\xi_1,\ldots,\xi_K}$ as the regular conditional distribution of $\xi$ given $(\xi_1,\ldots, \xi_K)$ under~$\pibar$, which exists thanks to~\citep[Theorem~10.2.2]{dudley_2002}. Thus, $(\pi,\pi_{\xi|\xi_1, \ldots,\xi_K})$ is feasible in~\eqref{eq:otbary-3}. By the law of total probability, $({ \QQ},\pibar)$ and $(\pi, \pi_{\xi|\xi_1, \ldots,\xi_K})$ attain the same objective function values. Conversely, select any $(\pi, \pi_{\xi|\xi_1, \ldots,\xi_K})$ feasible in~\eqref{eq:otbary-3}, and define $\pibar\in\mathcal P(\Xi^{K+1})$ through
\[
    \pibar(B) = \int_{\Xi^K}\int_{\Xi} \mathds 1_{(\xi,\xi_1,\ldots,\xi_K)\in B}\, {\rm d}\pi_{\xi|\xi_1,\ldots,\xi_K}(\xi|\xi_1, \ldots,\xi_K) \, {\rm d}\pi(\xi_1, \ldots,\xi_K)
\]
for every Borel set $B\subseteq\Xi^{K+1}$. In addition, set ${ \QQ}=P_{0\#}\pibar$, where $P_0:\Xi\times\Xi^K\rightarrow\Xi$ is defined through $P_{0}(\xi,\xi_1,\ldots,\xi_K)=\xi$. We thus have $\pibar\in\Pi({ \QQ},\mathbb{P}_1,\ldots,\mathbb{P}_K)$, implying that $({ \QQ},\pibar)$ is feasible in~\eqref{eq:otbary-2}. One also easily verifies that $(\pi, \pi_{\xi|\xi_1, \ldots,\xi_K})$ and $({ \QQ},\pibar)$ attain the same objective function values.

We are now ready to prove that problem~\eqref{eq:otbary-3} is equivalent to
\begin{align}\label{eq:otbary-4}\tag{BC-4}
    \min_{\pi\in\Pi(\mathbb{P}_1, \ldots,\mathbb{P}_K)} ~\int_{\Xi^{K}}\min_{\pi_{\xi}\in\mathcal{P}(\Xi)}\int_{\Xi}\sum_{k=1}^K\lambda_k c(\xi, \xi_k) \, {\rm d}\pi_{\xi}(\xi) \, {\rm d}\pi(\xi_1,\ldots,\xi_K).
\end{align}
Note first that~\eqref{eq:otbary-4} is well defined. Indeed the inner minimization problem over the distribution~$\pi_\xi$ in~\eqref{eq:otbary-4} can be reformulated as
\begin{equation}\label{eq:phi-eqv}
\min_{\pi_{\xi}\in\mathcal{P}(\Xi)}\int_{\Xi}\sum_{k=1}^K\lambda_k c(\xi, \xi_k)\,{\rm d}\pi_{\xi}(\xi)= \min_{\xi_0\in\Xi}\sum_{k=1}^K\lambda_k c(\xi_0, \xi_k) = \phi(\xi_1,\ldots,\xi_K),
\end{equation}
where the first equality holds because the integral with respect to $\pi_\xi$ is minimized by the Dirac distribution $\delta_{\Phi(\xi_1,\ldots,\xi_K)}$ that places unit mass on the minimizer $\Phi(\xi_1,\ldots,\xi_K)$ of the integrand, while the second equality follows from the definition of~$\phi$. By Assumption~\ref{ass:barycenter-existence}, the optimal value of the inner minimization problem in~\eqref{eq:otbary-4} is thus integrable under every $\pi\in\Pi(\mathbb{P}_1,\ldots,\mathbb{P}_K)$.

The equivalence of~\eqref{eq:otbary-3} and~\eqref{eq:otbary-4} follows from a standard interchangeability principle; see, {\em e.g.}, \cite[Proposition~2.1]{shapiro2017interchangeability}. As the inner problem in~\eqref{eq:otbary-4} is solved by a Dirac distribution, however, it is easier to prove this equivalence directly instead of leveraging an abstract result. To this end, assume first that $(\pi, \pi_{\xi|\xi_1, \ldots,\xi_K})$ is feasible in~\eqref{eq:otbary-3}. Hence, $\pi$ is feasible in~\eqref{eq:otbary-4}. Note also that, for every fixed $(\xi_1,\ldots, \xi_K)$, the probability distribution $\pi_{\xi|\xi_1, \ldots,\xi_K}(\cdot|\xi_1,\ldots,\xi_K)$ is feasible (but generally suboptimal) in the inner problem in~\eqref{eq:otbary-4}. Thus, the objective function value of $\pi$ in~\eqref{eq:otbary-4} is smaller than or equal to that of $(\pi, \pi_{\xi|\xi_1, \ldots,\xi_K})$ in~\eqref{eq:otbary-3}. Conversely, if~$\pi$ is feasible in~\eqref{eq:otbary-4}, then $(\pi, \pi_{\xi|\xi_1, \ldots,\xi_K})$ is feasible in~\eqref{eq:otbary-3}, where $\pi_{\xi|\xi_1, \ldots,\xi_K}$ is defined through $\pi_{\xi|\xi_1\ldots,\xi_K}(\cdot|\xi_1,\ldots,\xi_K)=\delta_{\Phi(\xi_1,\ldots,\xi_K)}(\cdot)$. Note that $\pi_{\xi|\xi_1\ldots,\xi_K}$ defined in this way satisfies all properties of a transition kernel because~$\Phi$ is measurable. By construction, the objective function values of $(\pi, \pi_{\xi|\xi_1, \ldots,\xi_K})$ in problem~\eqref{eq:otbary-4} and of $\pi$ in problem~\eqref{eq:otbary-3} match.

Substituting~\eqref{eq:phi-eqv} into~\eqref{eq:otbary-4} shows that~\eqref{eq:otbary-4} is equivalent to~\eqref{eq:barycenter:outer-min}. Hence, the original OT barycenter problem~\eqref{eq:barycenter} is equivalent to~\eqref{eq:barycenter:outer-min}, as well, and~\eqref{eq:barycenter:outer-min} inherits solvability from~\eqref{eq:barycenter}. Define now $\BARY=\Phi_{\#}\pi\opt$, where $\pi\opt$ is any minimizer of the multi-margin OT problem~\eqref{eq:barycenter:outer-min}. 
It remains to be shown that $\mathbb{P}^\star$ solves~\eqref{eq:barycenter}. Indeed, for any Borel set $B\subseteq\Xi$ we have 
\begin{align*}
\BARY(B)&=\Phi_{\#}\pi^\star(B)=\pi^\star(\Phi^{-1}(B)) \\
&=\int_{\Xi^K}\delta_{\Phi(\xi_1,\ldots,\xi_K)}(B)\,{\rm d}\pi^\star(\xi_1,\ldots,\xi_K)\\
&=\int_{\Xi^K}\int_B {\rm d}\pi\opt_{\xi|\xi_1\ldots,\xi_K}(\xi|\xi_1,\ldots,\xi_K) \,{\rm d}\pi^\star(\xi_1,\ldots,\xi_K)\\
& = \int_{B\times \Xi^K} {\rm d}\pibar^\star(\xi,\xi_1,\ldots,\xi_K) = P_{0\#}\pibar^\star(B),
\end{align*}
where the transition kernel $\pi\opt_{\xi|\xi_1\ldots,\xi_K}$ is defined through $\pi\opt_{\xi|\xi_1\ldots,\xi_K}(\cdot|\xi_1,\ldots,\xi_K)=\delta_{\Phi(\xi_1,\ldots,\xi_K)}(\cdot)$, and the joint distribution $\pibar^\star$ is obtained by combining the marginal distribution $\pi\opt$ with the conditional distribution $\pi\opt_{\xi|\xi_1\ldots,\xi_K}$. As the Borel set $B\subseteq \Xi$ was chosen arbitrarily, we have thus shown that $\BARY=P_{0\#}\pibar^\star$. The equivalences between the different optimization problems established above lead to the following conclusion. As $\pi^\star$ solves~\eqref{eq:barycenter:outer-min} and~\eqref{eq:otbary-4} by assumption, $(\pi\opt, \pi\opt_{\xi|\xi_1\ldots,\xi_K})$ solves~\eqref{eq:otbary-3}, $(\BARY, \pibar^\star)$ solves~\eqref{eq:otbary-2}, and $\BARY$ solves~\eqref{eq:barycenter}. In addition, it is now easy to see that if $\Phi$ is unique, then every solution of~\eqref{eq:barycenter} must be of the form $\BARY=\Phi_{\#}\pi\opt$ for some solution $\pi^\star$ of~\eqref{eq:barycenter:outer-min}.
\QEDA
\endproof

\proof{Proof of Corollary~\ref{cor:barycenterK2}}
The claim follows directly from Theorem~\ref{thm:barycenter:pushforward-reformulation}. Indeed, if $\lambda_1\geq\lambda_2$, then $\Phi(\xi_1,\xi_2)=\xi_1$ is a measurable selector of the $\arg\min$ multifunction of problem~\eqref{eq:barycenter:inner-min}. Thus, if $\pi^\star$ solves problem~\eqref{eq:barycenter:outer-min}, then $\xi_{1\#}\pi^\star=\mathbb{P}_1$ solves the OT barycenter problem~\eqref{eq:barycenter} thanks to Theorem~\ref{thm:barycenter:pushforward-reformulation}.
\QEDA
\endproof

\proof{Proof of Proposition~\ref{prop:moment-reducing-barycenter}}
The proof consists of two parts. In the first part we show that the multi-margin transportation plan~$\pi^\star$ from Theorem~\ref{thm:barycenter:pushforward-reformulation} that induces~$\mathbb P^\star$ has a variance-maximizing property. In the second part we exploit this variance-maximizing property to prove that $\widehat{\mathbb P}$ is biased. 

\paragraph{Part~I} Assume without loss of generality that $\sum_{k=1}^K \lambda_k=1$, and let $\pi\in\Pi(\mathbb{P}_1,\ldots,\mathbb P_K)$ be any multi-margin transportation plan feasible in~\eqref{eq:barycenter:outer-min}. Next, construct a joint distribution $\mathbb{Q}_\pi$ of the random vector $(\xi,\xi_1,\ldots,\xi_K)\in\Xi^{K+1}$ under which the marginal distribution of $(\xi_1,\ldots,\xi_K)$ is given by~$\pi$, and the distribution of~$\xi$ conditional on $(\xi_1,\ldots,\xi_K)$ is given by $\sum_{k=1}^K\lambda_k\delta_{\xi_k}$. Note that the marginal distribution $\sum_{k=1}^K \lambda_k\mathbb{P}_k$ and the variance $\mathbb{Q}_\pi\VAR(\xi)$ of $\xi$ under $\mathbb Q_\pi$ are actually independent of~$\pi$. Using the same notation as in Theorem~\ref{thm:barycenter:pushforward-reformulation}, one readily verifies that
\[
    \mathbb{Q}_\pi\VAR(\xi|\xi_1,\ldots,\xi_K)=\min_{\xi_0\in\Xi}\sum_{k=1}^K \lambda_k \|\xi_0 -\xi_k\|_2^2= \phi(\xi_1,\ldots,\xi_K)
\]
and that the minimization problem over $\xi_0$ is uniquely solved by $\mathbb{E}_{\mathbb{Q}_\pi}[\xi|\xi_1,\ldots,\xi_K]=\Phi(\xi_1,\ldots,\xi_K)$. By the law of total variance, we thus have
\begin{equation}
\label{eq:total-variance}
\begin{aligned}
    \mathbb{Q}_\pi\VAR(\xi) & = \mathbb{E}_{\mathbb{Q}_\pi}[\mathbb{Q}_\pi\VAR(\xi|\xi_1,\ldots,\xi_K)] + \mathbb{Q}_\pi\VAR(\mathbb{E}_{\mathbb{Q}_\pi}[\xi|\xi_1,\ldots,\xi_K])\\
    &= \mathbb{E}_{\mathbb{Q}_\pi}[\phi(\xi_1,\ldots,\xi_K)] + \mathbb{Q}_\pi\VAR(\Phi(\xi_1,\ldots,\xi_K)) \\ & = \mathbb{E}_{\pi}[\phi(\xi_1,\ldots,\xi_K)] + \Phi_{\#}\pi\VAR(\xi),
\end{aligned}
\end{equation}
where the last equality holds because $\pi$ is the marginal of $(\xi_1,\ldots,\xi_K)$ under~$\mathbb Q_\pi$, which in turn implies that the pushforward distribution $\Phi_{\#}\pi$ is the marginal of $\Phi(\xi_1,\ldots, \xi_K)$ under~$\pi$. As the minimizer function $\Phi$ is unique, Theorem~\ref{thm:barycenter:pushforward-reformulation} implies that the OT barycenter~$\mathbb P^\star$ must coincide with $\Phi_{\#}\pi^\star$ for some solution~$\pi^\star$ of the multi-margin OT problem~\eqref{eq:barycenter:outer-min}, which minimizes $\mathbb{E}_{\pi}[\phi(\xi_1,\ldots,\xi_K)]$ across all $\pi\in\Pi(\mathbb P_1,\ldots,\mathbb P_K)$. Recall now that $\mathbb{Q}_\pi\VAR(\xi)$ is independent of~$\pi$. We may thus conclude from~\eqref{eq:total-variance} that $\pi^\star$ must maximize the variance $\Phi_{\#}\pi\VAR(\xi)$ across all $\pi\in\Pi(\mathbb P_1,\ldots,\mathbb P_K)$.

\paragraph{Part~II} We now use the insights of Part~I to prove that $\widehat{\mathbb P}$ is biased. To this end, let $\widehat\pi\in\Pi(\widehat{\mathbb P}_1,\ldots, \widehat{\mathbb P}_K)$ be a multi-margin transportation plan with $\widehat{\mathbb P}=\Phi_{\#}\widehat \pi$. Note that $\widehat \pi$ exists thanks to Theorem~\ref{thm:barycenter:pushforward-reformulation}. Note also that the discrete distributions $\widehat{\mathbb P}$ and $\widehat\pi$ constitute random objects because they depend on the random samples underlying the empirical distributions $\widehat {\mathbb P}_1,\ldots, \widehat{\mathbb P}_K$. By \cite[Theorem~14.37]{rockafellar2009variational}, we may assume without loss of generality that both the locations as well as the probabilities of the atoms of $\widehat{\mathbb P}$ and $\widehat\pi$ depend measurably on these samples.\footnote{Strictly speaking, this is an assumption about $\widehat{\mathbb P}$ and should be mentioned in the proposition statement. As it can be imposed without loss of generality, we hide it in the proof for better readability of the proposition statement.} Next, define the distribution $\pi=\mathbb{E}[\widehat{\pi}]$ through $\pi(B)=\mathbb{E}[\widehat{\pi}(B)]$ for every Borel set $B\subseteq \Xi^K$, which is well-defined because $\widehat{\pi}(B)$ is a nonnegative measurable function of the samples. It is easy to verify that $\pi\in\Pi(\mathbb P_1,\ldots,\mathbb P_K)$, that is, $\pi$ is feasible in~\eqref{eq:barycenter:outer-min}. In summary, we thus have
\begin{align*}
    \mathbb{E}[\widehat{\mathbb{P}}\VAR(\xi)]&= \mathbb{E}[\Phi_{\#}\widehat \pi\VAR(\xi)] = \mathbb{E}\left[\widehat{\pi}\VAR(\Phi(\xi_1,\ldots,\xi_K))\right] \\& \leq \mathbb{E}[\widehat{\pi}]\VAR(\Phi(\xi_1,\ldots,\xi_K))=\pi\VAR(\Phi(\xi_1,\ldots,\xi_K))= \Phi_{\#} \pi\VAR(\xi)\\
    &< \Phi_{\#} \pi^\star\VAR(\xi) = \mathbb{P}^\star\VAR(\xi).
\end{align*}
Here, the first inequality follows from Jensen's inequality, which applies because the variance of any random vector is concave in that random vector's distribution. The second inequality follows from the variance-maximizing property of~$\pi^\star$ established in Part~I. This inequality is strict if $\pi^\star\neq\pi$. Indeed, the multi-margin OT plan~$\pi^\star$ is the unique minimizer of $\mathbb{E}_{\pi}[\phi(\xi_1,\ldots,\xi_K)]$ (by \cite[Theorem~4.1]{agueh2011barycenters}) and therefore also the unique maximizer of $\Phi_{\#}\pi\VAR(\xi)$ (by Part~I).

It remains to be shown that $\pi^\star\neq\pi=\mathbb{E}[\widehat\pi]$. By \cite[Theorem~4.1]{agueh2011barycenters} and \citep[Theorem~2.1]{gangbo1998optimal}, there exists a bijective transformation $T:\Xi\rightarrow\Xi$ with $T_{\#}\mathbb{P}_1=\mathbb{P}_2$ and~$\pi^\star(\xi_2=T(\xi_1))=1$. Select now two Borel sets $A,B\subseteq \Xi$ with $A\cap B=\emptyset$ and $\mathbb{P}_1(A),\, \mathbb{P}_1(
B)>0$. Such sets exist because $\mathbb P_1$ is continuous. We thus have
\[
    \pi^\star(\xi_1\in A,\; \xi_2\in T(B))= \pi^\star(\xi_1\in A,\; \xi_1\in B) = 0,
\]
where the first equality holds because $\pi^\star(\xi_2=T(\xi_1))=1$ and $T$ is bijective, while the second equality holds because $A\cap B=\emptyset$. Next, observe that the event in which~$\widehat{\mathbb{P}}_1$ is supported on~$A$ and~$\widehat{\mathbb{P}}_2$ is supported on~$T(B)$ has a strictly positive probability equal to $\mathbb{P}_1(A)^{N_1}\mathbb{P}_2(T(B))^{N_2}=\mathbb{P}_1(A)^{N_1} \mathbb{P}_1(B)^{N_2}>0$. In this event, we have $\widehat{\pi}(\xi_1\in A,\;\xi_2\in T(B))=1$ because~$\widehat{\mathbb{P}}_1$ and~$\widehat{\mathbb{P}}_2$ are the marginal distributions of~$\xi_1$ and~$\xi_2$ under~$\widehat{\pi}$, respectively. This reasoning implies that
\begin{align*}
        \mathbb{E}[\widehat{\pi}](\xi_1\in A,\; \xi_2\in T(B)) &= \mathbb{E}[\widehat{\pi}(\xi_1\in A,\; \xi_2\in T(B))]\geq \mathbb{P}_1(A)^{N_1} \mathbb{P}_2(T(B))^{N_2}>0.
\end{align*}
We thus have $\mathbb{E}[\widehat{\pi}](\xi_1\in A,\; \xi_2\in T(B))>\pi^\star(\xi_1\in A,\; \xi_2\in T(B))$, which implies that $\pi^\star\neq\mathbb{E}[\widehat\pi]$. 
\QEDA
\endproof

The proof of Theorem~\ref{th_strong_duality} relies on a decomposability result for multi-margin transportation plans.

\begin{lemma}[Multi-Margin Transportation Plans]
\label{lemma:split_pi}
If Assumption~\ref{ass:discrete-reference-distributions} holds, then for any $\pibar\in\Pi({ \QQ}, \widehat{\mathbb{P}}_1,\ldots,\widehat{\mathbb{P}}_K)$ there exist finite Borel measures $\nu_\mi\in\mathcal{M}_+(\Xi)$, $\alpha\in\mathcal A$, with $\pibar=\sum_{\mi\in\MI}\nu_\mi\otimes\delta_\mi$, where $\delta_\mi$ is a shorthand for the Dirac distribution on~$\Xi^K$ that concentrates unit mass at $(\widehat{\xi}_{1,\mi_1},\ldots, \widehat{\xi}_{K,\mi_K})$.
\end{lemma}

\proof{Proof of Lemma~\ref{lemma:split_pi}}
Fix an arbitrary $\pibar \in \Pi({ \QQ}, \widehat{\mathbb{P}}_1,\ldots,\widehat{\mathbb{P}}_K)$. Next, for any $\alpha\in\mathcal A$, construct a nonnegative measure $\nu_\mi\in\mathcal{M}_+(\Xi)$ by setting $\nu_\mi(B)=\pibar(\xi\in B, \;\xi_k=\wh\xi_{k,\alpha_k}\ \forall k=1,\ldots, K)$ for all Borel sets $B\subseteq\Xi$. By construction, $\nu_\mi$ thus constitutes a rescaled probability measure. As any probability measure is a Radon measure, we may conclude that $\nu_\mi$ is indeed a Radon measure for every $\mi\in\MI$. It remains to be shown that the measures $\nu_\alpha$, $\alpha\in\mathcal A$, constructed in this way satisfy $\pibar=\sum_{\mi\in\MI}\nu_\mi\otimes\delta_\mi$. To this end, note that for any Borel sets $B, B_1,\ldots, B_K\subseteq \Xi$ we have
\begin{align*}
    \pibar\left(\xi\in B,\ \xi_k\in B_k~\forall k\right)
    =&\sum_{\substack{\mi\in\MI:\\
    \wh\xi_{k,\mi_k}\in B_k\forall k}} \pibar\left(\xi\in B,\ \xi_k=\wh\xi_{k,\alpha_k}\ 
    \forall k\right) \\
    =&\sum_{\substack{\mi\in\MI:\\
    \wh\xi_{k,\mi_k}\in B_k\forall k}} \nu_\alpha(\xi\in B) 
    =\sum_{\mi\in\MI} \nu_\mi(\xi \in B) \delta_\mi(\xi_k\in B_k\ \forall k),
\end{align*}
where the first equality follows from the law of total probability, while the second and the third equalities exploit the definitions of~$\nu_\alpha$ and~$\delta_\alpha$, respectively. The identity $\pibar=\sum_{\mi\in\MI}\nu_\mi\otimes\delta_\mi$ now follows because the sets of the form $B\times B_1\times \cdots\times B_K$ generate the Borel $\sigma$-algebra on~$\Xi^{K+1}$.
\QEDA
\endproof

\proof{Proof of Theorem~\ref{th_strong_duality}.} 
Recall that $C({ \QQ}, \widehat{\mathbb{P}}_k)$ denotes the optimal value of a minimization problem over transportation plans and is solvable thanks to Assumption~\ref{ass:assumptions-c}; see \citep[Theorem~4.1]{villani2009OT}. Consequently, the constraint $C({ \QQ}, \widehat{\mathbb{P}}_k)\leq \varepsilon_k$ is equivalent to the requirement that there exists $\pi_k\in\Pi({ \QQ}, \widehat{\mathbb{P}}_k)$ with $\int_{\Xi^2} c(\xi, \xi_k) \, {\rm d}\pi_k(\xi, \xi_k)\leq \varepsilon_k$. This in turn implies that problem~\eqref{eq:oringal-problem:inner-max} is equivalent to
\begin{equation}\label{eq:original-developed-wasserstein}
\begin{array}{c@{\quad}l@{\quad}l}
  \sup & \displaystyle \int_{\Xi}\ell(\xi)\,{\rm d}{ \QQ}(\xi)\\
    \rm{s.t.} & { \QQ}\in\mathcal{P}(\Xi),\ \pi_k\in\Pi({ \QQ}, \widehat{\mathbb{P}}_k) & \forall k\in[K]\\
    & \displaystyle \int_{\Xi^2} c(\xi, \xi_k) \, {\rm d}\pi_k(\xi, \xi_k)\leq \varepsilon_k &\forall k\in[K].
\end{array}   
\end{equation}
Next, fix any distribution ${ \QQ}\in\mathcal{P}(\Xi)$ and transportation plans $\pi_k\in\Pi({ \QQ},\widehat{\mathbb{P}}_k)$ for $k\in[K]$. An iterative application of the Gluing Lemma~\citep[Chapter~1]{villani2009OT} then shows that there must exist a multi-margin transportation plan $\pibar\in\Pi({ \QQ}, \widehat{\mathbb{P}}_1,\ldots, \widehat{\mathbb{P}}_K)$ with $P_{0,k\#}\pibar= \pi_k$ for all $k\in[K]$, where $P_{0,k}:\Xi^{K+1}\rightarrow \Xi^2$ is defined as usual through $P_{0,k}(\xi,\xi_1,\ldots,\xi_K) = (\xi,\xi_k)$. Appending $\pibar$ as a decision variable to problem~\eqref{eq:original-developed-wasserstein} and enforcing the constraints that link $\pibar$ to $\pi_k$, we obtain
\begin{equation*}
    \begin{array}{c@{\quad}l@{\quad}l}
    \sup & \displaystyle \int_\Xi \ell(\xi) \, {\rm d}{ \QQ}(\xi)\\
    \rm{s.t.} &{ \QQ}\in\mathcal{P}(\Xi),\ \pibar\in\Pi({ \QQ}, \widehat{\mathbb{P}}_1, \ldots, \widehat{\mathbb{P}}_K),\ \pi_k\in\Pi({\QQ}, \widehat{\mathbb{P}}_k) &\forall k\in[K]\\
    & P_{0,k\#}\pibar=\pi_k &\forall k\in[K]\\
    & \displaystyle \int_{\Xi\times\Xi}c(\xi,\xi_k)\,{\rm d}\pi_k(\xi,\xi_k)\leq\varepsilon_k &\forall k\in[K].
    \end{array}
\end{equation*}
By the measure-theoretic change-of-variables formula, problem~\eqref{eq:oringal-problem:inner-max} is thus equivalent to
\begin{equation}
\label{eq:formulation_nu}
\begin{array}{c@{\quad}l@{\quad}l}
    \sup & \displaystyle \int_\Xi \ell(\xi) \, {\rm d}{ \QQ} (\xi)\\
    \rm{s.t.} &{ \QQ}\in\mathcal{P}(\Xi),\ \pibar\in\Pi({ \QQ}, \widehat{\mathbb{P}}_1, \ldots, \widehat{\mathbb{P}}_K)\\
    & \displaystyle \int_{\Xi^{K+1}}c(\xi,\xi_k)\,{\rm d}\pibar(\xi,\xi_1, \ldots, \xi_K)\leq\varepsilon_k &\forall k\in[K].
\end{array}    
\end{equation}
By Lemma~\ref{lemma:split_pi}, every multi-margin transportation plan $\pibar\in\Pi({ \QQ}, \widehat{\mathbb{P}}_1,\ldots,\widehat{\mathbb{P}}_K)$ can be represented as $\pibar=\sum_{\mi\in\MI}\nu_\mi\otimes\delta_\mi$ for some finite Borel measures $\nu_\alpha\in \mathcal{M}_{+}(\Xi)$, $\alpha\in\mathcal A$, where $\delta_\mi$ is a shorthand for the Dirac point mass at $(\widehat{\xi}_{1,\mi_1},\ldots, \widehat{\xi}_{K,\mi_K})$. Thus, problem~\eqref{eq:formulation_nu} is equivalent to 
\begin{equation}
\label{eq:dro_primal}
\begin{array}{c@{\quad}l@{\quad}l}
  \sup & \displaystyle \sum_{ \mi \in \MI} \int_{\Xi}\ell(\xi)\,{\rm d}\nu_{\mi}(\xi)  \\ 
\rm{s.t.} & \nu_{\mi} \in \mathcal{M}_{+}(\Xi) &\forall  \mi \in \MI\\
    & \displaystyle \sum_{\substack{\mi \in \MI:\\\mi_k=j}}\int_{\Xi}{\rm d}\nu_{\mi}(\xi)=p_{k,j} &\forall j\in[N_k],\ \forall k\in[K]  \\
    & \displaystyle \sum_{\mi\in \MI}\int_\Xi c(\xi,\wh \xi_{k,\mi_k}) \, {\rm d}\nu_{\mi}(\xi)\leq \varepsilon_k &\forall k\in[K],
\end{array}
\end{equation} 
which can be viewed as a semi-infinite linear program over $|\mathcal A|=\prod_{k=1}^K N_k$ Borel measures.

Assigning Lagrange multipliers $\gamma_{k,j}\in\mathbb R$ and $\lambda_k\in\mathbb R_+$ to the probability matching and OT constraints, respectively, the Lagrangian dual of problem~\eqref{eq:dro_primal} can be represented as
\begin{equation}
\label{eq:dro_dual1}
\begin{aligned}
 \inf_{\lambda \geq 0,\, \gamma} \;\sup_{\{\nu_{\mi}\}_{\alpha}\subseteq\mathcal{M}_+(\Xi)} \;& \sum_{k =1}^K\lambda_k \varepsilon_k+\sum_{k=1}^K\sum_{j=1}^{N_k}\gamma_{k, j}\,p_{k,j}-\sum_{k=1}^K\lambda_k\left( \sum_{\mi \in \MI}\int_\Xi c(\xi,\wh \xi_{k,\mi_k})\,{\rm d}\nu_{\mi}(\xi) \right)  \\
    & - \sum_{k=1}^K\sum_{j=1}^{N_k}\gamma_{k,j}\sum_{\substack{\mi \in \MI:\\\mi_k=j}}\int_\Xi {\rm d}\nu_{\mi}(\xi)
    + \sum_{\mi \in \MI}\int_\Xi \ell(\xi) \,{\rm d}\nu_{\mi}(\xi).
\end{aligned}
\end{equation}
Routine rearrangements can then be used to show that problem~\eqref{eq:dro_dual1} is equivalent to
\begin{equation*}
  \inf_{\lambda \geq 0,\,\gamma}\;\sum_{k=1}^K \varepsilon_k \lambda_k +\sum_{k=1}^K\sum_{j=1}^{N_k} p_{k,j}\, \gamma_{k,j} +\sum_{\mi \in \MI}\sup_{\nu_{\mi}\in \mathcal{M}_{+}(\Xi)}\int_\Xi\left(\ell(\xi)-\sum_{k=1}^K\lambda_k c(\xi,\wh \xi_{k, \mi_k}) -\sum_{k=1}^K\gamma_{k,\mi_k} \right){\rm d}\nu_{\mi}(\xi).
\end{equation*}
Note that the supremum over $\nu_\alpha$ evaluates to infinity unless the following robust constraint holds.
$$
\ell(\xi)-\sum_{k=1}^K\lambda_k c(\xi,\wh \xi_{k,\mi_k})- \sum_{k=1}^K\gamma_{k,\mi_k} \leq 0\quad\forall \xi\in\Xi
$$
If this robust constraint holds, on the other hand, then the supremum over $\nu_\mi$ vanishes. Hence, the dual problem~\eqref{eq:dro_dual1} is indeed equivalent to~\eqref{eq:dro_dual}.

It remains to be shown that strong duality holds. Note first that, as Assumptions~\ref{ass:assumptions-c} and~\ref{ass:growth-condition} hold for some~$p'\leq p$, problem~\eqref{eq:oringal-problem:inner-max} and its equivalent reformulation~\eqref{eq:dro_primal} cannot be unbounded, and problem~\eqref{eq:dro_dual} is guaranteed to be feasible. In addition, the optimal values of~\eqref{eq:dro_primal} and~\eqref{eq:dro_dual} are ostensibly concave in~$\eps$. Using now standard perturbation arguments, one can show that these functions coincide and that strong duality holds for all $\eps\notin\partial\mathcal E$. Specifically, one can generalize the strong duality result in~\cite[Theorem~4.5]{kuhn2024distributionally} from a single Borel measure to finitely many Borel measures with obvious minor modifications. Details are omitted for brevity.
\QEDA
\endproof

\proof{Proof of Proposition~\ref{prop:sparse-sol}}
Assumptions~\ref{ass:assumptions-c}, \ref{ass:discrete-reference-distributions} and~\ref{ass:growth-condition} imply via a straightforward generalization of \citep[Theorem~3]{yue2021linearWasserstein} that~\eqref{eq:oringal-problem:inner-max} is solvable. In addition, from the proof of Theorem~\ref{th_strong_duality} we know that~\eqref{eq:oringal-problem:inner-max} is equivalent to problem~\eqref{eq:formulation_nu}. By using the definition of $\Pi({\QQ}, \widehat{\mathbb{P}}_1, \ldots, \widehat{\mathbb{P}}_K)$ and eliminating the decision variable~{$\QQ$}, one then easily verifies that problem~\eqref{eq:formulation_nu} is equivalent to
\begin{equation}
\label{eq:sparsity:eq0}
\begin{array}{c@{\quad}l@{\quad}l}
    \sup & \displaystyle \int_{\Xi^{K+1}} \ell(\xi)\, {\rm d}\pibar(\xi,\xi_1,\ldots, \xi_K)\\
    \rm{s.t.} &\pibar\in\mathcal{M}_+(\Xi^{K+1})\\
    & \displaystyle \int_{\Xi^{K+1}} {\rm d}\pibar(\xi,\xi_1,\ldots,\xi_K)=1\\
    &\displaystyle \int_{\Xi^{K+1}} \mathds{1}_{\xi_k=\wh\xi_{k,j}}\, {\rm d}\pibar(\xi, \xi_{1},\ldots, \xi_{K}) = p_{k,j} &\forall j \in [N_k],\ \forall k \in [K] \\
    & \displaystyle \int_{\Xi^{K+1}} c(\xi,\xi_k) \, {\rm d}\pibar(\xi,\xi_1,\ldots, \xi_K)\leq\varepsilon_k &\forall k \in [K].
\end{array}    
\end{equation}
Note that problem~\eqref{eq:sparsity:eq0} optimizes over the convex cone of nonnegative Borel measures on~$\Xi^{K+1}$. The normalization constraint ensures that each feasible Borel measure~$\pibar$ is in fact a probability distribution. As the probabilities~$\{p_{k,j}\}_{j=1}^{N_k}$ associated with the atoms of the $k$-th reference distribution~$\widehat{\mathbb P}_k$ sum to~1, the probability constraints ensure together with the normalization constraint that any feasible probability distribution~$\pibar$ is supported on $Z=\Xi\times\widehat \Xi_1\times\cdots\times\widehat \Xi_K$, where $\widehat \Xi_k=\{\wh\xi_{k,j}\}_{j=1}^{N_k}$ for every $k\in[K]$. We may thus restrict the integration domains of all integrals in~\eqref{eq:sparsity:eq0} to~$Z$ without restricting the problem's feasible set. As $\sum_{j=1}^{N_k}p_{k,j} =1$ for all $k\in[K]$, it then becomes evident that the normalization constraint makes one of the probability constraints redundant for every~$k\in[K]$. In summary, these arguments imply that problem~\eqref{eq:sparsity:eq0} can be reformulated as
\begin{equation}\label{eq:sparsity:eq1}
\begin{array}{c@{\quad}l@{\quad}l}
    \sup& \displaystyle \int_{Z} \ell(\xi) \, {\rm d}\pibar(\xi,\xi_1,\ldots, \xi_K)\\
    \rm{s.t.} &\pibar\in\mathcal{M}_+(Z)\\
    & \displaystyle \int_{Z} {\rm d}\pibar(\xi,\xi_1,\ldots, \xi_K)=1\\
    &\displaystyle \int_{Z}\mathds{1}_{\xi_k=\wh\xi_{k,j}}\,{\rm d}\pibar(\xi,\xi_1, \ldots, \xi_{K}) = p_{k,j} &\forall j \in [N_k-1], \ \forall k \in [K]\\
    & \displaystyle \int_{Z} c(\xi,\xi_k)\,{\rm d}\pibar(\xi,\xi_1,\ldots, \xi_K)\leq\varepsilon_k &\forall k \in [K].
\end{array}    
\end{equation}
Consider now any maximizer~$\pibar^\star$ of problem~\eqref{eq:sparsity:eq1}, which exists because~\eqref{eq:oringal-problem:inner-max} is solvable by assumption and because~\eqref{eq:sparsity:eq1} is equivalent to~\eqref{eq:oringal-problem:inner-max}. Next, define $\varepsilon^\star_k=\int_{Z}c(\xi,\xi_k)\,{\rm d}\pibar^\star(\xi,\xi_1,\ldots, \xi_K)$, and note that $\varepsilon^\star_k\leq \varepsilon_k$ for every~$k\in[K]$. By construction, problem~\eqref{eq:sparsity:eq1} is thus equivalent to
\begin{equation}\label{eq:sparsity:eq2}
\begin{array}{c@{\quad}l@{\quad}l}
    \sup& \displaystyle \int_{Z} \ell(\xi) \, {\rm d}\pibar(\xi,\xi_1,\ldots, \xi_K)\\
    \rm{s.t.} &\pibar\in\mathcal{M}_+(Z)\\
    & \displaystyle \int_{Z} {\rm d}\pibar(\xi,\xi_1,\ldots, \xi_K)=1\\
    &\displaystyle \int_{Z}\mathds{1}_{\xi_k=\wh\xi_{k,j}}\,{\rm d}\pibar(\xi,\xi_1, \ldots, \xi_{K}) = p_{k,j} &\forall j \in [N_k-1], \ \forall k \in [K]\\
    & \displaystyle \int_{Z} c(\xi,\xi_k)\,{\rm d}\pibar(\xi,\xi_1,\ldots, \xi_K)=\varepsilon^\star_k &\forall k \in [K],
\end{array}    
\end{equation}
which is also solved by~$\pibar^\star$. Problem~\eqref{eq:sparsity:eq2} can be viewed as an infinite-dimensional linear program in standard form. Indeed, the decision variable~$\pibar$ is subject to a continuum of nonnegativity constraints as well as to a finite number of linear equality constraints. Standard-form linear programs in finitely many dimensions are solved by basic feasible solutions. Intuitively, in the context of the infinite-dimensional linear program~\eqref{eq:sparsity:eq2}, a basic feasible solution is a discrete Borel measure whose number of atoms equals the number $N=1+\sum_{k=1}^K N_k$ of linear equality constraints. This intuition can be formalized. Specifically, by \cite[Corollary~5 and Proposition~6(v)]{pinelis2016extreme} and \citep[Proposition~1]{yue2021linearWasserstein},
problem~\eqref{eq:sparsity:eq2} is equivalent to
\begin{equation}\label{eq:sparsity:eq3}
\begin{array}{c@{\quad}l@{\quad}l}
    \sup& \displaystyle \int_{Z} \ell(\xi) \, {\rm d}\pibar(\xi,\xi_1,\ldots, \xi_K)\\
    \rm{s.t.} &\pibar\in\mathcal{D}_N(Z)\\
    & \displaystyle \int_{Z} {\rm d}\pibar(\xi,\xi_1,\ldots, \xi_K)=1\\
    &\displaystyle \int_{Z}\mathds{1}_{\xi_k=\wh\xi_{k,j}}\,{\rm d}\pibar(\xi,\xi_1, \ldots, \xi_{K}) = p_{k,j} &\forall j \in [N_k-1], \ \forall k \in [K]\\
    & \displaystyle \int_{Z} c(\xi,\xi_k)\,{\rm d}\pibar(\xi,\xi_1,\ldots, \xi_K)=\varepsilon^\star_k &\forall k \in [K],
\end{array}    
\end{equation}
where~$\mathcal D_N(Z)$ denotes the set of all discrete Borel measures in~$\mathcal M_+(Z)$ that are supported on at most~$N$ points. It remains to be shown that the supremum of problem~\eqref{eq:sparsity:eq3} is attained. By using similar techniques as in \cite[\S~5]{yue2021linearWasserstein}, however, it is easy to show that the objective function and the feasible set of problem~\eqref{eq:sparsity:eq3} are weakly upper semicontinuous and weakly compact, respectively. The solvability of problem~\eqref{eq:sparsity:eq3} thus follows from Weierstrass' theorem.

In summary, we have shown that problems~\eqref{eq:sparsity:eq0}--\eqref{eq:sparsity:eq3} all admit a minimizer~$\pibar^\star$ supported on at most $N$ points. Hence, the marginal distribution $\mathbb P^\star$ of $\xi$ under $\pibar^\star$ solves the uncertainty quantification problem~\eqref{eq:oringal-problem:inner-max} and is a discrete distribution supported on at most $N$ points.
\QEDA
\endproof

\proof{Proof of Theorem~\ref{thm:nphard-informal}}
Throughout this proof we focus on OT barycenter and uncertainty quantification problems with quadratic transportation cost functions $c(\xi,\xi')=\|\xi-\xi'\|_2^2$ and full support sets~$\Xi=\mathbb R^d$. Computing the optimal value of the OT barycenter problem~\eqref{eq:barycenter} with uniform weights $\lambda_1=\cdots=\lambda_K=1$ is known to be NP-hard even if all distributions $\mathbb P_1,\ldots, \mathbb P_K$ are discrete \citep[Theorem~1.1]{altschuler2022wasserstein}. We can thus prove NP-hardness of the \textsc{UQ Feasibility} problem by reduction from the OT barycenter problem. To this end, fix any instance of the OT barycenter problem specified by~$K$ discrete distributions $\widehat{\mathbb{P}}_1,\ldots,\widehat{\mathbb{P}}_K$ with at most~$N$ atoms and uniform weights. Recall now that~$\mathcal E$ is the set of all radii~$\varepsilon=(\varepsilon_1,\ldots,\varepsilon_K)$ for which the uncertainty quantification problem~\eqref{eq:oringal-problem:inner-max} with reference distributions $\widehat{\mathbb{P}}_1,\ldots,\widehat{\mathbb{P}}_K$ and vanishing loss function~$\ell(\xi)=0$ is feasible. As $\lambda_k=1$ for every $k\in[K]$, \eqref{eq:barycenter} is equivalent~to
\begin{equation}
    \label{eq:nphard-problem}
    \min_{\varepsilon\in \mathcal E}\, \sum_{k=1}^K\varepsilon_k.
\end{equation}
Indeed, if~$\mathbb P$ is feasible in~\eqref{eq:barycenter}, then~$\varepsilon$ defined through $\varepsilon_k=C(\mathbb P,\widehat{\mathbb P}_k)$, $k\in[K]$, is feasible in~\eqref{eq:nphard-problem} with the same objective value. Thus, the minimum of~\eqref{eq:nphard-problem} is smaller or equal to that of~\eqref{eq:barycenter}. Conversely, if~$\varepsilon$ is feasible in~\eqref{eq:nphard-problem}, then the definition of~$\mathcal E$ implies that there exists~$\mathbb P\in\mathcal P(\Xi)$ with $C(\mathbb P,\widehat{\mathbb P}_k)\leq\varepsilon_k$ for every~$k\in[K]$, and~$\mathbb P$ has the same or a lower objective value than~$\sum_{k=1}^K\varepsilon_k$. Thus, the minimum of~\eqref{eq:barycenter} is smaller or equal to that of~\eqref{eq:nphard-problem}. Next, note that~$\mathcal E$ is convex and closed because the OT distance is convex and weakly lower semicontinuous in its arguments. Note also that the \textsc{UQ Feasibility} problem with fixed discrete distributions $\widehat{\mathbb{P}}_1,\ldots,\widehat{\mathbb{P}}_K$ serves as a separation oracle for~\eqref{eq:nphard-problem}. Indeed, for any~$\bar\varepsilon\in\mathbb R^K$, the \textsc{UQ Feasibility} problem either confirms that~$\bar\varepsilon \in\mathcal E$ (which happens if and only if~\eqref{eq:oringal-problem:inner-max} is feasible) or enables us to construct a halfspace~$\mathcal H$ that contains~$\mathcal E$ but not~$\bar\varepsilon$. If $\bar\varepsilon\notin \mathcal E$ and~$\bar\varepsilon_k<0$ for some~$k\in[K]$, then we can simply set $\mathcal H=\{\varepsilon\in\mathbb R^K:\varepsilon_k\geq 0\}$. If $\bar\varepsilon\notin \mathcal E$ and~$\bar\varepsilon_k\geq 0$ for every~$k\in[K]$, on the other hand, then we can set
\[
    \mathcal H= \left\{\varepsilon\in\mathbb R^K: \sum_{k=1}^K\varepsilon_k \lambda^\infty_k+\sum_{k=1}^K\sum_{j=1}^{N_k}p_{k,j}\gamma^\infty_{k,j} \geq 0 \right\},
\]
where~$(\lambda^\infty, \gamma^\infty)$ is the recession direction of the feasible set of~\eqref{eq:dro_dual} returned by the \textsc{UQ Feasibility} oracle. As the objective function of~\eqref{eq:dro_dual} with~$\varepsilon=\bar\varepsilon$ decreases along~$(\lambda^\infty, \gamma^\infty)$, it is clear that~$\bar\varepsilon\notin\mathcal H$. To show that~$\mathcal{E}\subseteq\mathcal{H}$, assume for the sake of contradiction that there exists~$\Tilde \varepsilon\in\mathcal{E}$ with~$\Tilde \varepsilon\not\in\mathcal{H}$. Thus, $(\lambda^0,\gamma^0)+t\cdot(\lambda^\infty,\gamma^\infty)$ is feasible in~\eqref{eq:dro_dual} for any~$t\geq 0$ and for any feasible solution~$(\lambda^0, \gamma^0)$ of~\eqref{eq:dro_dual}, and we have
$$
    \lim_{t\rightarrow\infty} \sum_{k=1}^K\Tilde{\varepsilon}_k (\lambda^0_k+t\lambda^\infty_k)+\sum_{k=1}^K\sum_{j=1}^{N_k}p_{k,j}(\gamma^0_{k,j}+t\gamma^\infty_{k,j})=-\infty.
$$
This implies that~\eqref{eq:dro_dual} is unbounded, which in turn implies via weak duality that~\eqref{eq:dro_primal} is infeasible. However, this contradicts our assumption that~$\Tilde\varepsilon\in\mathcal{E}$. We may thus conclude that~$\mathcal E\subseteq\mathcal H$.

Next, we construct a ball that is guaranteed to contain all minimizers of~\eqref{eq:nphard-problem}. To this end, note that every distribution on~$\mathbb R^d$ is feasible in the OT barycenter problem~\eqref{eq:barycenter}. For example, $\widehat{\mathbb P}_1$ is feasible and has a finite objective function value $R= \sum_{k=1}^K C(\widehat{\mathbb P}_1, \widehat{\mathbb P}_k)$, which provides an upper bound on the optimal value of~\eqref{eq:barycenter}. As~\eqref{eq:barycenter} and~\eqref{eq:nphard-problem} are equivalent, this implies that all minimizers of~\eqref{eq:nphard-problem} must be contained in the simplex $\{\varepsilon\in\mathbb R^K_+: \sum_{k=1}^K \varepsilon_k\leq R\}$, which in turn resides within the ball of radius~$R$ around~$0$. Note that since $C(\widehat{\mathbb P}_1, \widehat{\mathbb P}_k)$ represents the optimal value of a linear program with $\poly(d,K)$ variables and constraints, $R$ can be encoded with $\poly(d,K, \log U)$ bits, where $\log U$ denotes the bit complexity of each parameter needed to describe the reference distributions ({\em i.e.}, the probabilities as well as the locations of their atoms). Trivially, $\mathcal E$ further contains balls of any positive radius because $\varepsilon\in\mathcal E$ implies that $\varepsilon'\in\mathcal E$ for every~$\varepsilon'\geq \varepsilon$.

Given the separation oracle for~$\mathcal E$ constructed above as well as the ball of radius~$R$ that contains all minimizers of~\eqref{eq:nphard-problem}, we can compute the optimal value of~\eqref{eq:nphard-problem} to any absolute accuracy~$\delta>0$ with the ellipsoid algorithm using $\poly(d,K, \log U, \log \frac{1}{\delta})$ calls to the \textsc{UQ Feasibility} oracle \citep[Theorem~5.2.1]{ben2001lectures}. Recall that the OT barycenter problem~\eqref{eq:barycenter} can be reformulated as a multi-margin OT problem, which is equivalent to a linear program with~$\mathcal O(N^K)$ variables and constraints and with objective function coefficients that can be computed via $\poly(d,K,\log U)$ arithmetic operations (see Theorem~\ref{thm:barycenter:pushforward-reformulation}). Hence, the (equal) optimal values of~\eqref{eq:barycenter} and~\eqref{eq:nphard-problem} can be encoded with $\poly(d,K,N, \log U)$ bits. Selecting the accuracy~$\delta$ of the ellipsoid algorithm to satisfy $\log(1/\delta) = \poly(d,K,N,\log U)$, we can thus compute the exact optimal value of~\eqref{eq:nphard-problem} by rounding the output of the ellipsoid algorithm to the nearest multiple of~$\delta$. In summary, we have thus conceived an algorithm for computing the optimal value of the NP-hard OT barycenter problem~\eqref{eq:barycenter}, which requires only $\poly(d,K, N, \log U)$ calls to the \textsc{UQ Feasibility} oracle. Consequently, the \textsc{UQ Feasibility} oracle must be NP-hard, too.
\QEDA\endproof

Recall that, for ease of exposition, we have assumed throughout Section~\ref{sec:tractability} that $N_k=N$ for all~$k \in [K]$. We have also assumed that $\log U$ represents an upper bound on the number of bits needed to encode any radius~$\varepsilon_k$, any probability~$p_{k,j}$ and any component of~$\widehat{\xi}_{k,j}$ for all $j\in[N]$ and~$k\in[K]$.

Before proving Theorem~\ref{thm:dual-complexity}, we first construct two different separation oracles for the feasible set~$\mathcal D$ of problem~\eqref{eq:dro_dual}. A separation oracle for~$\mathcal D$ decides whether a given point $(\lambda, \gamma)\in\mathbb R^{K + KN}$ belongs to~$\mathcal D$ and, if not, constructs a closed halfspace~$\mathcal H\subseteq \mathbb R^{K + KN}$ with $(\lambda,\gamma) \notin\mathcal H$ and $\mathcal D \subseteq \mathcal H$. The following lemma constructs the first of the two separation oracles. This oracle runs in polynomial time provided that~$K$ is fixed.

\begin{lemma}[Polynomial-Time Separation Oracle for Fixed~$K$]
\label{lemma:expK-sep-oracle}
If Assumptions~\ref{ass:compl-subproblem}, \ref{ass:loss-complexity} and~\ref{ass:compl-costf}($i$) hold, then there is a separation oracle for~$\mathcal D$ that runs in time polynomial in~$d$, $N$, $\log U$ and the bit size of~$(\lambda, \gamma)$.
\end{lemma}

\proof{Proof of Lemma~\ref{lemma:expK-sep-oracle}}
We construct a separation oracle for~$\mathcal D$ as follows. Given a test point~$(\lambda,\gamma)$, we first check if~$\lambda\geq 0$. Otherwise, if~$\lambda_k<0$ for some~$k\in[K]$, the closed halfspace $\mathcal H=\{(\lambda', \gamma')\in\mathbb R^{K+KN}:\lambda'_k\geq 0\}$ contains~$\mathcal D$ but not~$(\lambda,\gamma)$. For any fixed multi-index~$\mi\in\MI$ we then use the generalized Moreau envelope oracle from Assumption~\ref{ass:compl-subproblem} to compute the supremum on the left hand side of the $\mi$-th robust constraint in~\eqref{eq:dro_dual}, which we denote here by~$\ell^\star_\alpha(\lambda)$. If $\ell^\star_\alpha(\lambda)=\infty$, then the oracle outputs a closed halfspace~$\Lambda \subseteq\mathbb R^K$ such that~$\lambda\notin \Lambda$ and $\lambda'\in \Lambda$ for every~$\lambda'\in\mathbb R^K$ with $\ell^\star(\lambda')<\infty$. Thus, the closed halfspace~$\mathcal H =\Lambda\times \mathbb R^{KN}$ contains~$\mathcal{D}$ but not~$(\lambda,\gamma)$. If~$\sum_{k=1}^K\gamma_{k,\mi_k}< \ell^\star_\alpha(\lambda) <\infty$, on the other hand, then the oracle outputs a maximizer~$\xi^\star\in\Xi$ for the optimization problem embedded in the $\mi$-th robust constraint in~\eqref{eq:dro_dual}. The closed halfspace
\[
    \mathcal H= \left \{(\lambda',\gamma')\in\mathbb R^{K+KN} \;:\; \ell(\xi^\star) - \sum_{k=1}^K \lambda_k'c(\xi^\star,\widehat\xi_{k,\mi_k}) \leq \sum_{k=1}^K \gamma'_{k,\mi_k} \right\}
\]
then contains every feasible solution~$(\lambda',\gamma')\in \mathcal{D}$ because $\ell(\xi^\star) - \sum_{k=1}^K \lambda_k'c(\xi^\star,\widehat\xi_{k,\mi_k}) \leq \ell^\star_\mi(\lambda')\leq \sum_{k=1}^K \gamma'_{k,\mi_k}$, where the first inequality holds because~$\xi^\star$ is feasible but not necessarily optimal in the embedded optimization problem associated with~$\lambda'$, and the second inequality follows from the definition of the feasible set~$\mathcal D$. On the other hand, $(\lambda,\gamma)\notin\mathcal H$ because $\ell(\xi^\star) - \sum_{k=1}^K \lambda_k c(\xi^\star,\widehat\xi_{k,\mi_k}) = \ell^\star_\alpha(\lambda)> \sum_{k=1}^K\gamma_{k,\mi_k}$, where the equality holds because~$\xi^\star$ is optimal in the embedded optimization problem associated with~$\lambda$. If~$\ell^\star_\alpha(\lambda) \leq \sum_{k=1}^K\gamma_{k,\mi_k}$, finally, then we have verified that the $\mi$-th robust constraint in~\eqref{eq:dro_dual} is satisfied. If $\lambda\geq 0$ and all $|\MI|=N^K$ robust constraints are satisfied, then we can conclude that~$(\lambda,\gamma)\in\mathcal D$. The procedure outlined above thus correctly decides whether~$(\lambda,\gamma)\in\mathcal D$ and, if not, outputs a halfspace~$\mathcal H$ that contains~$\mathcal D$ but not~$(\lambda, \gamma)$.

The runtime of the proposed procedure is dominated by the~$N^K$ calls to the Moreau envelope oracle. By Assumptions~\ref{ass:compl-subproblem}, \ref{ass:loss-complexity} and~\ref{ass:compl-costf}($i$), the runtime of each call grows polynomially with~$d$, $K$, $N$, $\log U$ and the bit size of~$(\lambda,\gamma)$. Thus, if~$K$ is fixed, the oracle runs in time polynomial in~$d$, $N$, $\log U$ and the bit size of~$(\lambda,\gamma)$.
\QEDA
\endproof

Perhaps surprisingly, one can construct a second oracle that runs in time polynomial in~$K$ if~$d$ is fixed. This construction relies on the following lemma inspired by \citep[Lemma~17]{altschuler2021barycenterComplexity}.

\begin{lemma}[Constraint Elimination]
\label{lemma:complexity-prepa:restriction}
If Assumptions~\ref{ass:assumptions-c} and~\ref{ass:compl-costf}(ii) hold and $\lambda_k\geq 0$ for all $k\in[K]$, then
$$
    \max_{\mi\in\MI}\; \sup_{\xi\in\Xi} \; \ell(\xi)-\sum_{k=1}^K \left(\lambda_k c(\xi,\wh \xi_{k, \mi_k})+\gamma_{k,\mi_k} \right)
    ~=~ \max_{\mi\in \MI'}\;\sup_{\xi\in\Xi} \; \ell(\xi)-\sum_{k=1}^K \left(\lambda_k c(\xi,\wh \xi_{k, \mi_k})+\gamma_{k,\mi_k} \right),
$$
where $\MI'=\{\mi\in\MI:\mathcal C_\mi\neq\emptyset\}$ is defined in terms of the cells $\mathcal C_\mi=\cap_{k=1}^K \Xi_{k,\mi_k}$, $\mi\in\MI$, with 
\[
        \Xi_{k,i}=\left\{\xi\in\Xi \;:\;\lambda_kc(\xi, \xi_{k,i})+\gamma_{k,i}<\lambda_k c(\xi, \xi_{k,j})+\gamma_{k,j} ~~ \forall j\in[N],\; j\neq i \right\}\quad \forall k\in[K],~\forall i\in[N].
\]
\end{lemma}
\proof{Proof of Lemma~\ref{lemma:complexity-prepa:restriction}}
Since $\MI'\subseteq \MI$, the maximum over~$\MI'$ is smaller than or equal to the maximum over~$\MI$. To prove the reverse inequality, note that the maximum over~$\MI$ can be written~as
\begin{align*}
    \sup_{\xi\in \Xi} \; \max_{\mi\in\MI} ~\ell(\xi)-\sum_{k=1}^K \left(\lambda_k c(\xi, \wh \xi_{k,\mi_k}) + \gamma_{k,\mi_k}\right)
    =&\max_{\mi'\in \MI'} \; \sup_{\xi\in \mathcal C_{\mi'}} \; \max_{\mi\in\MI} ~ \ell(\xi)-\sum_{k=1}^K \left( \lambda_k c(\xi, \wh \xi_{k,\mi_k}) + \gamma_{k,\mi_k} \right) \\
    =&\max_{\mi'\in \MI'} \; \sup_{\xi\in \mathcal C_{\mi'}} ~ \ell(\xi) -\sum_{k=1}^K \left( \lambda_k c(\xi, \wh \xi_{k,\mi'_k}) + \gamma_{k,\mi'_k}\right)\\
    \leq&\max_{\mi'\in \MI'} \; \sup_{\xi\in \Xi} ~\ell(\xi) -\sum_{k=1}^K \left(\lambda_k c(\xi, \wh \xi_{k,\mi'_k}) + \gamma_{k,\mi'_k}\right),
\end{align*}
where the first equality holds because $\cup_{\mi'\in\MI'}\text{cl}(\mathcal C_{\mi'}) = \Xi$ thanks to Assumption~\ref{ass:compl-costf}($ii$) and because~$\ell$ and~$c$ are upper and lower semicontinuous, respectively (see Assumption~\ref{ass:assumptions-c}). The second equality exploits the construction of the cells~$\mathcal C_{\mi'}$, $\mi'\in\MI'$, and the inequality holds because~$\mathcal C_{\mi'}\subseteq \Xi$ for every~$\mi'\in\MI'$.
\QEDA
\endproof

Lemma~\ref{lemma:complexity-prepa:restriction} implies that, among all~$N^K$ robust constraints in~\eqref{eq:dro_dual}, it suffices to retain only those indexed by~$\mi\in\MI'$. The others can be eliminated without affecting the solution of~\eqref{eq:dro_dual}. By Assumption~\ref{ass:compl-costf}($iii$), the reduced index set can be computed in time~${\rm poly}(N, K, \log U)$ for any fixed dimension~$d$. This implies in particular that $|\MI'| = {\rm poly}(N, K, \log U)$, which is much smaller than $|\MI|= N^K$ for large~$K$. The following lemma constructs the second separation oracle, which runs in polynomial time provided that~$d$ is fixed.

\begin{lemma}[Polynomial-Time Separation Oracle for Fixed~$d$]
\label{lemma:expd-sep-oracle}
If Assumptions~\ref{ass:compl-subproblem}, \ref{ass:loss-complexity} and~\ref{ass:compl-costf} hold, then there is a separation oracle for~$\mathcal D$ that runs in time polynomial in~$K$, $N$, $\log U$ and the bit size of~$(\lambda, \gamma)$.
\end{lemma}

\proof{Proof of Lemma~\ref{lemma:expd-sep-oracle}}
The desired separation oracle is constructed as in the proof of Lemma~\ref{lemma:expK-sep-oracle}. However, the Moreau envelope oracle from Assumption~\ref{ass:compl-subproblem} is only called for robust constraints indexed by some~$\mi\in \MI'$. As $|\MI'| = {\rm poly}(N, K, \log U)$ and because of the assumed efficiency of the Moreau envelope oracle, the time needed for all $|\MI'|$ oracle calls grows polynomially with~$d$, $K$, $N$, $\log U$ and the bit size of~$(\lambda,\gamma)$. Also, by Assumption~\ref{ass:compl-costf}($iii$), the reduced index set~$\MI'$ can be computed in time $\poly(K,N,\log U)$. The overall runtime of the constructed separation oracle thus grows polynomially with~$K$, $N$, $\log U$ and the bit size of~$(\lambda,\gamma)$. 
\QEDA
\endproof

All known algorithms for computing the reduced index set~$\MI'$ have a runtime that grows either exponentially with~$d$ or with~$K$; see the discussion in \citep{altschuler2021barycenterComplexity}.

The proof of Theorem~\ref{thm:dual-complexity} follows now almost immediately from Lemmas~\ref{lemma:expK-sep-oracle} and~\ref{lemma:expd-sep-oracle}.

\proof{Proof of Theorem~\ref{thm:dual-complexity}}
The assumptions of the theorem imply via Lemmas~\ref{lemma:expK-sep-oracle} and \ref{lemma:expd-sep-oracle} that the feasible set of~\eqref{eq:dro_dual} admits a separation oracle that runs in time $\poly(d,N, \log U)$ and $\poly(K, N, \log U)$, respectively. Thus, problem~\eqref{eq:dro_dual} can be solved efficiently with the ellipsoid algorithm \citep[Theorem~5.2.1]{ben2001lectures}. The initial ellipsoid may be set to the ball of radius~$R$ around~$0$, which can be computed in time $\poly(d,K, N, \log U)$. As the transportation cost function is nonnegative, it is clear that if~$(\lambda,\gamma)$ is feasible in~\eqref{eq:dro_dual}, then any~$(\lambda',\gamma')\geq (\lambda,\gamma)$ is also feasible in~\eqref{eq:dro_dual}. Hence, we may assume without loss of generality that the intersection of the initial ellipsoid and the feasible set of problem~\eqref{eq:dro_dual} contains some ball of radius~$1$. The ellipsoid algorithm with the first (second) separation oracle thus computes the optimal value of~\eqref{eq:dro_dual} to within any accuracy~$\delta$ in time $\poly(d,N, \log U,\log \frac{1}{\delta})$ ($\poly(K, N, \log U,\log \frac{1}{\delta})$).
\QEDA
\endproof

\proof{Proof of Proposition~\ref{prop:conc_ineq_balls}}
By the definition of the $p$-Wasserstein balls, we have
\begin{align*}
    \PP_{[K]} \left( \mathbb{P} \in \bigcap_{k=1}^{K} \mathbb{B}^p_{\eps_k}(\widehat{\mathbb{P}}_k) \right) & = \PP_{[K]} \left( W_p(\PP, \widehat \PP_k)\leq \eps_k~\forall k\in[K]  \right) \\
    & \geq 1 - \sum_{k=1}^K \PP^{N_k}_k\left( W_p(\PP, \widehat \PP_k) > \eps_k\right) \\
    & \geq 1 - \sum_{k=1}^K \PP^{N_k}_k\left( W_p(\PP,\PP_k)+W_p(\PP_k, \widehat \PP_k) > \eps_k\right) \geq 1- \sum_{k=1}^K \beta_k,
\end{align*} 
where the first inequality exploits De~Morgan's laws and the union bound, whereas the second inequality follows from the triangle inequality. The third inequality holds thanks to
Lemma~\ref{lem:concentration-empirical-FG}, which applies because $\eps_k-W_p(\PP,\PP_k)\geq \eps_k-r_k\geq \eps(\beta_k,N_k)$ for all~$k\in[K]$. Thus, the claim follows.
\QEDA
\endproof

\proof{Proof of Proposition~\ref{lemma:bayesian-fixed_N}}
Bayes' rule implies that
\begin{align}
\label{eq:bayes-rule}
    \begin{split}
    & \mu\left(\left. W_p(\mathbb\PP_1, \widehat\PP_k)> \eps_k \right| W_p( \widehat\PP_1,\widehat\PP_k)\leq\widehat r_k \right) \\ &\hspace{3cm} = 
    \frac{\mu\left(\left. W_p( \widehat\PP_1,\widehat\PP_k)\leq\widehat r_k \right| W_p(\mathbb\PP_1, \widehat\PP_k)> \eps_k\right) \cdot \mu \left(W_p( \PP_1, \widehat\PP_k)> \eps_k\right)
    }{ \mu \left(W_p( \widehat\PP_1,\widehat\PP_k)\leq\widehat r_k\right) }.
    \end{split}
\end{align}
The evidence in the denominator of the above expression is strictly positive for otherwise it would have been ($\mu$-almost surely) impossible to observe~$\widehat r_k$. The likelihood term in the numerator satisfies
\begin{align*}
    \mu\left(\left. W_p( \widehat\PP_1,\widehat\PP_k)\leq\widehat r_k \right| W_p(\mathbb\PP_1, \widehat\PP_k)> \eps_k\right) & \leq \mu\left(\left. W_p(\PP_1, \widehat\PP_k) - W_p(\PP_1,\widehat\PP_1) \leq\widehat r_k \right| W_p(\mathbb\PP_1, \widehat\PP_k)> \eps_k\right)\\
    & \leq \mu\left(\left. W_p(  \PP_1,\widehat\PP_1) > \eps_k - \widehat r_k \right| W_p(\mathbb\PP_1, \widehat\PP_k)> \eps_k\right)\\
    & = \mathbb E_\mu\left[ \left. \mu\left( \left. W_p(  \PP_1,\widehat\PP_1) > \eps_k - \widehat r_k \right| \PP_1,\ldots,\PP_K \right)  \right| W_p(\mathbb\PP_1, \widehat\PP_k)> \eps_k \right]\\
     & \leq \beta(\eps_k - \widehat r_k,N_1),
\end{align*}
where the equality follows from the law of total expectation and the observation that~$\widehat \PP_1$ and~$\widehat \PP_k$ are independent under~$\mu$ conditional on any fixed realizations of the reference distributions~$\PP_{k'}$, $k'\in[K]$. The last inequality follows from Lemma~\ref{lem:concentration-empirical-FG}, which applies because, conditional on any fixed realizations of the reference distributions, $\widehat \PP_1$ is governed by the distribution~$\PP_1^{N_1}$ under~$\mu$ and because the reference distributions satisfy Assumption~\ref{light_tailed_assumption} $\mu$-almost surely. Similarly, the prior probability in the numerator of~\eqref{eq:bayes-rule} can be recast as
\begin{align*}
    \mu \left(W_p( \PP_1, \widehat\PP_k)> \eps_k\right) &\leq \mu \left( W_p( \PP_1, \PP_k) + W_p( \PP_k, \widehat\PP_k) > \eps_k\right)\\ 
    &= \int_0^\infty \mu \left( \left. W_p( \PP_1, \PP_k) + W_p( \PP_k, \widehat\PP_k) > \eps_k\right| W_p(\PP_1,\PP_k) = r_k \right) {\rm d}F_k(r_k)\\
    &= \int_0^{\infty} \mu \left( \left. W_p( \PP_k, \widehat\PP_k) > \eps_k-r_k\right| W_p(\PP_1,\PP_k) = r_k \right) {\rm d}F_k(r_k).
\end{align*}
If $\eps_k\geq r_k$, then the integrand of the last expression can be further simplified to
\begin{align*}
    & \mu \left( \left. W_p( \PP_k, \widehat\PP_k) > \eps_k-r_k\right| W_p(\PP_1,\PP_k) = r_k \right) \\
    & \quad = \mathbb E_\mu\left[ \left. \mu \left( \left. W_p( \PP_k, \widehat\PP_k) > \eps_k-r_k\right| \PP_1,\ldots,\PP_K\right)  \right| W_p(\PP_1,\PP_k) = r_k \right]  \leq \beta(\eps_k-r_k,N_k),
\end{align*}
where the first equality exploits the law of total expectation and the trivial observation that $W_p(\PP_1,\PP_k)$ reduces to a constant when we condition on~$\PP_{k'}$, $k'\in[K]$. The inequality follows from Lemma~\ref{lem:concentration-empirical-FG}, which applies because, conditional on any fixed realizations of the reference distributions, $\widehat \PP_k$ is governed by the distribution~$\PP_k^{N_k}$ under~$\mu$ and because the reference distributions satisfy Assumption~\ref{light_tailed_assumption} $\mu$-almost surely. Combining the above estimates shows that the prior satisfies
\[
    \mu \left(W_p( \PP_1, \widehat\PP_k)> \eps_k\right) \leq \int_0^{\eps_k} \beta(\eps_k-r_k,N_k) \, {\rm d}F_k(r_k) + 1-F_k(\eps_k).
\]
Substituting our estimates for the likelihood and the prior into~\eqref{eq:bayes-rule} finally yields
\begin{align*}
    & \mu\left(\left. W_p(\mathbb\PP_1, \widehat\PP_k)> \eps_k \right| W_p( \widehat\PP_1,\widehat\PP_k)\leq\widehat r_k \right) \leq 
    \frac{\beta(\eps_k- \widehat r_k,N_1) \cdot \left(\int_0^{\eps_k} \beta(\eps_k-r_k,N_k) \, {\rm d}F_k(r_k)+1-F_k(\eps_k) \right)}{ \mu \left(W_p( \widehat\PP_1,\widehat\PP_k)\leq\widehat r_k\right) }.
\end{align*}
The last expression coincides with~$\beta_{F_k}(\eps_k,\widehat r_k,N_1,N_k)$. We may therefore conclude that
\begin{align*}
    \mu\left(\left.\PP_1 \in \mathbb{B}^p_{\eps_k}(\widehat \PP_k) \right| W_p( \widehat\PP_1,\widehat\PP_k)\leq\widehat r_k \right) & =1- \mu\left(\left. W_p(\mathbb\PP_1, \widehat\PP_k)> \eps_k \right| W_p( \widehat\PP_1,\widehat\PP_k)\leq\widehat r_k\right) \\ & \geq 1-\beta_{F_k}(\eps_k,\widehat r_k,N_1,N_k),
\end{align*}
and thus the claim follows.
\QEDA
\endproof

\section{Relationships between Scenarios in Section~\ref{sec:stat_guarantees}}
\label{sec:scenario-relationships}
Note first that scenario~1 can be obtained as a special case of scenario~2 if the number of training samples from each source tends to infinity. Next, it is easy to see that scenario~2 emerges as a special case of scenario~5 under Dirac priors, that is, if~$F_k(\eps_k)=0$ for~$\eps_k<r_k$ and~$F_k(\eps_k)=1$ for~$\eps_k\geq r_k$, $k\in[K]$. Similarly, scenario~4 is obtained from scenario~5 by a Bayesian update after observing samples from~$\PP_1$. Finally, scenario~3 emerges as a special case of scenario~4 as the priors become less informative. The information available to the decision-maker in the different scenarios is summarized in Table~\ref{tab:scenarios}.

\begin{table}[]
    \centering\scriptsize{
    \begin{tabular}{c||c|c|c}
    \toprule
         Scenario & Access to $\mathbb{P}$ & Access to $\mathbb{P}_k$  & Information on distribution shift \\
         \hline\hline
         1 & \xmark &\cmark &$W_p(\mathbb{P},\mathbb{P}_k)\leq r_k$\\
         2 & \xmark &samples from $\mathbb{P}_k$ &$W_p(\mathbb{P},\mathbb{P}_k)\leq r_k$\\
         3 &samples from $\mathbb{P}$ &samples from $\mathbb{P}_k$ &\xmark\\
         4 &samples from $\mathbb{P}$ &samples from $\mathbb{P}_k$ &$W_p(\mathbb{P},\mathbb{P}_k)\sim F_k$\\
         5&\xmark&samples from $\mathbb{P}_k$&$W_p(\mathbb{P},\mathbb{P}_k)\sim F_k$\\
         \bottomrule
    \end{tabular}}
    \caption{Comparison of the different scenarios}
    \label{tab:scenarios}
\end{table}

\section{Additional Experiments}
\label{sec:additional-experiments}
This appendix presents additional numerical experiments.  Section~\ref{sec:portfolio-optimization_p2WDRO} repeats the experiment of Section~\ref{sec:synthetic-data} using $2$-Wasserstein distances. In addition,  Section~\ref{sec:assortment-optimization} evaluates the benefits of multi-source DRO in assortment optimization, while Section~\ref{sec:computation-time} compares the computational effort required to solve multi-source DRO problems with that of single-source DRO problems based on a Wasserstein barycenter.

\subsection{Additional Backtest on Synthetic Data}
\label{sec:portfolio-optimization_p2WDRO}
In all numerical experiments in Section~\ref{sec:experiments}, distributional proximity is measured using the $1$-Wasserstein distance. Accordingly, ambiguity sets and distributional barycenters are defined with respect to this distance, and the radii of the individual Wasserstein balls in multi-source DRO are calibrated as affine functions of the $1$-Wasserstein distance between the empirical source and target distributions. We now repeat the experiment from Section~\ref{sec:synthetic-data} using the $2$-Wasserstein distance instead. Thus, ambiguity sets and barycenters are defined with respect to the $2$-Wasserstein distance, while the radii of the individual Wasserstein balls in multi-source DRO are calibrated as affine functions of the corresponding $2$-Wasserstein distance between the empirical source and target distributions. We use $N_1=15$ samples from the target distribution and $N_2=60$ samples from the source distribution and set the search grid for the Wasserstein radii to $M\triangleq \{0.001, 0.005, 0.01, 0.05, 0.1\}$, while keeping all other model parameters identical to those in Section~\ref{sec:synthetic-data}.

Figure~\ref{fig:portfolio-results:synthetic_performance-W2} reports the results of this experiment, which are qualitatively similar to those of the original experiment in Section~\ref{sec:synthetic-data} (see Figure~\ref{fig:portfolio-results:synthetic_performance}). In particular, the multi-source DRO portfolios continue to achieve the highest expected returns and Sharpe ratios for most values of the cutoff threshold~$C$.

\begin{figure}[htbp]
     \centering
     \begin{subfigure}[b]{0.3\textwidth}
         \centering
         \includegraphics[width=\textwidth]{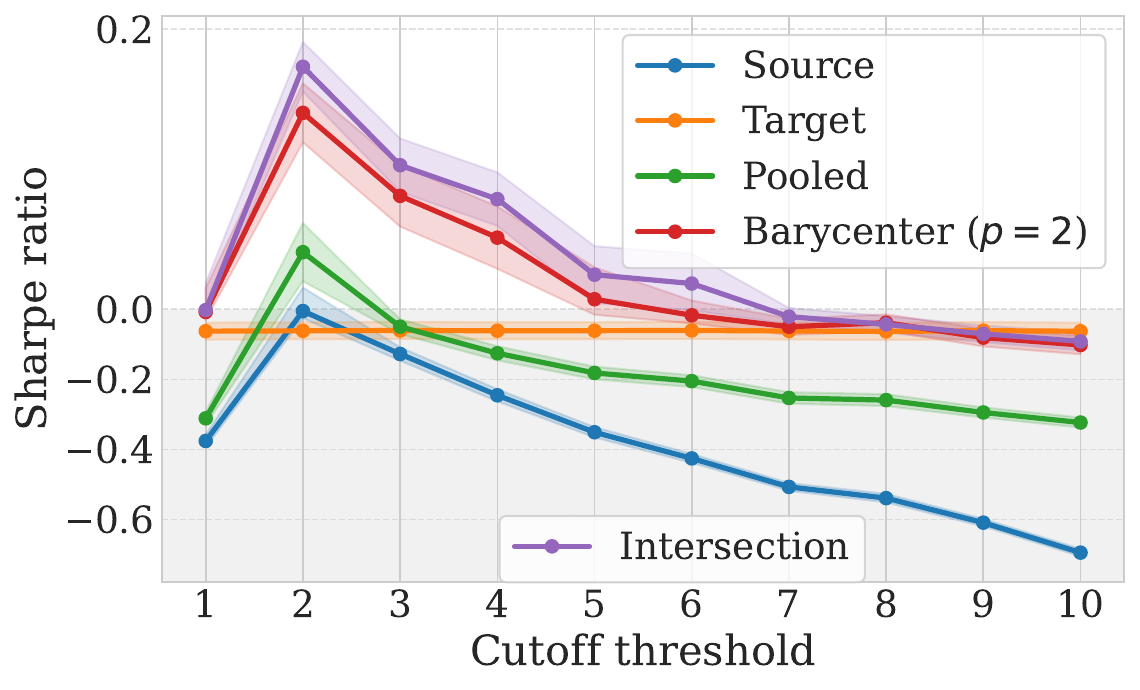}
         \caption{Sharpe Ratio}
     \end{subfigure}
     \begin{subfigure}[b]{0.3\textwidth}
         \centering
         \includegraphics[width=\textwidth]{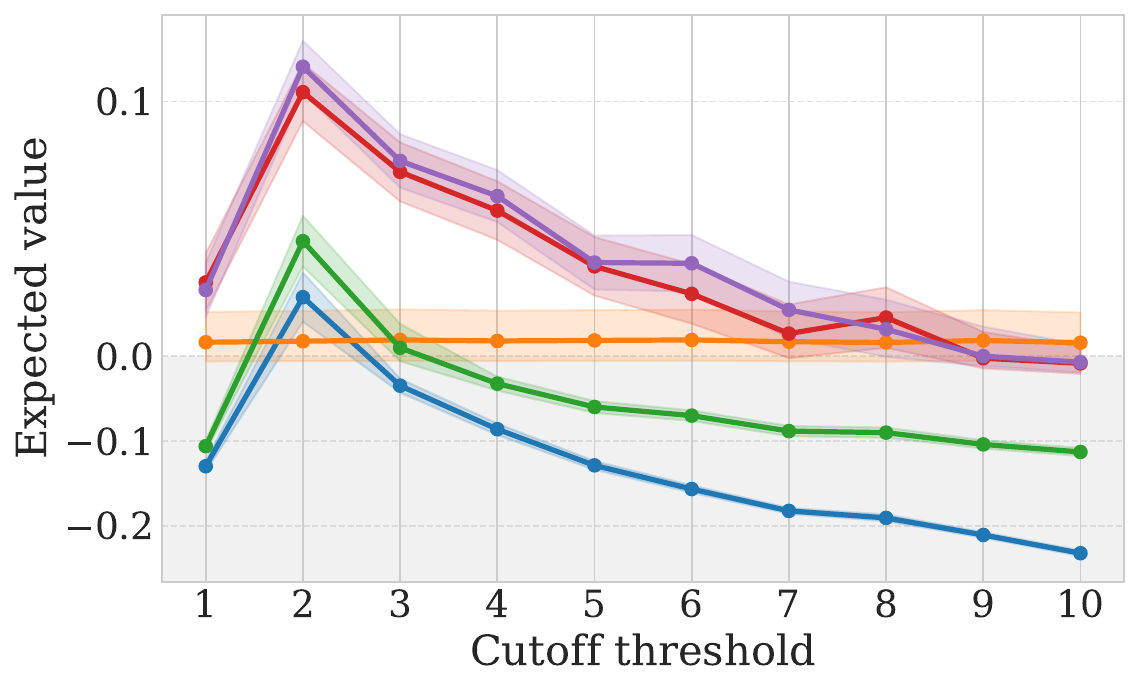}
         \caption{Expected Return}
     \end{subfigure}
     \begin{subfigure}[b]{0.3\textwidth}
         \centering
         \includegraphics[width=\textwidth]{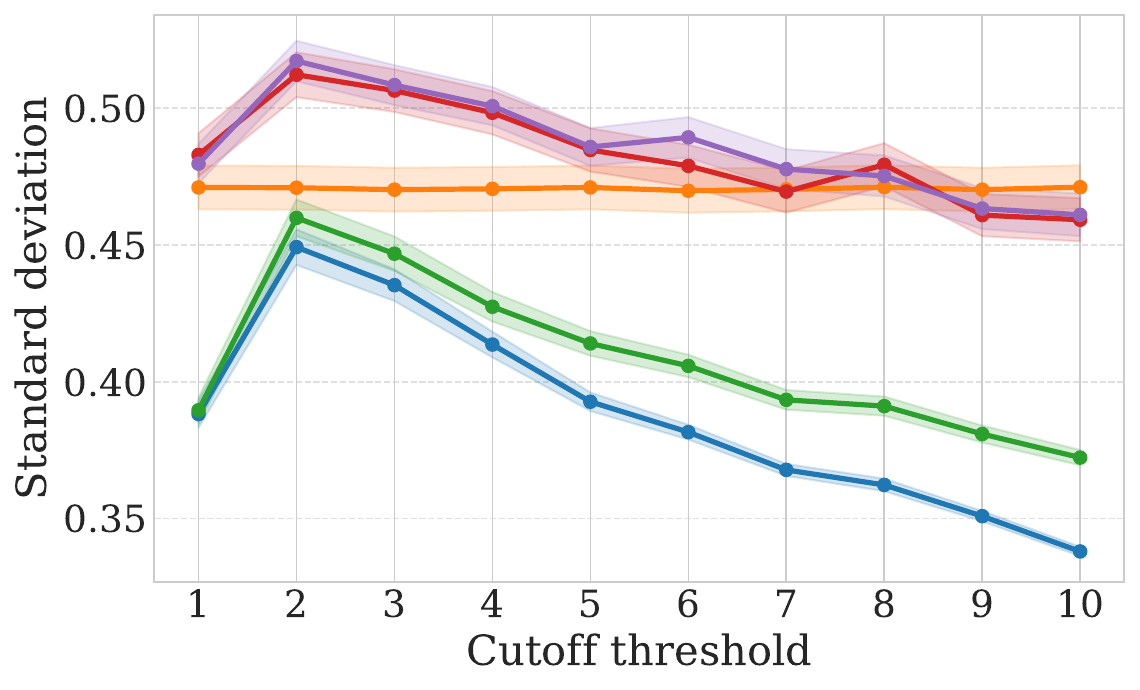}
         \caption{Return Standard Deviation}
     \end{subfigure}
     \begin{subfigure}[b]{0.3\textwidth}
    \centering
    \end{subfigure}
    \caption{Out-of-sample performance of optimal portfolios on synthetic data (mean (std.\ error) over 200~replications). For readability, the negative range is shown on a compressed linear scale in panels (a) and (b).}
    \label{fig:portfolio-results:synthetic_performance-W2}
\end{figure}

\subsection{Assortment Optimization}
\label{sec:assortment-optimization}
A firm operates two stores in different geographical regions, each of which sells~$d$ different products. For simplicity we assume that the products are perfect substitutes. The firm aims to open a third store in a new region. Due to space constraints, only a subset of~$B\leq d$ products can be offered in the new store. If $p_i\geq 0$ denotes the (given) price and $\xi_i\geq 0$ the (uncertain) demand of the $i$-th product for every $i\in[d]$, then the firm faces the following stylized assortment optimization problem.
\begin{equation}
\label{eq:assortment_problem}
  \max_{\theta \in \{ 0,1 \}^d} \left\{ \mathbb{E}_{\mathbb{P}} \left[\sum_{i=1}^d  \theta_i p_i \xi_i \right] \;:\; \sum_{i=1}^d\theta_i \leq B \right\}
\end{equation}
This problem seeks a subset of at most~$B$ products, encoded by the binary vector~$\theta$, that maximizes expected revenue. The expectation is evaluated with respect to the joint demand distribution~$\PP$ in the region of the new store. We henceforth assume that~$\PP$ is unknown but close to the $K=2$ demand distributions~$\PP_1$ and~$\PP_2$ in the regions of the existing stores. While~$\PP_1$ and~$\PP_2$ are unknown, too, the firm has access to respective empirical distributions~$\widehat\PP_1$ and~$\widehat\PP_2$ and can thus solve an instance of~\eqref{eq:original-problem} corresponding to~\eqref{eq:assortment_problem}. The ambiguity set is defined as the intersection of two 1-Wasserstein balls induced by the transportation cost function $c(\xi,\xi')=\|\xi-\xi'\|_1$. By Corollary~\ref{cor:pw-affine}, this instance of~\eqref{eq:original-problem} can be reformulated as a mixed-integer linear program. 

In our experiments we set~$d=10$ and~$B=3$, and we define $p_i=0.01\times i$ for all $i\in[d]$. The training, validation and test datasets are constructed from historical product demand data from Kaggle (\url{https://www.kaggle.com/datasets/felixzhao/productdemandforecasting}). This dataset contains order quantities for multiple products in three different geographical regions designated by~J, S and~C, respectively. Among the products that are available in all three regions, we select the 10 most popular ones in terms of average demand. Kernel density estimates of the marginal demand distributions of all products are shown in Figure~\ref{fig:product-demands}. We assume that the new store is located in region~J and that the two existing stores are located in regions~S and~C. This choice is motivated by the observation that the demand distributions in region~J are sometimes more similar to those in region~S (see products~4, 6, 7, 9 and~10 in Figure~\ref{fig:product-demands}) and sometimes more similar to those in region~C (see products~1, 2, 3, 5 and~8 in Figure~\ref{fig:product-demands}). Hence, none of the two source distributions provides a consistently better approximation for~$\PP$, and thus it makes sense to use data from both sources.

\begin{figure}{}
\centering
\includegraphics[width=\textwidth]{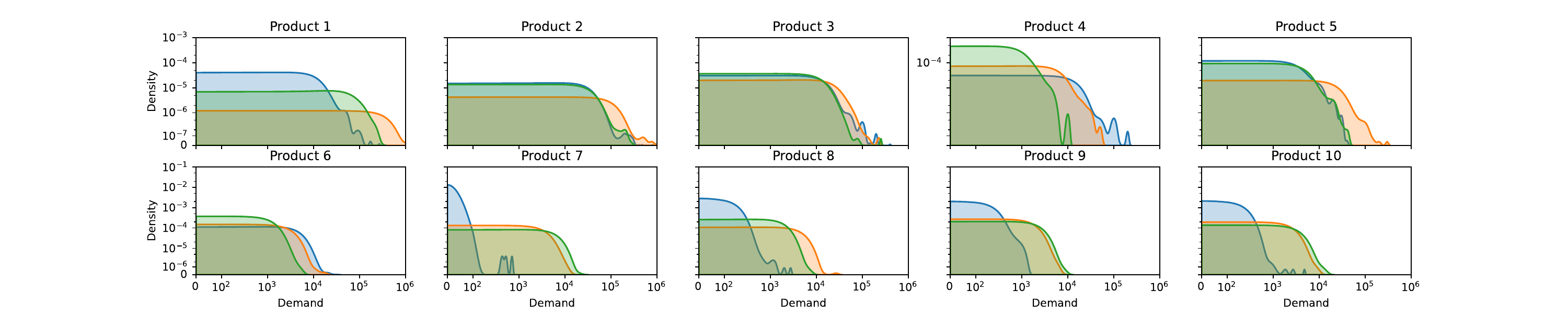}
\caption{Kernel density estimates of the marginal demand distributions in region J (blue), S ({orange}) and C (green).}
\label{fig:product-demands}
\end{figure}
As in Section~\ref{sec:synthetic-data}, we compare the optimal solutions of the multi-source DRO problem against those of different single-source DRO baselines. We use the same parametrization for the radii of the Wasserstein balls and the same search grids for all hyperparameters ($\lambda$ and~$m$ in multi-source DRO and~$\eps$ in single-source DRO) as in Section~\ref{sec:synthetic-data}. In all cases, we select hyperparameters from the underlying search grids that maximize the average revenue across ten independent validation samples. The validation samples are drawn from the same distribution as the training samples. For methods that use data from~$\PP_1$ as well as~$\PP_2$, we select five validation samples from each distribution.

We compute the expected revenues of all optimal assortments on a dataset of 3{,}372 test samples.  Table~\ref{tab:msci-assortment-results} reports the means and standard errors of these expected revenues across 50~independent simulation runs for different choices of the source sample sizes~$N_1$ and~$N_2$. We observe that multi-source DRO consistently outperforms all baselines. The dominance of multi-source DRO is linked to the structure of the two source distributions, none of which provides a globally better approximation for the target distribution. Thus, using both data sources is better than using only one. In addition, multi-source DRO makes better use of the two data sources than traditional DRO with a Wasserstein ball centered at a mixture or at a Wasserstein barycenter of the empirical source distributions.

\begin{table}\caption{Out-of-sample expected revenue of optimal product portfolios (mean (std.\ error) over 50~replications)}
    \centering
    {\scriptsize
\begin{tabular}{rr|lllllll}
\toprule
$N_1$& $N_2$ &        Single-s.\ DRO on &        Single-s.\ DRO on&        Single-s.\ DRO on&    Single-s.\ DRO on&          Multi-s.\ DRO \\
   & & distribution $\widehat\PP_1$ & distribution $\widehat\PP_2$ & pooled data &  barycenter &      &     \\
\midrule
\multirow{4}{*}{25} &  25 & 485.82\ (41.38) & 424.02\ (31.78) & 466.95\ (31.39) & 465.45\ (38.88) & $\boldsamewidth{507.25\ (18.35)}$ \\
  & 50 & $\boldsamewidth{499.72\ (39.18)}$ & 436.07\ (36.59) & 482.39\ (28.48) & 483.29\ (41.97) & 489.50\ (28.62) \\
  & 75 & 490.15\ (37.08) & 460.22\ (26.53) & 474.18\ (31.37) & 483.54\ (34.28) & $\boldsamewidth{501.91\ (18.86)}$ \\
  & 100 & 511.38\ (33.67) & 429.77\ (41.68) & 476.26\ (25.17) & 482.80\ (39.65) & $\boldsamewidth{521.38\ (17.45)}$ \\
\cline{1-8}
  \multirow[t]{4}{*}{50} & 25 & 482.97\ (33.25) & 469.61\ (31.32) & 495.78\ (15.07) & 474.53\ (31.33) & $\boldsamewidth{515.70\ (16.19)}$ \\
   & 50 & 483.40\ (27.96) & 450.81\ (31.01) & 486.00\ (15.10) & 474.99\ (29.04) & $\boldsamewidth{506.45\ (14.28)}$ \\
  & 75 & 459.22\ (34.86) & 448.90\ (29.64) & 471.02\ (21.68) & 457.32\ (33.39) & $\boldsamewidth{508.48\ (21.13)}$ \\
   & 100 & 480.09\ (31.86) & 464.28\ (30.30) & 480.90\ (20.34) & 483.12\ (25.40) & $\boldsamewidth{514.36\ (16.51)}$ \\
\cline{1-8}
  \multirow[t]{4}{*}{75} & 25 & 472.61\ (27.54) & 455.59\ (29.32) & 481.00\ (17.28) & 466.65\ (30.37) & $\boldsamewidth{510.46\ (13.53)}$ \\
   & 50 & 468.62\ (28.61) & 456.25\ (30.65) & 476.89\ (19.65) & 462.41\ (30.28) & $\boldsamewidth{513.03\ (16.80)}$ \\
  & 75 & 478.82\ (27.72) & 447.56\ (29.20) & 497.11\ (15.06) & 462.17\ (31.01) & $\boldsamewidth{511.79\ (13.22)}$ \\
  & 100 & 470.48\ (24.51) & 445.45\ (31.54) & 479.66\ (17.03) & 451.81\ (27.77) & $\boldsamewidth{506.35\ (17.90)}$ \\
\cline{1-8}
  \multirow[t]{4}{*}{100} & 25 & 492.03\ (22.29) & 472.98\ (28.43) & 498.45\ (15.02) & 491.31\ (19.52) & $\boldsamewidth{520.69\ (10.99) }$\\
  & 50 & 480.81\ (24.30) & 471.65\ (28.21) & 483.77\ (14.06) & 482.63\ (20.46) & $\boldsamewidth{505.78\ (16.60)}$ \\
  & 75 & 472.79\ (21.61) & 473.64\ (29.72) & 486.44\ (14.47) & 479.74\ (17.67) & $\boldsamewidth{514.36\ (16.51)}$ \\
  & 100 & 472.15\ (26.71) & 458.26\ (30.04) & 493.11\ (15.02) & 474.90\ (27.56) & $\boldsamewidth{517.13\ (11.58)}$ \\
\bottomrule
\end{tabular}
}
    \label{tab:msci-assortment-results}
\end{table}

\subsection{Runtime Comparison}
\label{sec:computation-time}

In the last experiment we investigate an instance of the stochastic optimization problem~\eqref{eq:sp} with~$\Theta=\R^d$, $\Xi=\R^d$ and $\ell(\theta,\xi)=\|\xi-\theta\|_1$, which seeks the componentwise median of the target distribution~$\PP$. Note that the loss function can be expressed as a pointwise maximum of~$2^d$ many affine functions, that is, 
$$
    \ell(\theta,\xi)=\max_{m\in\{-1,1\}^d} m^\top(\xi-\theta).
$$
Loss functions of this type arise, for instance, in multivariate newsvendor problems with symmetric cost. 
We assume as usual that the target distribution is unknown but that~$N$ samples from each of~$K$ source distributions are available. We draw these samples independently from the respective source distributions~$\mathbb{P}_k=\mathcal{N}(\mu_k,I_d)$, where the mean vectors~$\mu_k$ are sampled from the uniform distribution on~$[-3,3]^d$ for all~$k\in[K]$.

Based on the given information we can construct an instance of the multi-source DRO problem~\eqref{eq:original-problem}, whose ambiguity set is given by the intersection of 2-Wasserstein balls around the $K$ empirical source distributions. This means that we set the transportation cost function to~$c(\xi_1,\xi_2)=\|\xi_1-\xi_2\|_2^2$. To ensure that the intersection is nonempty, we first compute $\widehat r_k=W_2(\widehat{\mathbb{P}}_1,\widehat{\mathbb{P}}_k)$ for all~$k=2,\ldots,K$. Then, we set the radii of the 2-Wasserstein balls to~$r_1=\frac{1}{K-1}\sum_{k=2}^{K}\widehat r_k$ and to~$r_k=1.3 \widehat r_k$ for all~$k=2,\ldots,K$. By Corollary~\ref{cor:pw-affine}, the inner maximization problem in~\eqref{eq:original-problem} can be reformulated as a minimization problem, which can be merged with the outer minimization problem over~$\theta$. The resulting monolithic optimization problem constitutes a second-order cone program with $\mathcal O(N^K2^d)$ decision variables and $\mathcal O(N^K2^d)$ constraints. 



The given information can also be used to construct an instance of the single-source DRO problem~\eqref{eq:dro}, whose ambiguity set is given by a single 2-Wasserstein ball of radius $r = \max_{k\in[K]} W_2(\widehat\PP_k, \widehat \PP)$ around a 2-Wasserstein barycenter~$\widehat\PP$ of the $K$ empirical source distributions. To construct~$\widehat\PP$, we solve problem~\eqref{eq:barycenter} with $\lambda_k=1/K$, $\PP_k=\widehat\PP_k$ and $C({ \QQ},\PP_k)=W_2^2({ \QQ},\widehat \PP_k)$ for every~$k\in[K]$. This can be done as explained in Theorem~\ref{thm:barycenter:pushforward-reformulation} by solving the multi-margin OT problem~\eqref{eq:barycenter:outer-min}, which constitutes a linear program in $N^K$ decision variables and $NK$ constraints. Note that as~$c(\xi_1,\xi_2)=\|\xi_1-\xi_2\|_2^2$, the functions $\Phi(\xi_1,\ldots, \xi_K)$ and $\phi(\xi_1,\ldots, \xi_K)$ evaluate to the mean vector and the trace of the covariance matrix of the points $\xi_1,\ldots,\xi_K$, respectively. Next, we can use Corollary~\ref{cor:pw-affine} with~$K=1$ to reformulate the single-source DRO problem~\eqref{eq:dro} as a monolithic second-order cone program. As $\widehat\PP$ has at most $(N-1)K+1$ atoms \cite[Lemma~2.2]{altschuler2021barycenterComplexity}, this optimization problem accommodates $\mathcal O(NK2^d)$ decision variables and $\mathcal O(NK2^d)$ constraints.

Figure~\ref{fig:run_times_barycenterDRO} reports the solution times of the second-order cone reformulations of the multi-source and single-source DRO models for different values of~$d$, $K$, and~$N$. All optimization problems are solved with MOSEK~11.1.11 in Python~3.14.3 on a Linux server using up to~14 threads and 770\,GB of RAM. The relative optimality gap is set to~$0.1\%$, with presolve enabled and feasibility and complementarity tolerances set to~$10^{-3}$. All other parameters remain at their default values. Problems with~$K \leq 3$ can be solved within minutes, but runtimes increase rapidly with~$K$, highlighting the need for improved scalability in multi-source DRO. Possible approaches include dimensionality reduction via principal component analysis~\citep{cheng2018DRO-PCA} and sample reduction through clustering~\citep{wang2025mean}. Moreover, developing scalable stochastic gradient methods for large-scale, potentially non-convex multi-source DRO problems with large~$K$ remains an important direction, building on single-source DRO techniques such as those in~\citep{sinha18certifying}.

\begin{figure}{}
\centering
\includegraphics[width=\textwidth]{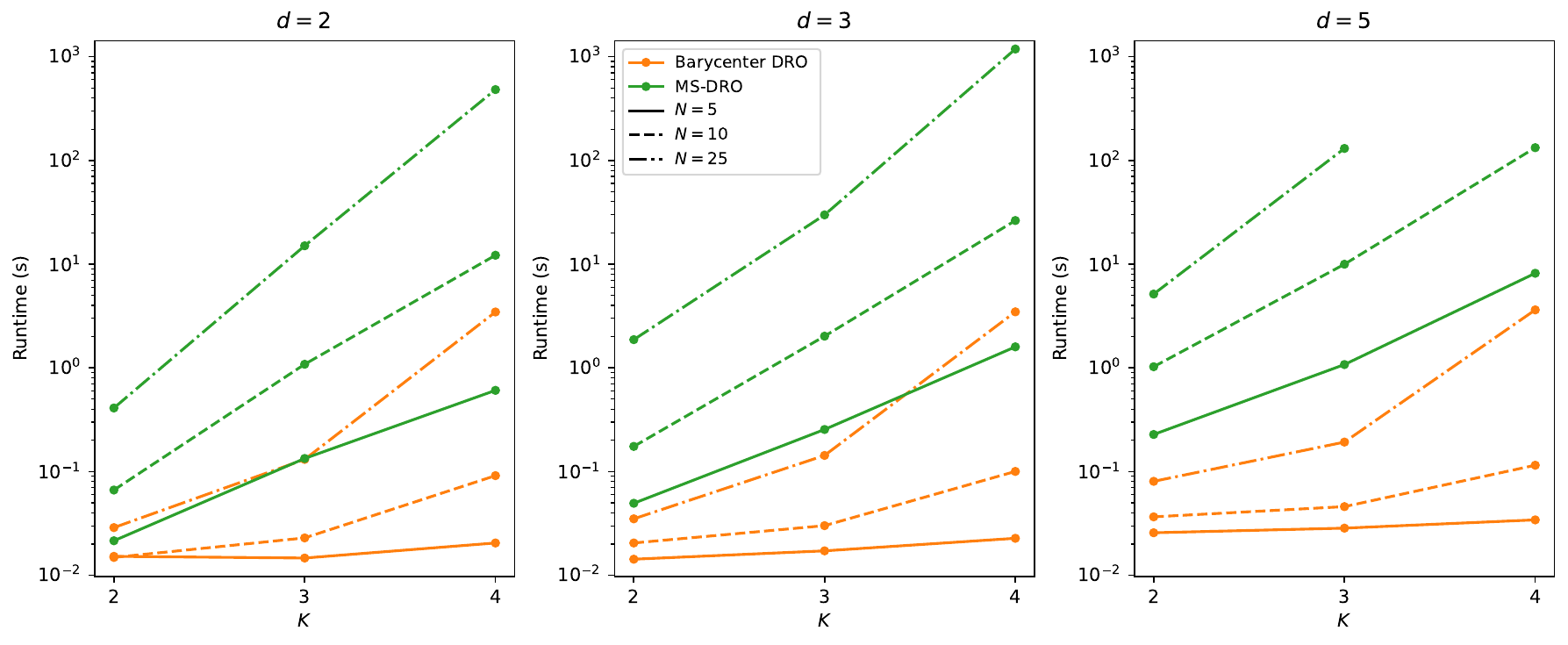}
\caption{Runtimes of the second-order cone reformulations of the multi-source DRO (MS-DRO) and barycenter-based single-source DRO (Barycenter DRO) models across different configurations of~$d$, $K$, and~$N$.}
\label{fig:run_times_barycenterDRO}
\end{figure}

\end{APPENDICES}

\end{document}